%% file: AarhusII.tex
\def\draft{n}

\documentclass[12pt]{amsart}
\usepackage{amssymb,epic,eepic,amscd,xr}
\input macros.tex
\input defs.tex

\begin{document}

\title[The \AA{}rhus integral II: Invariance and Universality]
  {
  The \AA{}rhus integral of rational homology 3-spheres II: Invariance and
  universality
}

\author[Bar-Natan]{Dror~Bar-Natan}
\address{Institute of Mathematics\\
        The Hebrew University\\
        Giv'at-Ram, Jerusalem 91904\\
        Israel}
\email{drorbn@math.huji.ac.il}

\author[Garoufalidis]{Stavros~Garoufalidis}
\address{Department of Mathematics\\
        Brandeis University\\
        Waltham MA 02254-9110\\
        USA}
\curraddr{Department of Mathematics\\
  Harvard University\\
  Cambridge MA 02138\\
  USA
}
\email{stavros@math.harvard.edu}

\author[Rozansky]{Lev~Rozansky}
\address{Department of Mathematics, Statistics, and Computer Science\\
        University of Illinois at Chicago\\
        Chicago IL 60607-7045\\
        USA}
\curraddr{Department of Mathematics\\
  Yale University\\
  10 Hillhouse Avenue\\
  P.O. Box 208283 \\
  New Haven, CT 06520-8283\\
  USA
}
\email{rozansky@math.yale.edu}

\author[Thurston]{Dylan~P.~Thurston}
\address{Department of Mathematics\\
        University of California at Berkeley\\
        Berkeley CA 94720-3840\\
        USA}
\email{dpt@math.berkeley.edu}

\thanks{This paper is available electronically at
  {\tt http://www.ma.huji.ac.il/\~{}drorbn}, at \newline
  {\tt http://jacobi.math.brown.edu/\~{}stavrosg}, and at {\tt
    http://xxx.lanl.gov/abs/math/9801049}.
}

\dedicatory{To appear in {\em Selecta Mathematica.}}
\date{This edition: Feb.~7,~1999; \ \ First edition: January~11, 1998.}

\begin{abstract}
We continue the work started in~\cite{Aarhus:I}, and prove the
invariance and universality in the class of finite type invariants of
the object defined and motivated there, namely the \Arhus{} integral of
rational homology 3-spheres. Our main tool in proving invariance is a
translation scheme that translates statements in multi-variable
calculus (Gaussian integration, integration by parts, etc.) to
statements about diagrams. Using this scheme the straight-forward
``philosophical'' calculus-level proofs of~\cite{Aarhus:I} become
straight-forward honest diagram-level proofs here. The universality
proof is standard and utilizes a simple ``locality'' property of the
Kontsevich integral.
\end{abstract}

\maketitle

\tableofcontents

\input{intro}
\input{calculus}
\input{Invariance}
\input{Universality}
\input{odds}

\par\noindent
{\bf Acknowledgement: } The seeds leading to this work were planted when
the four of us (as well as Le, Murakami (H\&J), Ohtsuki, and many other
like-minded people) were visiting \Arhus, Denmark, for a special semester
on geometry and physics, in August 1995. We wish to thank the organizers,
J.~Dupont, H.~Pedersen, A.~Swann and especially J.~Andersen for their
hospitality and for the stimulating atmosphere they created.  We also wish
to thank N.~Habegger, A.~Haviv, M.~Hutchings, T.~Q.~T.~Le, A.~Referee,
N.~Reshetikhin and S.~Yitzchaik for additional remarks and suggestions,
the Center for Discrete Mathematics and Theoretical Computer Science at
the Hebrew University for financial support, and the Volkswagen-Stiftung
(RiP-program in Oberwolfach) for their hospitality and financial support.

\input{refs}

\end{document}

%% file: macros.tex
\theoremstyle{plain}

\newtheorem{theorem}{Theorem}
\newtheorem{proposition}{Proposition}[section]

\newtheorem{lemma}[proposition]{Lemma}
\newtheorem{corollary}[proposition]{Corollary}

\theoremstyle{definition}
\newtheorem{definition}[proposition]{Definition}

\newtheorem{problem}[proposition]{Problem}

\theoremstyle{remark}

\newtheorem{exercise}[proposition]{Exercise}
\newtheorem{hint}[proposition]{Hint}
\newtheorem{remark}[proposition]{Remark}

\def\printname#1{
	\if\draft y
		\smash{\makebox[0pt]{\hspace{-0.5in}
			\raisebox{8pt}{\tt\tiny #1}}}
	\fi
}

\newcommand{\mathmode}[1]{$#1$}

\newlength{\standardunitlength}
\catcode`\@=11
\long\def\@makecaption#1#2{%
    \vskip 10pt
    \setbox\@tempboxa\hbox{
      \small\sf{\bfcaptionfont #1. }\ignorespaces #2}%
    \ifdim \wd\@tempboxa >\captionwidth {%
        \rightskip=\@captionmargin\leftskip=\@captionmargin
        \unhbox\@tempboxa\par}%
      \else
        \hbox to\hsize{\hfil\box\@tempboxa\hfil}%
    \fi}
\font\bfcaptionfont=cmssbx10 scaled \magstephalf
\newdimen\@captionmargin\@captionmargin=2\parindent
\newdimen\captionwidth\captionwidth=\hsize
\catcode`\@=12

\def\lbl#1{\label{#1}\printname{#1}}

\def\qed{{\hfill\text{$\Box$}}}

\newlength{\globalparindent}
\newcommand{\udot}{{\mathaccent\cdot\cup}}

%% file: defs.tex
\advance \topmargin by -\headheight
\advance \topmargin by -\headsep
\evensidemargin \oddsidemargin
\newcommand{\Aa}{\text{\sl\AA}}
\newcommand{\Arhus}{\AA rhus}

\newcommand{\BCH}{\text{BCH}}

\newcommand{\Gint}{\int^{FG}}

\newcommand{\QQ}{\mathbb Q}
\newcommand{\RR}{\mathbb R}

\newcommand{\Zcheck}{\check{Z}}
\newcommand{\calA}{\mathcal A}

\newcommand{\calB}{\mathcal B}

\newcommand{\frakg}{\mathfrak g}

\renewcommand{\qed}{~\hfill$\square$}

\newcommand\strutf[2]{\,^{#1}\,^{\hbox to2em{\hrulefill}}\,_{#2}}
\newcommand\strutn[2]{\,^{#1}\!\!\frown^{#2}}
\newcommand\strutu[2]{\,_{#1}\!\!\smile_{#2}}

\newcommand\mright{\overrightarrow{m}}
\newcommand\Upsright{\overrightarrow{\Upsilon}}

\newcommand\apply{\mathbin\flat}

%% file: intro.tex
\section{Introduction} This paper is the second in a four-part series on
``the \Arhus{} integral of rational homology 3-spheres''. In the first
part of this series, \cite{Aarhus:I}, we gave the definition of a
diagram-valued invariant $\Aa$ of ``regular pure tangles'', pure
tangles whose linking matrix is non-singular,\footnote{
  A precise definition of regular pure tangles appears
  in~\cite[definition~\ref{I-def:RPT}]{Aarhus:I}. It is a good idea to
  have~\cite{Aarhus:I} handy while reading this paper, as many of the
  definitions introduced and explained there will only be repeated here
  in a very brief manner.
}
and gave ``philosophical'' reasons why $\Aa$ should descend to an
invariant of regular links (framed links with non-singular linking
matrix), and as such satisfy the Kirby relations and hence descend further
to an invariant of rational homology 3-spheres. Very briefly, we have
defined the pre-normalized \Arhus{} integral $\Aa_0$ to be the composition

\[ \begin{CD}
  \Aa_0:\left\{\parbox{0.5in}{\begin{center}
    regular pure tangles
  \end{center}}\right\}=RPT
  @>\Zcheck>{\parbox{1.1in}{\scriptsize\begin{center}
    the \cite{LeMurakami2Ohtsuki:3Manifold} version of the Kontsevich integral
  \end{center}}}>
  \calA(\uparrow_X)
  @>\sigma>{\parbox{0.4in}{\scriptsize\begin{center}
    formal PBW
  \end{center}}}>
  \calB(X)
  @>\Gint>{\parbox{0.65in}{\scriptsize\begin{center}
    formal Gaussian integration
  \end{center}}}>
  \calA(\emptyset).
\end{CD} \]

In this formula,
\begin{itemize}
\item $RPT$ denotes the set of regular pure tangles whose components are
  marked by the elements of some finite set $X$
  (see~~\cite[definition~\ref{I-def:RPT}]{Aarhus:I}).
\item $\Zcheck$ denotes the Kontsevich integral normalized as
  in~\cite{LeMurakami2Ohtsuki:3Manifold}
  (check~\cite[definition~\ref{I-def:Zcheck}]{Aarhus:I} for the adaptation to
  pure tangles).
\item $\calA(\uparrow_X)$ denotes the completed graded space of chord
  diagrams for $X$-marked pure tangles modulo the usual $4T/STU$
  relations (see~\cite[definition~\ref{I-def:TangleDiagrams}]{Aarhus:I}).
\item $\sigma$ denotes the diagrammatic version of the Poincare-Birkhoff-Witt
  theorem (defined as in~\cite{Bar-Natan:Vassiliev, Bar-Natan:Homotopy},
  but normalized slightly differently, as
  in~\cite[definition~\ref{I-def:sigma}]{Aarhus:I}).
\item $\calB(X)$ denotes the completed graded space of $X$-marked
  uni-trivalent diagrams as in~\cite{Bar-Natan:Vassiliev, Bar-Natan:Homotopy}
  and~\cite[definition~\ref{I-def:UnitrivalentDiagrams}]{Aarhus:I}.
\item $\calA(\emptyset)$ denotes the completed graded space of manifold
  diagrams as in~\cite[definition~\ref{I-def:ManifoldDiagrams}]{Aarhus:I}.
\item $\Gint$ is a new ingredient, first introduced
  in~\cite[definition~\ref{I-def:Gint}]{Aarhus:I}, called ``formal
  Gaussian integration''. In a sense explained there and developed
  further here, it is a diagrammatic analogue of the usual notion of
  Gaussian integration.
\end{itemize}

Our main challenge in this paper is to prove that $\Aa_0$ descends to
an invariant of links which is invariant under the second Kirby move.
As it turns out, this depends heavily on understanding properties of
formal Gaussian integration, which are all analogues of properties of
standard integration over Euclidean spaces.  We develop the necessary
machinery in section~\ref{sec:calculus} of this article, and then in
section~\ref{sec:Invariance} we move on and use this machinery to prove
two of our main results, proposition~\ref{prop:KirbyII} and
theorem~\ref{thm:Invariance}:

\begin{proposition} \lbl{prop:KirbyII}
The regular pure tangle invariant $\Aa_0$ descends to an invariant
of regular links and as such it is insensitive to orientation flips (of
link components) and invariant under the second Kirby move.
\end{proposition}

\begin{definition} Let $U_\pm$ be the unknot with framing $\pm 1$, and let
$\sigma_+$ ($\sigma_-$) be the number of positive (negative) eigenvalues of
the linking matrix of a regular link $L$. Let the \Arhus{} integral
$\Aa(L)$ of $L$ be
\begin{equation} \lbl{eq:Renormalization}
  \Aa(L)=\Aa_0(U_+)^{-\sigma_+}\Aa_0(U_-)^{-\sigma_-}\Aa_0(L),
\end{equation}
with all products and powers taken using the disjoint union product of
$\calA(\emptyset)$.
\end{definition}

\begin{theorem} \lbl{thm:Invariance}
$\Aa$ is invariant under orientation flips and under both Kirby moves,
and hence (\cite{Kirby:Calculus}) it is an invariant of rational homology
3-spheres.
\end{theorem}

Our second goal in this article is to prove that $\Aa$ is a
universal Ohtsuki invariant, and hence that all $\QQ$-valued
finite-type invariants of integer homology spheres are compositions of
$\Aa$ with linear functionals on $\calA(\emptyset)$. We present all
relevant definitions and proofs in section~\ref{sec:Universality}
below.

%% file: calculus.tex
\section{Formal diagrammatic calculus} \lbl{sec:calculus}

In this section we study the theory of formal Gaussian integration,
along with a neighboring theory of formal differential operators. The
idea is that monomials can be represented by vertices of certain
valencies, and differentiation (almost always) and integration
(at least in the case of Gaussian integration) are given by
combinatorial formulas that can be viewed as manipulations done on
certain kinds of diagrams built out of these vertices. This extracts
some parts of good old elementary calculus, and replaces algebraic
manipulations by a diagrammatic calculus. Now forget the interpretation
of diagrams as functions and operators, and you will be left with a
formal theory of diagrams in which there are formal diagrammatic
analogs of various calculus operations and of certain theorems from
classical calculus.

This diagrammatic theory is more general than what we need for this
paper; it is not restricted to the diagrams (and relations) that make
up the spaces that we use often, such as $\calA$ and $\calB$.  We are
sure such a general formal diagrammatic theory was described many times
before and we make no claims of originality. This theory is implicit in
many discussions of Feynman diagrams in physics texts, but we are not
aware of a good reference that does everything that we need the way we
need it.  Hence in this section we describe in some detail that part of
the general theory that we will use in the later sections.

\subsection{The general setup.} \lbl{subsec:GeneralSetup}
Our basic objects are diagrams with some internal structure (that we
mostly do not care about), and some number of ``legs'', outward
pointing edges that end in a univalent vertex. The legs are labeled by
a vector space, or by a variable that lives in the
dual of that vector space. We think (for the purpose of the analogy
with standard calculus) of such a diagram as representing a tensor
in the tensor product of the spaces labeled next to its legs. We assume
that legs that are labeled the same way are interchangeable, meaning
that our diagrams represent symmetric tensors whenever labels are
repeating. Symmetric tensors can be identified with polynomials on the
dual:
\[
  \printname{DiagramPolynomial}
	\setlength{\unitlength}{0.5\standardunitlength}
	\begin{array}{c}  \hspace{-1.7mm}
        	\raisebox{-8pt}{\input draws/DiagramPolynomial.tex }
        	\hspace{-1.9mm}
	\end{array}

  \mapsto\quad
  \parbox{1.85in}{
     a polynomial of degree $3$ in $x$ and degree $2$ in $y$.
  }
\]

We also allow legs labeled by dual spaces or by dual variables. By
convention, the variable dual to $x$ is denoted $\partial_x$. Just as
the symmetric algebra $S(V^\star)$ can be regarded as a space of constant
coefficient differential operators acting on $S(V)$, diagrams labeled
by dual variables represent differential operators:
\begin{eqnarray*}
  \printname{DiagramOperator}
	\setlength{\unitlength}{0.5\standardunitlength}
	\begin{array}{c}  \hspace{-1.7mm}
        	\raisebox{-8pt}{\input draws/DiagramOperator.tex }
        	\hspace{-1.9mm}
	\end{array}

  & \mapsto & \qquad 
  \parbox{2.2in}{
     A second order differential operator acting on functions of $y$, with
     a coefficient cubic in $x$.
  } \\
  & & \\
  \printname{CCOperator}
	\setlength{\unitlength}{0.5\standardunitlength}
	\begin{array}{c}  \hspace{-1.7mm}
        	\raisebox{-8pt}{\input draws/CCOperator.tex }
        	\hspace{-1.9mm}
	\end{array}

  & \mapsto & \qquad
  \parbox{2.2in}{
    A fourth order constant coefficient differential operator.
  }
\end{eqnarray*}

We assume in addition that each diagram has an ``internal degree'', some
non-negative half integer ($0,\frac12,1,\frac32,\ldots$), associated with
it. It is to be thought of as the degree in some additional (small or
formal) parameter $\hbar$ that the whole theory depends upon. That is,
the polynomials and differential operators that we imitate also depend
on some additional parameter $\hbar$.

We also consider weighted sums of diagrams (representing not
necessarily homogeneous polynomials and differential operators), and
even infinite weighted sums of diagrams provided either their internal
degree grows to infinity or their number of legs grows to infinity.
These infinite sums represent power series (in $\hbar$ and/or in the
variables labeled on the legs) and/or infinite order differential
operators.

Finally, we allow some ``internal relations'' between the diagrams
involved. That is, sometimes we mod out the spaces of diagrams involved
by relations, such as the $IHX$ and $AS$ relation, that do not touch
the external legs and the internal degree of a diagram. All operations
that we will discuss below only involve the external legs and/or the
internal degree, and so they will be well-defined even after moding out
by such internal relations.

We then consider some operations on such diagrams. The operation of
adding diagrams (whose output is simply the formal sum of the summands)
corresponds to additions of polynomials or operators. The operation of
disjoint union of diagrams (adding their internal degrees), extended
bilinearly to sums of diagrams, corresponds to multiplying polynomials
and/or composing differential operators (at least in the constant
coefficients case, where one need not worry about the order of
composition). Once summation and multiplication are available, one can
define exponentiation and other analytic functions using power series
expansions.

The most interesting operation we consider is the operation of
contraction (or ``gluing'').  Two tensors, one in, say, $V^\star\otimes
W$ and the other in, say, $V\otimes Z$, can be contracted, and the
result is a new tensor in $W\otimes Z$. The graphical analog of this
operation is the fusion of two diagrams along a pair (or pairs) of legs
labeled by dual spaces or variables (while adding their internal
degrees). In the case of legs labeled by dual variables, the calculus
meaning of the fusion operation is the pairing of a derivative with a
linear function. The laws of calculus dictate that when a differential
operator $D$ acts on a monomial $f$, the result is the sum of all
possible ways of pairing the derivatives in $D$ with the factors of $f$.
Hence if $D$ is a diagram representing a differential operator (i.e.,
it has legs labeled $\partial_x$, $\partial_y$, etc.) and $f$ is a
diagram representing a function (legs labeled $x,y,\ldots$), we define
\[ D\apply f =
  \left(\parbox{2.95in}{
    sum of all ways of gluing all legs labeled $\partial_x$ on $D$
    with some or all legs labeled $x$ on $f$ (and same for $\partial_y$
    and $y$, etc.)
  }\right).
\]
(This sum may be $0$ if there are, say, more legs labeled $\partial_x$
on $D$ than legs labeled $x$ on $f$). For example,
\[ 
  \def\px{{\partial_x}}
  \printname{DfExample}
	\setlength{\unitlength}{0.5\standardunitlength}
	\begin{array}{c}  \hspace{-1.7mm}
        	\raisebox{-8pt}{\input draws/DfExample.tex }
        	\hspace{-1.9mm}
	\end{array}
.
\]
(In this figure 4-valent vertices are not real, but just artifacts of
the planar projection). If this were calculus and the spaces involved
were one-dimensional, we would call the above formula a proof that
$x(\partial_x)^2\,x^3y^2=6x^2y^2$.

\begin{exercise} Show that Leibnitz's formula, $D\apply(fg)=(D\apply
f)g+f(D\apply g)$ holds in our context, whenever $D$ is a first order
differential operator.
\end{exercise}

\begin{hint} Multiplication is disjoint union. $D$ can connect to the
disjoint union of $f$ and $g$ either by connecting to a leg of $f$, or by
connecting to a leg of $g$.
\end{hint}

\begin{exercise} \lbl{ex:ExponentialLiebnitz}
Prove the exponential Leibnitz's formula, $(\exp D)\apply\prod
f_i=\prod(\exp D)\apply f_i$, where $D$ is first order in $x$ and has
no coefficients proportional to $x$ (i.e., where $D$ has one leg labeled
$\partial_x$ and no legs labeled $x$).
\end{exercise}

Below we will need at some technical points an extension of this exercise
to the case when the operators involved are not necessarily first order.
The result we need is a bit difficult to formulate, and doing so precisely
would take us too far aside. But nevertheless, the result is rather easy
to understand in ``chemical'' terms, in which diagrams are replaced by
molecules and exponentiations are replaced by substance-filled containers.
Notice that the exponentiation of some object ${\mathcal O}$ is the
sum $\sum_k{\mathcal O}^k/k!$ of all ways of taking ``many'' unordered
copies of ${\mathcal O}$, so it can be thought of ``taking a big container
filled with (copies of) the molecule ${\mathcal O}$''.

A ``homogeneous reaction'' (in chemistry) is a reaction in which a
homogeneous mixture $A$ of mutually inert reactants is mixed with another
homogeneous mixture $B$ of mutually inert reactants, allowing reactions
to occur and products to be produced. The result of such a reaction is
homogeneous mixture of substances, each of which produced by some allowed
reaction between one (or many) of the reactants in $A$ and one (or many)
of the reactants in $B$.

In our context, the ``mixture'' $A$ is the exponential $\exp\sum\alpha_i
f_i$ of some linear combination of (diagrams representing) functions.
The mixture $B$ is the exponential $\exp\sum \beta_jD_j$ of some sum of
(diagrams representing) mutually inert differential operators $D_j$. That
is, all of the differentiations in the $D_j$'s must act trivially on all
of the coefficients of the $D_j$'s. That is, the diagrams occurring in the
$D_j$'s have ``coefficient legs'' labeled by some set of variables $X$ and
``differentiations legs'' labeled by the dual variables to some disjoint
set of variables $Y$. Computing $B\apply A$ is in some sense analogous
to mixing $A$ and $B$ and allowing them to react.  The result is some
``mixture'' (exponential of a sum) of compounds produced by reactions in
which the legs in some number of the diagrams in $B$ are glued to some of
the legs in some number of the diagrams in $A$. These compounds come with
weights (``densities'') that are (up to minor combinatorial factors) the
products of the densities $\alpha_i$ and $\beta_j$ of their ingredients.

\begin{exercise} \lbl{ex:HomogeneousReactions}
Understand $\left(\exp\sum \beta_jD_j\right)\apply\left(\exp\sum\alpha_i
f_i\right)$ in the above terms.
\end{exercise}

We sometimes consider relabeling operations, where one takes (say) all
legs labeled $x$ in a given diagram $f$ and replaces the $x$ labels
by, say, $y$'s, calling the result $D/(x\to y)$. This corresponds to a
simple change of variable in standard calculus. We wish to allow more
complicated linear reparametrizations as well, but for that we need
to add a bit to the rules of the game.  The added rule is that we also
allow labels that are linear combinations of the basic labels (such as
$x+y$), with the additional provision that the resulting diagrams are
multi-linear in the labels (so a diagram with a leg labeled $x+y$ and
another leg labeled $z+w$ is set equal to a sum of four diagrams labeled
$(x,z)$, $(x,w)$, $(y,z)$, and $(y,w)$). Now reparametrizations such as
$x\to\alpha+\beta$, $y\to\alpha-\beta$ make sense.

\begin{exercise} \lbl{ex:Reparametrization} Show that the operation of
reparametrization is compatible with the application of a differential
operator to a function, as in standard calculus. For instance, in
standard calculus the change of variables $x\to\alpha+\beta$,
$y\to\alpha-\beta$ implies an inverse change for partial derivatives:
$\partial_x\to(\partial_\alpha+\partial_\beta)/2$,
$\partial_y\to(\partial_\alpha-\partial_\beta)/2$. Show that the same
holds in the diagrammatic context:
\[ (D\apply f)\left/\left(
    \begin{array}{c} x\to\alpha+\beta \\ y\to\alpha-\beta \end{array}
  \right)\right. =
  \left(D\left/\left(\begin{array}{c}
    \partial_x\to(\partial_\alpha+\partial_\beta)/2 \\
    \partial_y\to(\partial_\alpha-\partial_\beta)/2
  \end{array}\right)\right.\right)
  \apply \left(f\left/\left(
    \begin{array}{c} x\to\alpha+\beta \\ y\to\alpha-\beta \end{array}
  \right)\right.\right).
\]
\end{exercise}

\subsection{Formal Gaussian integration.} \lbl{subsec:FGI} In standard
calculus, Gaussians are the exponentials of non-degenerate quadratics, and
Gaussian integrals are the integrals of such exponentials multiplied by
polynomials or appropriately convergent power series. Such integrals can be
evaluated using the technique of Feynman diagrams
(see~\cite[appendix]{Aarhus:I}).  The diagrammatic analogs of these
definitions and procedures are described below.

To proceed we must have available the diagrammatic building
blocks for quadratics. These are what we call struts. They
are lines labeled on both ends, they come in an assortment of
forms (figure~\ref{fig:AssortedStruts}) and they satisfy simple
composition laws (figure~\ref{fig:StrutLaws}), that imply that
the strut $\strutf{x}{\partial_x}$ acts like the identity on
legs labeled by $x$, and that the struts $\strutn{x}{x}$ and
$\strutu{\partial_x}{\partial_x}$ are inverses of each other and their
composition is $\strutf{x}{\partial_x}$. Struts always have internal
degree $0$. To make convergence issues simpler below, we assume that
there are no diagrams of internal degree $0$ other than the struts,
and that for any $j\geq0$ there is only a finite number of strutless
diagrams (diagrams none of whose connected components are struts) with
internal degree $\leq j$.

\begin{figure}[htpb]
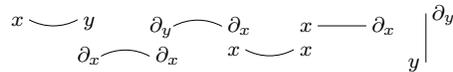

\def\px{{\partial_x}}
\def\py{{\partial_y}}
\[ \printname{AssortedStruts}
	\setlength{\unitlength}{0.5\standardunitlength}
	\begin{array}{c}  \hspace{-1.7mm}
        	\raisebox{-8pt}{\input draws/AssortedStruts.tex }
        	\hspace{-1.9mm}
	\end{array}
 \]
\caption{An assortment of struts.}
\lbl{fig:AssortedStruts}
\end{figure}

\begin{figure}[htpb]
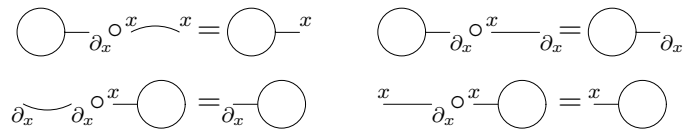

\def\px{{\partial_x}}
\[ \printname{StrutLaws}
	\setlength{\unitlength}{0.5\standardunitlength}
	\begin{array}{c}  \hspace{-1.7mm}
        	\raisebox{-8pt}{\input draws/StrutLaws.tex }
        	\hspace{-1.9mm}
	\end{array}
 \]
\caption{The laws governing strut compositions. The informal notation
$\circ$ means: glue the two adjoining legs.}
\lbl{fig:StrutLaws}
\end{figure}

\begin{definition} Let $X$ be a finite set of variables. A quadratic
$Q$ in the variables in $X$ is a sum of diagrams made of struts whose
ends are labeled by these variables:
\[ Q = \sum_{x,y\in X}l_{xy}\strutn{x}{y}, \]
where the matrix $(l_{xy})$ is symmetric.  Such a quadratic is
``non-degenerate'' if $(l_{xy})$ is invertible. In that case, the
inverse quadratic is the sum
\[ Q^{-1} = \sum_{x,y\in X}l^{xy}\strutu{\partial_x}{\partial_y}, \]
where $(l^{xy})$ denotes the inverse matrix of $(l_{xy})$.
\end{definition}

\begin{definition} An infinite combination of diagrams of the form
\[ G = P\cdot\exp Q/2 \]
is said to be Gaussian with respect to the variables in $X$ if $Q$ is a
quadratic (in those variables) and $P$ is $X$-substantial, meaning that
the diagrams in $P$ have no components which are struts both of whose
ends are labeled by members of $X$. Notice that $P$ and $Q$ are
determined by $G$. In particular, the matrix $\Lambda=(l_{xy})$ of the
coefficients of $Q$ is determined by $G$. We call it the ``covariance
matrix'' of $G$.
\end{definition}

\begin{definition} We say that a Gaussian $G=P\exp Q/2$ is
non-degenerate, or integrable, if $Q$ is non-degenerate. In such a
case, we define the formal Gaussian integral of $G$ with respect to $X$
to be
\begin{equation} \lbl{eq:Gint}
  \Gint P\cdot\exp Q/2\,dX
  = \left\langle \exp -Q^{-1}/2,\,P\right\rangle_X
  = \left(\left(\exp-Q^{-1}/2\right)\apply P\right)
    \left/\left(
      \begin{array}{c} x\to 0 \\ \forall x\in X \end{array}
    \right)\right.,
\end{equation}
where
\[
  \langle D,P \rangle_X :=
  \left(\parbox{2.2in}{
    sum of all ways of gluing the $\partial_x$-marked legs of
    $D$ to the $x$-marked legs of $P$, for all $x\in X$
  }\right).
  \qquad
  \left(\parbox{1.9in}{\footnotesize
    This sum can be non-zero only if the number of $\partial_x$-marked
    legs of $D$ is equal to the number of $x$-marked legs of $P$ for
    all $x\in X$.
  }\right)
\]
(Compare with~\cite[definition~\ref{I-def:Gint} and
equation~\eqref{I-eq:CompactForm}]{Aarhus:I}). The fact that $P$ is
$X$-substantial guarantees that for any given internal degree and any
given number of legs, the computation of the Gaussian integral is
finite.
\end{definition}

Below we need to know some things about the relation between
differentiation and integration. If a certain infinite order differential
operator $D$ contains too many struts, then the computation of (even a
finite part of) $D\apply G$ may be infinite, or else, the result may be
non-Gaussian and thus outside of our theory of integration. Both problems
do not occur if $D$ is ``$X$-substantial'', defined below:

\begin{definition} Let $X$ be a set of variables and $D$ an
differential operator. We say that $D$ is $X$-substantial if it contains no
struts both of whose ends are labeled by members of $X$ or their duals.
\end{definition}

\subsection{Invariance under parity transformations.} It is useful to know
that formal Gaussian integrals, just like their real counterparts, are
invariant under negation of one of the variables:

\begin{proposition} \lbl{prop:Negation} Let $G=P\exp Q/2$ be integrable
with respect to $X$, let $y\in X$, and let $G'=P'\exp Q'/2$ be
$G/(y\to-y)$. Then $\Gint G\,dX = \Gint G'\,dX$.
\end{proposition}

\begin{proof} One easily verifies that $P'=P/(y\to-y)$, $Q'=Q/(y\to-y)$ and
${Q'}^{-1}={Q'}^{-1}/(\partial_y\to-\partial_y)$. Thus in each
($y\leftrightarrow\partial_y$)-gluing in the computation of $\Gint G'\,dX$
two signs get flipped (relative to the computation of $\Gint G\,dX$).
Overall we flip an even number of signs, meaning, no signs at all.
\end{proof}

\begin{exercise} \lbl{ex:ReparametrizeGint}
More generally, verify that Gaussian integration $\Gint$ is compatible
with linear reparametrizations such as in
exercise~\ref{ex:Reparametrization}.
\end{exercise}

\begin{hint} If the reparametrization matrix is $M$, then dual variables
are acted on by $M^{-1}$. Each gluing in~\eqref{eq:Gint} is between one
variable and one dual variable, and so occurrences of $M$ and of $M^{-1}$
come in pairs.
\end{hint}

\subsection{Iterated integration.} \lbl{subsec:IteratedIntegration}
The classical Fubini theorem says that whenever all integrals involved are
well defined, integration over a product space is equivalent to integration
over one factor followed by integration over the other. We seek a similar
iterated integration identity for formal Gaussian integrals.

Let $G=P\exp Q/2$ be a non-degenerate Gaussian with respect to a
set of variables $Z$, and let $Z=X\mathbin\udot Y$ be a decomposition of
the set of variables into two disjoint subsets. Write the covariance
matrix $\Lambda$ of $G$ and its inverse $\Lambda^{-1}$ as block
matrices with respect to this decomposition, taking the variables in
$X$ first and the variables in $Y$ later:
\[ \Lambda=\begin{pmatrix} A & B \\ B^T & C \end{pmatrix},
  \qquad
  \Lambda^{-1} = \begin{pmatrix} D & E \\ E^T & F \end{pmatrix}.
\]
The blocks $A$, $C$, $D$, and $F$ are symmetric, and the fact that
$\Lambda$ and $\Lambda^{-1}$ are inverses implies the following
identities:
\begin{alignat}{2} \lbl{eq:ABCDEFRelations}
  AD+BE^T &= I_X,	& AE+BF &= 0,	\\
  B^TD+CE^T &= 0,	& B^TE+CF &= I_Y. \notag
\end{alignat}

\begin{proposition} \lbl{prop:IteratedIntegral}
If the block $A$ is invertible then the formula
\begin{equation} \lbl{eq:IteratedIntegral}
  \Gint G\,dZ = \Gint\left(\Gint G\,dX\right)\,dY
\end{equation}
makes sense and holds.
\end{proposition}

\begin{proof} The left hand side of this formula is not problematic, and
directly from the definition of formal Gaussian integration and from the
decomposition of $\Lambda^{-1}$ into blocks, we find that it equals
\[ \left\langle\exp-\frac12\left(
    \sum_{x,x'\in X}D^{xx'}\strutu{\partial_x}{\partial_{x'}}
    + 2\sum_{x\in X,\,y\in Y}E^{xy}\strutu{\partial_x}{\partial_y}
    + \sum_{y,y'\in Y}F^{yy'}\strutu{\partial_y}{\partial_{y'}}
  \right),\ P\right\rangle.
\]
We usually suppress the summation symbols, getting
\begin{equation} \lbl{eq:LHS}
  \left\langle\exp-\frac12\left(
    D^{xx'}\strutu{\partial_x}{\partial_{x'}}
    + 2E^{xy}\strutu{\partial_x}{\partial_y}
    + F^{yy'}\strutu{\partial_y}{\partial_{y'}}
  \right),\ P\right\rangle.
\end{equation}

We need to prove that the right hand side
of~\eqref{eq:IteratedIntegral} is well defined and
equals~\eqref{eq:LHS}.  Let us start with the inner integral. Rewriting
the integrand in the form
\[
  \left(
    P\exp\frac12\left(2B_{xy}\strutn{x}{y}+C_{yy'}\strutn{y}{y'}\right)
  \right)
  \exp\frac12A_{xx'}\strutn{x}{x'},
\]
we find that the integrand is Gaussian with respect to $X$ with
covariance matrix $A$. This matrix was assumed to be invertible, and
hence the inner integral $G'$ is defined and equals (denoting the
inverse of $A$ by $\bar{A}$, and using~\eqref{eq:Gint})
\[
  G'=\left.\left(
    \exp-\frac12\left(\bar{A}^{xx'}\strutu{\partial_x}{\partial_{x'}}\right)
    \apply
    \left(
      P\exp\frac12\left(2B_{xy}\strutn{x}{y}+C_{yy'}\strutn{y}{y'}\right)
    \right)
  \right)\right/(x\to 0).
\]
Using exercise~\ref{ex:HomogeneousReactions} and suppressing the
automatic evaluation at $x=0$, this becomes
\[
  \left(
    \exp-\frac12\left(
      \bar{A}^{xx'}\strutu{\partial_x}{\partial_{x'}}
      + 2\bar{A}^{xx_1}B_{x_1y}\strutf{y}{\partial_x}
    \right)\apply P
  \right)
  \exp\frac12\left(
    C_{yy'}\strutn{y}{y'}
    - B^T_{y'x'}\bar{A}^{x'x}B_{xy}\strutn{y}{y'}
  \right).
\]

We are now ready to evaluate the $dY$ integral
of $G'$. In the above formula $G'$ is already
written in the required format $P'\exp\frac12Q'$, with
$Q=C_{yy'}\strutn{y}{y'}-B^T_{y'x'}\bar{A}^{x'x}B_{xy}\strutn{y}{y'}$.
Thus the covariance matrix is $\Lambda'=C-B^T\bar{A}B$. The
relations~\eqref{eq:ABCDEFRelations} imply that $\Lambda'$ is invertible,
with inverse $F$. Thus the integral with respect to $Y$ of $G'$ is
(suppressing the evaluation at $y=0$)
\[
  \left(\exp-\frac12
    F^{yy'}\strutu{\partial_y}{\partial_{y'}}
  \right)
  \apply
  \left(
    \exp-\frac12\left(
      \bar{A}^{xx'}\strutu{\partial_x}{\partial_{x'}}
      + 2\bar{A}^{xx_1}B_{x_1y}\strutf{y}{\partial_x}
    \right)\apply P
  \right).
\]
Again using exercise~\ref{ex:HomogeneousReactions}, this is
\[ =
  \exp-\frac12\left(
    \begin{array}{c}
      \bar{A}^{xx'}\strutu{\partial_x}{\partial_{x'}}
      +\ \bar{A}^{xx_1}B_{x_1y}F^{yy'}B^T_{y'x'_1}\bar{A}^{x'_1x'}
        \strutu{\partial_x}{\partial_{x'}} \\
      \rule{0pt}{16pt}
      -\ 2\bar{A}^{xx'}B_{x'y'}F^{y'y}\strutu{\partial_y}{\partial_x}
      +\ F^{yy'}\strutu{\partial_y}{\partial_{y'}}
    \end{array}
  \right) \apply P.
\]
Giving names to the coefficients and switching to matrix-talk, we find that
this is
\[ 
  \exp-\frac12\left(
    L^{xx'}\strutu{\partial_x}{\partial_{x'}}
    + 2M^{xy}\strutu{\partial_x}{\partial_y}
    + F^{yy'}\strutu{\partial_y}{\partial_{y'}}
  \right)\apply P,
\]
with $L=\bar{A}+\bar{A}BFB^T\bar{A}^T$ and $M=-\bar{A}BF$. A second
look at the relations~\eqref{eq:ABCDEFRelations} reveals that $L=D$ and
$M=E$, proving that the last formula is equal to~\eqref{eq:LHS}, as
required.
\end{proof}

\subsection{Integration by parts.} \lbl{subsec:IBP}
Let $D$ be a diagram representing a differential operator with respect
to the variable $z$ (that is, it may have ``differentiation legs''
labeled $\partial_{z}$, ``coefficient legs'' labeled ${z}$, and
possibly other legs labeled by other variables). Assume $D$ has $l$
differentiation legs labeled $\partial_{z}$ (``$D$ is of order $l$'')
and $k$ coefficient legs labeled $z$.

\begin{definition} The ``divergence'' $\operatorname{div}_{z}D$ of $D$
with respect to ${z}$ is the result of ``applying $D$ to its own
coefficients''. That is,
\[ \operatorname{div}_{z}D=\begin{cases}
  0 & \text{if }l>k, \\
  \left(\parbox{3.1in}{\begin{center}
    sum of all ways of attaching all legs labeled $\partial_{z}$ with
    some or all legs labeled ${z}$
  \end{center}}\right), & \text{if }l\leq k.
  \end{cases}
\]
(Compare with the standard definition of the divergence of a vector field,
where each derivative ``turns back'' and acts on its own coefficient).
\end{definition}

In standard calculus, the following proposition is an easy consequence of
integration by parts:

\begin{proposition} \lbl{prop:Divergence} Let $X$ be a set of
variables, and let $z\in X$. If $G$ is a non-degenerate 
Gaussian with respect to $X$ and $D$ is an $X$-substantial operator of
order $l$, then
\begin{equation} \lbl{eq:DivergenceFormula}
  \Gint D\apply G\,dX
  = (-1)^l\Gint(\operatorname{div}_{z}D)G\,dX.
\end{equation}
\end{proposition}

\begin{proof} Write $G=P\exp Q/2$ with $Q=\sum_{x,y\in
X}l_{xy}\strutn{x}{y}$. Let us pick one leg marked $\partial_{z}$ in
$D$, put a little asterisk ($\ast$) on it, and follow it throughout the
computation of the left hand side of~\eqref{eq:DivergenceFormula}.
First, in computing $D\apply G$, the special leg gets glued either to
one of the legs in $P$, or to one of the legs in $\exp Q/2$. The result
looks something like
\[
  \def\px{{\partial_{z}}}
  \def\ysum{{\displaystyle\sum_{y\in X}l_{zy}}}
  D\apply G = \left(\printname{DG}
	\setlength{\unitlength}{0.5\standardunitlength}
	\begin{array}{c}  \hspace{-1.7mm}
        	\raisebox{-8pt}{\input draws/DG.tex }
        	\hspace{-1.9mm}
	\end{array}
\right)\exp Q/2.
\]

The next step is integration. The factor $\exp Q/2$ is removed, and
struts labeled and weighted by the negated inverse covariance matrix are
glued in. Of particular interest is the strut $\strutu{y}{y'}$ glued
to the marked leg in the right term. It comes with a coefficient like
$-l^{yy'}$ from the negated inverse covariance matrix, which multiplies
the coefficient $l_{zy}$ already in place. The summation over $y$
evaluates matrix multiplication of a matrix and its inverse, and we
find that $y'=z$ and the overall coefficient is $-1$.  The other end
of this strut is glued to some leg (with a label in $X$), either on $P$
or on $D$.  Following that, all other negated inverse covariance gluings
are performed. The result looks something like

{
  \def\px{{\partial_{z}}}
  \begin{align*}
    \Gint D\apply G & dX
    = \left(
        \exp-\frac12\sum_{x,y\in X}l^{xy}\strutu{\partial_x}{\partial_y}
      \right) \apply \\
    & \left(\printname{GintDG}
	\setlength{\unitlength}{0.5\standardunitlength}
	\begin{array}{c}  \hspace{-1.7mm}
        	\raisebox{-8pt}{\input draws/GintDG.tex }
        	\hspace{-1.9mm}
	\end{array}
\right).
  \end{align*}
}
In this formula the first term cancels the second, and we are left only
with the third. But the same argument can be made for all legs marked
$\partial_{z}$ in $D$, and hence in left hand side integral in
equation~\eqref{eq:DivergenceFormula} they all have to ``turn back'' and
differentiate a coefficient of $D$. Counting signs, this is precisely the
right hand side of equation~\eqref{eq:DivergenceFormula}.
\end{proof}

\subsection{Our formal universe.} \lbl{sebsec:OurUniverse}
Below we apply the formalism and techniques developed in this section in
the case where the diagrams are $X$-marked uni-trivalent diagrams modulo the
$AS$ and $IHX$ relations, for some label set $X$. The ``internal degree''
of a diagram is half the number of internal (trivalent) vertices it has,
and the struts are simply the internal degree $0$ diagrams ---
uni-trivalent diagrams that have no internal vertices. In the Kontsevich integral,
the coefficients of struts measure linking numbers between the components
marked at their ends (or self linkings, if the two ends are marked the
same way).

The spaces of uni-trivalent diagrams that we consider are Hopf algebras,
and the formal linear combinations of uni-trivalent diagrams that we
take as inputs come from evaluating the Kontsevich integral, whose
values are always grouplike. Thus our inputs are always exponentials.
Splitting away the struts, we find that they are always Gaussian with
the linking matrix of the underlying link as the covariance matrix. If
we stick to pure tangles whose linking matrix is non-singular, our inputs
are always integrable.

%% file: draws/DiagramPolynomial.tex
%
\begingroup\makeatletter\ifx\SetFigFont\undefined
\def\x#1#2#3#4#5#6#7\relax{\def\x{#1#2#3#4#5#6}}%
\expandafter\x\fmtname xxxxxx\relax \def\y{splain}%
\ifx\x\y   
\gdef\SetFigFont#1#2#3{%
  \ifnum #1<17\tiny\else \ifnum #1<20\small\else
  \ifnum #1<24\normalsize\else \ifnum #1<29\large\else
  \ifnum #1<34\Large\else \ifnum #1<41\LARGE\else
     \huge\fi\fi\fi\fi\fi\fi
  \csname #3\endcsname}%
\else
\gdef\SetFigFont#1#2#3{\begingroup
  \count@#1\relax \ifnum 25<\count@\count@25\fi
  \def\x{\endgroup\@setsize\SetFigFont{#2pt}}%
  \expandafter\x
    \csname \romannumeral\the\count@ pt\expandafter\endcsname
    \csname @\romannumeral\the\count@ pt\endcsname
  \csname #3\endcsname}%
\fi
\fi\endgroup
{\renewcommand{\dashlinestretch}{30}
\begin{picture}(1909,1413)(0,-10)
\path(1650,1302)(1050,702)
\path(1050,702)(1650,102)
\path(1050,702)(150,1302)
\path(1050,702)(150,702)
\path(1050,702)(150,102)
\put(0,1302){\makebox(0,0)[lb]{\smash{{\mathmode{\scriptstyle x}}}}}
\put(0,702){\makebox(0,0)[lb]{\smash{{\mathmode{\scriptstyle x}}}}}
\put(1725,1227){\makebox(0,0)[lb]{\smash{{\mathmode{\scriptstyle y}}}}}
\put(1725,27){\makebox(0,0)[lb]{\smash{{\mathmode{\scriptstyle y}}}}}
\put(0,102){\makebox(0,0)[lb]{\smash{{\mathmode{\scriptstyle x}}}}}
\end{picture}
}

%% file: draws/DiagramOperator.tex
%
\begingroup\makeatletter\ifx\SetFigFont\undefined
\def\x#1#2#3#4#5#6#7\relax{\def\x{#1#2#3#4#5#6}}%
\expandafter\x\fmtname xxxxxx\relax \def\y{splain}%
\ifx\x\y   
\gdef\SetFigFont#1#2#3{%
  \ifnum #1<17\tiny\else \ifnum #1<20\small\else
  \ifnum #1<24\normalsize\else \ifnum #1<29\large\else
  \ifnum #1<34\Large\else \ifnum #1<41\LARGE\else
     \huge\fi\fi\fi\fi\fi\fi
  \csname #3\endcsname}%
\else
\gdef\SetFigFont#1#2#3{\begingroup
  \count@#1\relax \ifnum 25<\count@\count@25\fi
  \def\x{\endgroup\@setsize\SetFigFont{#2pt}}%
  \expandafter\x
    \csname \romannumeral\the\count@ pt\expandafter\endcsname
    \csname @\romannumeral\the\count@ pt\endcsname
  \csname #3\endcsname}%
\fi
\fi\endgroup
{\renewcommand{\dashlinestretch}{30}
\begin{picture}(2494,1422)(0,-10)
\path(1650,1311)(1050,711)
\path(1050,711)(1650,111)
\path(1050,711)(150,1311)
\path(1050,711)(150,711)
\path(1050,711)(150,111)
\put(0,1311){\makebox(0,0)[lb]{\smash{{\mathmode{\scriptstyle x}}}}}
\put(0,711){\makebox(0,0)[lb]{\smash{{\mathmode{\scriptstyle x}}}}}
\put(0,111){\makebox(0,0)[lb]{\smash{{\mathmode{\scriptstyle x}}}}}
\put(1725,1236){\makebox(0,0)[lb]{\smash{{\mathmode{\scriptstyle \partial_y}}}}}
\put(1725,36){\makebox(0,0)[lb]{\smash{{\mathmode{\scriptstyle \partial_y}}}}}
\end{picture}
}

%% file: draws/CCOperator.tex
%
\begingroup\makeatletter\ifx\SetFigFont\undefined
\def\x#1#2#3#4#5#6#7\relax{\def\x{#1#2#3#4#5#6}}%
\expandafter\x\fmtname xxxxxx\relax \def\y{splain}%
\ifx\x\y   
\gdef\SetFigFont#1#2#3{%
  \ifnum #1<17\tiny\else \ifnum #1<20\small\else
  \ifnum #1<24\normalsize\else \ifnum #1<29\large\else
  \ifnum #1<34\Large\else \ifnum #1<41\LARGE\else
     \huge\fi\fi\fi\fi\fi\fi
  \csname #3\endcsname}%
\else
\gdef\SetFigFont#1#2#3{\begingroup
  \count@#1\relax \ifnum 25<\count@\count@25\fi
  \def\x{\endgroup\@setsize\SetFigFont{#2pt}}%
  \expandafter\x
    \csname \romannumeral\the\count@ pt\expandafter\endcsname
    \csname @\romannumeral\the\count@ pt\endcsname
  \csname #3\endcsname}%
\fi
\fi\endgroup
{\renewcommand{\dashlinestretch}{30}
\begin{picture}(1969,846)(0,-10)
\put(675.000,612.000){\arc{1200.000}{6.2832}{9.4248}}
\path(450,612)(675,12)(900,612)
\put(0,687){\makebox(0,0)[lb]{\smash{{\mathmode{\scriptstyle \partial_x}}}}}
\put(375,687){\makebox(0,0)[lb]{\smash{{\mathmode{\scriptstyle \partial_x}}}}}
\put(825,687){\makebox(0,0)[lb]{\smash{{\mathmode{\scriptstyle \partial_x}}}}}
\put(1200,687){\makebox(0,0)[lb]{\smash{{\mathmode{\scriptstyle \partial_x}}}}}
\end{picture}
}

%% file: draws/DfExample.tex
\begingroup\makeatletter\ifx\SetFigFont\undefined
\def\x#1#2#3#4#5#6#7\relax{\def\x{#1#2#3#4#5#6}}%
\expandafter\x\fmtname xxxxxx\relax \def\y{splain}%
\ifx\x\y   
\gdef\SetFigFont#1#2#3{%
  \ifnum #1<17\tiny\else \ifnum #1<20\small\else
  \ifnum #1<24\normalsize\else \ifnum #1<29\large\else
  \ifnum #1<34\Large\else \ifnum #1<41\LARGE\else
     \huge\fi\fi\fi\fi\fi\fi
  \csname #3\endcsname}%
\else
\gdef\SetFigFont#1#2#3{\begingroup
  \count@#1\relax \ifnum 25<\count@\count@25\fi
  \def\x{\endgroup\@setsize\SetFigFont{#2pt}}%
  \expandafter\x
    \csname \romannumeral\the\count@ pt\expandafter\endcsname
    \csname @\romannumeral\the\count@ pt\endcsname
  \csname #3\endcsname}%
\fi
\fi\endgroup
{\renewcommand{\dashlinestretch}{30}
\begin{picture}(11514,1723)(0,-10)
\path(150,926)(450,926)
\path(450,926)(750,1151)
\path(450,926)(750,701)
\path(1651,1226)(2101,926)
\path(1651,926)(2101,926)
\path(1651,626)(2101,926)
\path(2101,926)(2551,1151)
\path(2101,926)(2551,701)
\path(4051,1375)(4351,1375)
\path(4351,1375)(4651,1600)
\path(4351,1375)(4651,1150)
\path(4953,1675)(5403,1375)
\path(4953,1375)(5403,1375)
\path(4953,1075)(5403,1375)
\path(5403,1375)(5853,1600)
\path(5403,1375)(5853,1150)
\path(4051,475)(4351,475)
\path(4351,475)(4651,700)
\path(4351,475)(4651,250)
\path(4953,775)(5403,475)
\path(4953,475)(5403,475)
\path(4953,175)(5403,475)
\path(5403,475)(5853,700)
\path(5403,475)(5853,250)
\path(6753,1375)(7053,1375)
\path(7053,1375)(7353,1600)
\path(7053,1375)(7353,1150)
\path(7655,1675)(8105,1375)
\path(7655,1375)(8105,1375)
\path(7655,1075)(8105,1375)
\path(8105,1375)(8555,1600)
\path(8105,1375)(8555,1150)
\path(6753,475)(7053,475)
\path(7053,475)(7353,700)
\path(7053,475)(7353,250)
\path(7655,775)(8105,475)
\path(7655,475)(8105,475)
\path(7655,175)(8105,475)
\path(8105,475)(8555,700)
\path(8105,475)(8555,250)
\path(4649,1602)(4949,1377)
\path(4649,1152)(4949,1077)
\path(7349,1602)(7649,1677)
\path(7349,1152)(7649,1077)
\path(10049,1602)(10349,1677)
\path(10049,1152)(10349,1377)
\path(4649,702)(4949,177)
\path(4649,252)(4949,477)
\path(7349,702)(7649,177)
\path(7349,252)(7649,777)
\path(10049,702)(10349,477)
\path(10049,252)(10349,777)
\path(9453,1375)(9753,1375)
\path(9753,1375)(10053,1600)
\path(9753,1375)(10053,1150)
\path(10355,1675)(10805,1375)
\path(10355,1375)(10805,1375)
\path(10355,1075)(10805,1375)
\path(10805,1375)(11255,1600)
\path(10805,1375)(11255,1150)
\path(9453,475)(9753,475)
\path(9753,475)(10053,700)
\path(9753,475)(10053,250)
\path(10355,775)(10805,475)
\path(10355,475)(10805,475)
\path(10355,175)(10805,475)
\path(10805,475)(11255,700)
\path(10805,475)(11255,250)
\path(149,477)	(192.729,424.564)
	(224.000,402.000)

\path(224,402)	(275.598,392.204)
	(338.300,393.293)
	(400.102,398.735)
	(449.000,402.000)

\path(449,402)	(497.898,398.735)
	(559.700,393.293)
	(622.402,392.204)
	(674.000,402.000)

\path(674,402)	(705.271,424.564)
	(749.000,477.000)

\path(1649,477)	(1692.729,424.564)
	(1724.000,402.000)

\path(1724,402)	(1763.434,390.889)
	(1809.996,385.673)
	(1861.185,384.993)
	(1914.500,387.488)
	(1967.440,391.796)
	(2017.504,396.558)
	(2062.191,400.413)
	(2099.000,402.000)

\path(2099,402)	(2135.809,400.413)
	(2180.496,396.558)
	(2230.560,391.796)
	(2283.500,387.488)
	(2336.815,384.993)
	(2388.004,385.673)
	(2434.566,390.889)
	(2474.000,402.000)

\path(2474,402)	(2505.271,424.564)
	(2549.000,477.000)

\put(0,851){\makebox(0,0)[lb]{\smash{{\mathmode{\scriptstyle x}}}}}
\put(750,1076){\makebox(0,0)[lb]{\smash{{\mathmode{\scriptstyle \px}}}}}
\put(750,626){\makebox(0,0)[lb]{\smash{{\mathmode{\scriptstyle \px}}}}}
\put(374,27){\makebox(0,0)[lb]{\smash{{\mathmode{D}}}}}
\put(2999,852){\makebox(0,0)[lb]{\smash{{\mathmode{=}}}}}
\put(1501,1151){\makebox(0,0)[lb]{\smash{{\mathmode{\scriptstyle x}}}}}
\put(1501,851){\makebox(0,0)[lb]{\smash{{\mathmode{\scriptstyle x}}}}}
\put(1501,551){\makebox(0,0)[lb]{\smash{{\mathmode{\scriptstyle x}}}}}
\put(2626,1151){\makebox(0,0)[lb]{\smash{{\mathmode{\scriptstyle y}}}}}
\put(2626,701){\makebox(0,0)[lb]{\smash{{\mathmode{\scriptstyle y}}}}}
\put(3449,402){\makebox(0,0)[lb]{\smash{{\mathmode{+}}}}}
\put(2024,27){\makebox(0,0)[lb]{\smash{{\mathmode{f}}}}}
\put(6150,1302){\makebox(0,0)[lb]{\smash{{\mathmode{+}}}}}
\put(6150,402){\makebox(0,0)[lb]{\smash{{\mathmode{+}}}}}
\put(8925,1302){\makebox(0,0)[lb]{\smash{{\mathmode{+}}}}}
\put(8925,402){\makebox(0,0)[lb]{\smash{{\mathmode{+}}}}}
\put(3901,1300){\makebox(0,0)[lb]{\smash{{\mathmode{\scriptstyle x}}}}}
\put(5928,1600){\makebox(0,0)[lb]{\smash{{\mathmode{\scriptstyle y}}}}}
\put(5928,1150){\makebox(0,0)[lb]{\smash{{\mathmode{\scriptstyle y}}}}}
\put(3901,400){\makebox(0,0)[lb]{\smash{{\mathmode{\scriptstyle x}}}}}
\put(5928,700){\makebox(0,0)[lb]{\smash{{\mathmode{\scriptstyle y}}}}}
\put(5928,250){\makebox(0,0)[lb]{\smash{{\mathmode{\scriptstyle y}}}}}
\put(6603,1300){\makebox(0,0)[lb]{\smash{{\mathmode{\scriptstyle x}}}}}
\put(8630,1600){\makebox(0,0)[lb]{\smash{{\mathmode{\scriptstyle y}}}}}
\put(8630,1150){\makebox(0,0)[lb]{\smash{{\mathmode{\scriptstyle y}}}}}
\put(6603,400){\makebox(0,0)[lb]{\smash{{\mathmode{\scriptstyle x}}}}}
\put(8630,700){\makebox(0,0)[lb]{\smash{{\mathmode{\scriptstyle y}}}}}
\put(8630,250){\makebox(0,0)[lb]{\smash{{\mathmode{\scriptstyle y}}}}}
\put(4799,702){\makebox(0,0)[lb]{\smash{{\mathmode{\scriptstyle x}}}}}
\put(4799,1602){\makebox(0,0)[lb]{\smash{{\mathmode{\scriptstyle x}}}}}
\put(7499,1302){\makebox(0,0)[lb]{\smash{{\mathmode{\scriptstyle x}}}}}
\put(7499,402){\makebox(0,0)[lb]{\smash{{\mathmode{\scriptstyle x}}}}}
\put(10199,1002){\makebox(0,0)[lb]{\smash{{\mathmode{\scriptstyle x}}}}}
\put(10199,102){\makebox(0,0)[lb]{\smash{{\mathmode{\scriptstyle x}}}}}
\put(9303,1300){\makebox(0,0)[lb]{\smash{{\mathmode{\scriptstyle x}}}}}
\put(11330,1600){\makebox(0,0)[lb]{\smash{{\mathmode{\scriptstyle y}}}}}
\put(11330,1150){\makebox(0,0)[lb]{\smash{{\mathmode{\scriptstyle y}}}}}
\put(9303,400){\makebox(0,0)[lb]{\smash{{\mathmode{\scriptstyle x}}}}}
\put(11330,700){\makebox(0,0)[lb]{\smash{{\mathmode{\scriptstyle y}}}}}
\put(11330,250){\makebox(0,0)[lb]{\smash{{\mathmode{\scriptstyle y}}}}}
\put(1198,852){\makebox(0,0)[lb]{\smash{{\mathmode{\apply}}}}}
\end{picture}
}

%% file: draws/AssortedStruts.tex
\begingroup\makeatletter\ifx\SetFigFont\undefined
\def\x#1#2#3#4#5#6#7\relax{\def\x{#1#2#3#4#5#6}}%
\expandafter\x\fmtname xxxxxx\relax \def\y{splain}%
\ifx\x\y   
\gdef\SetFigFont#1#2#3{%
  \ifnum #1<17\tiny\else \ifnum #1<20\small\else
  \ifnum #1<24\normalsize\else \ifnum #1<29\large\else
  \ifnum #1<34\Large\else \ifnum #1<41\LARGE\else
     \huge\fi\fi\fi\fi\fi\fi
  \csname #3\endcsname}%
\else
\gdef\SetFigFont#1#2#3{\begingroup
  \count@#1\relax \ifnum 25<\count@\count@25\fi
  \def\x{\endgroup\@setsize\SetFigFont{#2pt}}%
  \expandafter\x
    \csname \romannumeral\the\count@ pt\expandafter\endcsname
    \csname @\romannumeral\the\count@ pt\endcsname
  \csname #3\endcsname}%
\fi
\fi\endgroup
{\renewcommand{\dashlinestretch}{30}
\begin{picture}(5573,786)(0,-10)
\put(2325.000,-10.500){\arc{1275.000}{4.2224}{5.2023}}
\put(525.000,1189.500){\arc{1275.000}{1.0808}{2.0608}}
\put(1425.000,-385.500){\arc{1275.000}{4.2224}{5.2023}}
\put(3225.000,814.500){\arc{1275.000}{1.0808}{2.0608}}
\path(3825,552)(4425,552)
\path(5175,702)(5175,102)
\put(900,552){\makebox(0,0)[lb]{\smash{{\mathmode{\scriptstyle y}}}}}
\put(1800,102){\makebox(0,0)[lb]{\smash{{\mathmode{\scriptstyle \px}}}}}
\put(2700,477){\makebox(0,0)[lb]{\smash{{\mathmode{\scriptstyle \px}}}}}
\put(3600,177){\makebox(0,0)[lb]{\smash{{\mathmode{\scriptstyle x}}}}}
\put(4500,477){\makebox(0,0)[lb]{\smash{{\mathmode{\scriptstyle \px}}}}}
\put(825,102){\makebox(0,0)[lb]{\smash{{\mathmode{\scriptstyle \px}}}}}
\put(1725,477){\makebox(0,0)[lb]{\smash{{\mathmode{\scriptstyle \py}}}}}
\put(5250,627){\makebox(0,0)[lb]{\smash{{\mathmode{\scriptstyle \py}}}}}
\put(4950,27){\makebox(0,0)[lb]{\smash{{\mathmode{\scriptstyle y}}}}}
\put(0,552){\makebox(0,0)[lb]{\smash{{\mathmode{\scriptstyle x}}}}}
\put(3600,477){\makebox(0,0)[lb]{\smash{{\mathmode{\scriptstyle x}}}}}
\put(2700,177){\makebox(0,0)[lb]{\smash{{\mathmode{\scriptstyle x}}}}}
\end{picture}
}

%% file: draws/StrutLaws.tex
\begingroup\makeatletter\ifx\SetFigFont\undefined
\def\x#1#2#3#4#5#6#7\relax{\def\x{#1#2#3#4#5#6}}%
\expandafter\x\fmtname xxxxxx\relax \def\y{splain}%
\ifx\x\y   
\gdef\SetFigFont#1#2#3{%
  \ifnum #1<17\tiny\else \ifnum #1<20\small\else
  \ifnum #1<24\normalsize\else \ifnum #1<29\large\else
  \ifnum #1<34\Large\else \ifnum #1<41\LARGE\else
     \huge\fi\fi\fi\fi\fi\fi
  \csname #3\endcsname}%
\else
\gdef\SetFigFont#1#2#3{\begingroup
  \count@#1\relax \ifnum 25<\count@\count@25\fi
  \def\x{\endgroup\@setsize\SetFigFont{#2pt}}%
  \expandafter\x
    \csname \romannumeral\the\count@ pt\expandafter\endcsname
    \csname @\romannumeral\the\count@ pt\endcsname
  \csname #3\endcsname}%
\fi
\fi\endgroup
{\renewcommand{\dashlinestretch}{30}
\begin{picture}(8421,1529)(0,-10)
\put(1799.000,644.500){\arc{1275.000}{4.2224}{5.2023}}
\put(450.000,869.500){\arc{1275.000}{1.0808}{2.0608}}
\put(2999,1207){\ellipse{600}{600}}
\put(374,1207){\ellipse{600}{600}}
\put(1875,307){\ellipse{600}{600}}
\put(3375,307){\ellipse{600}{600}}
\put(7499,1207){\ellipse{600}{600}}
\put(4874,1207){\ellipse{600}{600}}
\put(6375,307){\ellipse{600}{600}}
\put(7875,307){\ellipse{600}{600}}
\path(3299,1207)(3599,1207)
\path(674,1207)(974,1207)
\path(1275,307)(1575,307)
\path(2775,307)(3075,307)
\path(7799,1207)(8099,1207)
\path(5174,1207)(5474,1207)
\path(5775,307)(6075,307)
\path(7275,307)(7575,307)
\path(4648,307)(5248,307)
\path(6598,1207)(5998,1207)
\put(3599,1282){\makebox(0,0)[lb]{\smash{{\mathmode{\scriptstyle x}}}}}
\put(1199,1132){\makebox(0,0)[lb]{\smash{{\mathmode{\circ}}}}}
\put(2324,1132){\makebox(0,0)[lb]{\smash{{\mathmode{=}}}}}
\put(974,982){\makebox(0,0)[lb]{\smash{{\mathmode{\scriptstyle \px}}}}}
\put(2099,1282){\makebox(0,0)[lb]{\smash{{\mathmode{\scriptstyle x}}}}}
\put(1424,1282){\makebox(0,0)[lb]{\smash{{\mathmode{\scriptstyle x}}}}}
\put(2324,232){\makebox(0,0)[lb]{\smash{{\mathmode{=}}}}}
\put(975,232){\makebox(0,0)[lb]{\smash{{\mathmode{\circ}}}}}
\put(1200,382){\makebox(0,0)[lb]{\smash{{\mathmode{\scriptstyle x}}}}}
\put(750,82){\makebox(0,0)[lb]{\smash{{\mathmode{\scriptstyle \px}}}}}
\put(0,82){\makebox(0,0)[lb]{\smash{{\mathmode{\scriptstyle \px}}}}}
\put(2625,82){\makebox(0,0)[lb]{\smash{{\mathmode{\scriptstyle \px}}}}}
\put(5699,1132){\makebox(0,0)[lb]{\smash{{\mathmode{\circ}}}}}
\put(6824,1132){\makebox(0,0)[lb]{\smash{{\mathmode{=}}}}}
\put(5474,982){\makebox(0,0)[lb]{\smash{{\mathmode{\scriptstyle \px}}}}}
\put(5924,1282){\makebox(0,0)[lb]{\smash{{\mathmode{\scriptstyle x}}}}}
\put(6824,232){\makebox(0,0)[lb]{\smash{{\mathmode{=}}}}}
\put(5475,232){\makebox(0,0)[lb]{\smash{{\mathmode{\circ}}}}}
\put(5700,382){\makebox(0,0)[lb]{\smash{{\mathmode{\scriptstyle x}}}}}
\put(5250,82){\makebox(0,0)[lb]{\smash{{\mathmode{\scriptstyle \px}}}}}
\put(4573,382){\makebox(0,0)[lb]{\smash{{\mathmode{\scriptstyle x}}}}}
\put(7198,382){\makebox(0,0)[lb]{\smash{{\mathmode{\scriptstyle x}}}}}
\put(6598,982){\makebox(0,0)[lb]{\smash{{\mathmode{\scriptstyle \px}}}}}
\put(8098,982){\makebox(0,0)[lb]{\smash{{\mathmode{\scriptstyle \px}}}}}
\end{picture}
}

%% file: draws/DG.tex
\begingroup\makeatletter\ifx\SetFigFont\undefined
\def\x#1#2#3#4#5#6#7\relax{\def\x{#1#2#3#4#5#6}}%
\expandafter\x\fmtname xxxxxx\relax \def\y{splain}%
\ifx\x\y   
\gdef\SetFigFont#1#2#3{%
  \ifnum #1<17\tiny\else \ifnum #1<20\small\else
  \ifnum #1<24\normalsize\else \ifnum #1<29\large\else
  \ifnum #1<34\Large\else \ifnum #1<41\LARGE\else
     \huge\fi\fi\fi\fi\fi\fi
  \csname #3\endcsname}%
\else
\gdef\SetFigFont#1#2#3{\begingroup
  \count@#1\relax \ifnum 25<\count@\count@25\fi
  \def\x{\endgroup\@setsize\SetFigFont{#2pt}}%
  \expandafter\x
    \csname \romannumeral\the\count@ pt\expandafter\endcsname
    \csname @\romannumeral\the\count@ pt\endcsname
  \csname #3\endcsname}%
\fi
\fi\endgroup
{\renewcommand{\dashlinestretch}{30}
\begin{picture}(8043,1539)(0,-10)
\put(6530.000,657.895){\arc{1260.210}{4.2162}{5.2086}}
\path(450,1363)(1050,1363)(1050,163)
	(450,163)(450,1363)
\path(450,1063)(150,1063)
\path(2551,1512)(3151,1512)(3151,12)
	(2551,12)(2551,1512)
\path(450,763)(150,763)
\path(450,463)(150,463)
\path(1050,912)(1200,912)
\path(1050,312)(1200,312)
\path(2401,161)(2551,161)
\path(2401,1062)(2551,1062)
\path(1201,1212)(2401,1363)
\path(1051,1212)(1201,1212)
\path(2402,1363)(2552,1363)
\path(5330,1363)(5930,1363)(5930,163)
	(5330,163)(5330,1363)
\path(5330,1063)(5030,1063)
\path(7431,1512)(8031,1512)(8031,12)
	(7431,12)(7431,1512)
\path(5330,763)(5030,763)
\path(5330,463)(5030,463)
\path(5930,1213)(6226,1213)
\path(5930,312)(6080,312)
\path(5930,912)(6076,912)
\path(7282,1362)(7432,1362)
\path(7281,161)(7431,161)
\path(7281,1062)(7431,1062)
\put(3300,688){\makebox(0,0)[lb]{\smash{{\mathmode{+\ysum}}}}}
\put(600,688){\makebox(0,0)[lb]{\smash{{\mathmode{D}}}}}
\put(2700,688){\makebox(0,0)[lb]{\smash{{\mathmode{P}}}}}
\put(2401,387){\makebox(0,0)[lb]{\smash{{\mathmode{\vdots}}}}}
\put(1126,987){\makebox(0,0)[lb]{\smash{{\mathmode{\scriptstyle \ast\px}}}}}
\put(1500,688){\makebox(0,0)[lb]{\smash{{{\rm\scriptsize more}}}}}
\put(1125,462){\makebox(0,0)[lb]{\smash{{\mathmode{\vdots}}}}}
\put(0,1138){\makebox(0,0)[lb]{\smash{{\mathmode{\scriptstyle z}}}}}
\put(0,838){\makebox(0,0)[lb]{\smash{{\mathmode{\scriptstyle z}}}}}
\put(0,538){\makebox(0,0)[lb]{\smash{{\mathmode{\scriptstyle w}}}}}
\put(2250,1213){\makebox(0,0)[lb]{\smash{{\mathmode{\scriptstyle z}}}}}
\put(5476,688){\makebox(0,0)[lb]{\smash{{\mathmode{D}}}}}
\put(7576,688){\makebox(0,0)[lb]{\smash{{\mathmode{P}}}}}
\put(6004,987){\makebox(0,0)[lb]{\smash{{\mathmode{\scriptstyle \ast\px}}}}}
\put(6005,462){\makebox(0,0)[lb]{\smash{{\mathmode{\vdots}}}}}
\put(7276,387){\makebox(0,0)[lb]{\smash{{\mathmode{\vdots}}}}}
\put(6376,688){\makebox(0,0)[lb]{\smash{{{\rm\scriptsize more}}}}}
\put(4876,1138){\makebox(0,0)[lb]{\smash{{\mathmode{\scriptstyle z}}}}}
\put(4876,838){\makebox(0,0)[lb]{\smash{{\mathmode{\scriptstyle z}}}}}
\put(4876,538){\makebox(0,0)[lb]{\smash{{\mathmode{\scriptstyle w}}}}}
\put(6301,1288){\makebox(0,0)[lb]{\smash{{\mathmode{\scriptstyle z}}}}}
\put(6751,1288){\makebox(0,0)[lb]{\smash{{\mathmode{\scriptstyle y}}}}}
\put(7126,1213){\makebox(0,0)[lb]{\smash{{\mathmode{\scriptstyle z}}}}}
\put(1350,388){\makebox(0,0)[lb]{\smash{{{\rm\scriptsize activity}}}}}
\put(6225,388){\makebox(0,0)[lb]{\smash{{{\rm\scriptsize activity}}}}}
\end{picture}
}

%% file: draws/GintDG.tex
\begingroup\makeatletter\ifx\SetFigFont\undefined
\def\x#1#2#3#4#5#6#7\relax{\def\x{#1#2#3#4#5#6}}%
\expandafter\x\fmtname xxxxxx\relax \def\y{splain}%
\ifx\x\y   
\gdef\SetFigFont#1#2#3{%
  \ifnum #1<17\tiny\else \ifnum #1<20\small\else
  \ifnum #1<24\normalsize\else \ifnum #1<29\large\else
  \ifnum #1<34\Large\else \ifnum #1<41\LARGE\else
     \huge\fi\fi\fi\fi\fi\fi
  \csname #3\endcsname}%
\else
\gdef\SetFigFont#1#2#3{\begingroup
  \count@#1\relax \ifnum 25<\count@\count@25\fi
  \def\x{\endgroup\@setsize\SetFigFont{#2pt}}%
  \expandafter\x
    \csname \romannumeral\the\count@ pt\expandafter\endcsname
    \csname @\romannumeral\the\count@ pt\endcsname
  \csname #3\endcsname}%
\fi
\fi\endgroup
{\renewcommand{\dashlinestretch}{30}
\begin{picture}(11797,1720)(0,-10)
\put(3376,688){\makebox(0,0)[lb]{\smash{{\mathmode{-}}}}}
\put(7203,688){\makebox(0,0)[lb]{\smash{{\mathmode{-}}}}}
\path(450,1363)(1050,1363)(1050,163)
	(450,163)(450,1363)
\path(450,1063)(150,1063)
\path(2551,1512)(3151,1512)(3151,12)
	(2551,12)(2551,1512)
\path(450,763)(150,763)
\path(450,463)(150,463)
\path(1050,912)(1200,912)
\path(1050,312)(1200,312)
\path(2401,161)(2551,161)
\path(2401,1062)(2551,1062)
\path(1201,1212)(2401,1363)
\path(1051,1212)(1201,1212)
\path(2402,1363)(2552,1363)
\path(4278,1363)(4878,1363)(4878,163)
	(4278,163)(4278,1363)
\path(4278,1063)(3978,1063)
\path(6379,1512)(6979,1512)(6979,12)
	(6379,12)(6379,1512)
\path(4278,763)(3978,763)
\path(4278,463)(3978,463)
\path(4878,912)(5028,912)
\path(4878,312)(5028,312)
\path(6229,161)(6379,161)
\path(6229,1062)(6379,1062)
\path(5029,1212)(6229,1363)
\path(4879,1212)(5029,1212)
\path(6230,1363)(6380,1363)
\path(9084,1363)(9684,1363)(9684,163)
	(9084,163)(9084,1363)
\path(9084,1063)(8784,1063)
\path(11185,1512)(11785,1512)(11785,12)
	(11185,12)(11185,1512)
\path(9084,763)(8784,763)
\path(9084,463)(8784,463)
\path(9684,1213)(9829,1213)
\path(9684,312)(9834,312)
\path(9684,912)(9830,912)
\path(11036,1362)(11186,1362)
\path(11035,161)(11185,161)
\path(11035,1062)(11185,1062)
\path(9829,1213)	(9902.182,1251.729)
	(9949.233,1287.295)
	(9973.667,1323.213)
	(9979.000,1363.000)

\path(9979,1363)	(9967.095,1410.125)
	(9940.893,1450.756)
	(9903.886,1485.425)
	(9859.566,1514.665)
	(9811.426,1539.010)
	(9762.956,1558.992)
	(9717.650,1575.144)
	(9679.000,1588.000)

\path(9679,1588)	(9636.563,1601.251)
	(9591.441,1613.646)
	(9543.812,1625.187)
	(9493.855,1635.873)
	(9441.750,1645.705)
	(9387.676,1654.681)
	(9331.811,1662.802)
	(9274.334,1670.069)
	(9215.426,1676.480)
	(9155.264,1682.037)
	(9094.028,1686.739)
	(9031.896,1690.586)
	(8969.049,1693.578)
	(8905.664,1695.715)
	(8841.922,1696.998)
	(8778.000,1697.425)
	(8714.078,1696.998)
	(8650.336,1695.715)
	(8586.951,1693.578)
	(8524.104,1690.586)
	(8461.972,1686.739)
	(8400.736,1682.037)
	(8340.574,1676.480)
	(8281.666,1670.069)
	(8224.189,1662.802)
	(8168.324,1654.681)
	(8114.250,1645.705)
	(8062.145,1635.873)
	(8012.188,1625.187)
	(7964.559,1613.646)
	(7919.437,1601.251)
	(7877.000,1588.000)

\path(7877,1588)	(7838.939,1575.350)
	(7795.064,1559.696)
	(7748.363,1540.331)
	(7701.823,1516.544)
	(7658.429,1487.626)
	(7621.169,1452.869)
	(7593.031,1411.563)
	(7577.000,1363.000)

\path(7577,1363)	(7577.333,1309.044)
	(7600.100,1257.074)
	(7648.817,1201.817)
	(7684.006,1171.308)
	(7727.000,1138.000)

\path(7612.768,1185.446)(7727.000,1138.000)(7648.535,1233.620)
\put(600,688){\makebox(0,0)[lb]{\smash{{\mathmode{D}}}}}
\put(2700,688){\makebox(0,0)[lb]{\smash{{\mathmode{P}}}}}
\put(2401,387){\makebox(0,0)[lb]{\smash{{\mathmode{\vdots}}}}}
\put(1125,462){\makebox(0,0)[lb]{\smash{{\mathmode{\vdots}}}}}
\put(1126,987){\makebox(0,0)[lb]{\smash{{\mathmode{\scriptstyle \ast\px}}}}}
\put(1500,688){\makebox(0,0)[lb]{\smash{{{\rm\scriptsize more}}}}}
\put(2250,1213){\makebox(0,0)[lb]{\smash{{\mathmode{\scriptstyle z}}}}}
\put(0,1138){\makebox(0,0)[lb]{\smash{{\mathmode{\scriptstyle z}}}}}
\put(0,838){\makebox(0,0)[lb]{\smash{{\mathmode{\scriptstyle z}}}}}
\put(0,538){\makebox(0,0)[lb]{\smash{{\mathmode{\scriptstyle w}}}}}
\put(4428,688){\makebox(0,0)[lb]{\smash{{\mathmode{D}}}}}
\put(6528,688){\makebox(0,0)[lb]{\smash{{\mathmode{P}}}}}
\put(6229,387){\makebox(0,0)[lb]{\smash{{\mathmode{\vdots}}}}}
\put(4953,462){\makebox(0,0)[lb]{\smash{{\mathmode{\vdots}}}}}
\put(4954,987){\makebox(0,0)[lb]{\smash{{\mathmode{\scriptstyle \ast\px}}}}}
\put(9230,688){\makebox(0,0)[lb]{\smash{{\mathmode{D}}}}}
\put(11330,688){\makebox(0,0)[lb]{\smash{{\mathmode{P}}}}}
\put(9759,462){\makebox(0,0)[lb]{\smash{{\mathmode{\vdots}}}}}
\put(11030,387){\makebox(0,0)[lb]{\smash{{\mathmode{\vdots}}}}}
\put(7877,913){\makebox(0,0)[lb]{\smash{{{\rm\scriptsize glue}}}}}
\put(7727,1213){\makebox(0,0)[lb]{\smash{{\mathmode{\scriptstyle \px}}}}}
\put(9976,1438){\makebox(0,0)[lb]{\smash{{\mathmode{\scriptstyle \ast}}}}}
\put(5326,688){\makebox(0,0)[lb]{\smash{{{\rm\scriptsize more}}}}}
\put(10126,688){\makebox(0,0)[lb]{\smash{{{\rm\scriptsize more}}}}}
\put(10876,1213){\makebox(0,0)[lb]{\smash{{\mathmode{\scriptstyle z}}}}}
\put(8626,1138){\makebox(0,0)[lb]{\smash{{\mathmode{\scriptstyle z}}}}}
\put(8626,838){\makebox(0,0)[lb]{\smash{{\mathmode{\scriptstyle z}}}}}
\put(8626,538){\makebox(0,0)[lb]{\smash{{\mathmode{\scriptstyle w}}}}}
\put(6076,1213){\makebox(0,0)[lb]{\smash{{\mathmode{\scriptstyle z}}}}}
\put(3826,1138){\makebox(0,0)[lb]{\smash{{\mathmode{\scriptstyle z}}}}}
\put(3826,838){\makebox(0,0)[lb]{\smash{{\mathmode{\scriptstyle z}}}}}
\put(3826,538){\makebox(0,0)[lb]{\smash{{\mathmode{\scriptstyle w}}}}}
\put(1350,388){\makebox(0,0)[lb]{\smash{{{\rm\scriptsize activity}}}}}
\put(5175,388){\makebox(0,0)[lb]{\smash{{{\rm\scriptsize activity}}}}}
\put(9975,388){\makebox(0,0)[lb]{\smash{{{\rm\scriptsize activity}}}}}
\end{picture}
}

%% file: Invariance.tex
\section{The Invariance Proof} \lbl{sec:Invariance}

Let us start with an easy warm-up:

\begin{proposition} \lbl{prop:Insensitive}
$\Aa_0$ is insensitive to orientation flips.
\end{proposition}

\begin{proof} Flipping the orientation of the component labeled $x$ in
some pure tangle $L$ acts on $\sigma\Zcheck(L)$ by flipping the sign
of all uni-trivalent diagrams that have an odd number of $x$-marked legs
(check~\cite[section~7.2]{Bar-Natan:Vassiliev} for the case of knots; the
case of pure tangles is the same). Namely, it acts by the substitution
$x\to-x$. Now use proposition~\ref{prop:Negation}.
\end{proof}

\subsection{$\Aa_0$ descends to regular links.} \lbl{subsec:DescendsToLinks}
Our first real task is to show that if two regular pure tangles have the
same closures then they have the same pre-normalized \Arhus{} integral
and hence the pre-normalized \Arhus{} integral $\Aa_0$ descends to
regular links. We first extend the definition of $\Aa_0$ to some
larger class of ``closable'' objects (definition~\ref{def:DML}),
the class of regular dotted Morse links. We then show that $\Aa_0$
descends from that class to links (proposition~\ref{prop:MoveDots}),
and finally that regular pure tangles ``embed'' in regular dotted Morse
links (proposition~\ref{prop:RPTisRDML}). Taken together, these two
propositions imply that $\Aa_0$ descends to regular links also from
regular pure tangles.

We should note that the \Arhus{} integral can be defined and all of its
properties can be proven fully within the class of regular dotted Morse
links, and that this is essentially what we do in this paper. The only
reasons we also work with regular pure tangles are reasons of elegance.

\begin{definition} \lbl{def:DML}
A dotted Morse link $L$ is a link embedded in $\RR^3_{xyt}$ so that the
third Euclidean coordinate $t$ is a Morse function on it, together with
a dot marked on each component.  We assume that the components of $L$
are labeled by the elements of some label set $X$.  Notice that we do
not divide by isotopies. The ``closure'' of a dotted Morse link is
the ($X$-marked) link obtained by forgetting the dots and dividing
by isotopies. These definitions have obvious framed counterparts.
\end{definition}

\begin{remark} Why so ugly a definition? Because all other choices are
even worse. We have to ``dot'' the link components because we want
the Kontsevich integral to be valued in $\calA(\uparrow_X)$
(see below). But then we have to give up isotopy invariance at the time
slices of the dots, and it is simpler to give it up altogether. See
also the comment about q-tangles/non-associative tangles above
proposition~\ref{prop:RPTisRDML}.
\end{remark}

The framed Kontsevich integral $Z$, as well as the variant $\Zcheck$
due to Le, H.~Murakami, J.~Murakami, and
Ohtsuki~\cite{LeMurakami2Ohtsuki:3Manifold}, both have obvious
definitions in the case of framed dotted Morse links. The new bit is
that each component has dot marked on it, which can serve as a cutting
mark for scissors. In other words, every component can be regarded as a
directed line, and thus the images of $Z$ and $\Zcheck$ are in
$\calA(\uparrow_X)$. But now we can compose $\Zcheck$ with
$\sigma$ and then with $\Gint$, and we find that the pre-normalized \Arhus{}
integral $\Aa_0$ can also be defined on regular dotted Morse links
(framed dotted Morse links with a non-singular linking matrix).

\begin{proposition} \lbl{prop:MoveDots}
$\Aa_0$ descends from regular dotted Morse links to regular links.
\end{proposition}

\begin{proof} The usual invariance argument for the Kontsevich integral
(see~\cite{Kontsevich:Vassiliev, Bar-Natan:Vassiliev}) applies also in
the case of (framed) dotted Morse links, provided the time slices of
the dots are frozen. So the only thing we need to prove is that $\Aa_0$
is invariant under sliding the dots along a component; once this is
done, the frozen time slices melt and we have complete invariance.

A different way of saying that a dot moves on a framed dotted Morse
link $L$ is saying that we have two dots on one component (say $z$),
cutting it into two subcomponents $x$ and $y$. Each time we ignore one
of the dots and compute $\Zcheck$, getting two results $G_1$ and $G_2$,
and we wish to compare the integrals of $G_1$ and $G_2$. Alternatively,
we can keep both dots on the $z$ component and compute $\Zcheck$ in the
usual way, only cutting the resulting chord diagrams open at both dots,
getting a result $G$ in the space\footnote{
  The notation means: pure tangle diagrams whose skeleton components
  are labeled by the symbols $x$, $y$, and some additional $n-1$ symbols
  in some set $E$ of ``Extra variables''. Below we will use variations of
  this notation with no further comment.
}
$\calA(\uparrow_x\uparrow_y\uparrow_E)$. From $G$ both $G_1$ and
$G_2$ can be recovered by attaching the components $x$ and $y$
in either of the two possible orders. This process is made precise
in definition~\ref{def:mxyz} below, and the fact that $\Gint\sigma
G_1=\Gint\sigma G_2$ follows from the ``cyclic invariance lemma''
(lemma~\ref{lemma:CyclicSymmetry}) below. We only need to comment that
$G$ is group-like like any evaluation of the Kontsevich integral.
\end{proof}

\begin{definition} \lbl{def:mxyz}
Let
$
  \mright^{xy}_z:\calA(\uparrow_x\uparrow_y\uparrow_E)
  \to\calA(\uparrow_z\uparrow_E)
$
be the map described in figure~\ref{fig:mxyz}. The map $\mright^{yx}_z$
is the same, only with the roles of $x$ and $y$ interchanged.
\end{definition}

\begin{figure}[htpb]
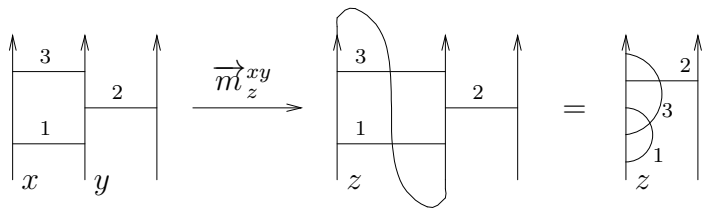

\[ \printname{mxyz}
	\setlength{\unitlength}{0.75\standardunitlength}
	\begin{array}{c}  \hspace{-1.7mm}
        	\raisebox{-8pt}{\input draws/mxyz.tex }
        	\hspace{-1.9mm}
	\end{array}
 \]
\caption{
  The map $\protect\mright^{xy}_z$ in the case
  $n=2$: Connect the strands labeled $x$ and $y$ in a diagram
  in $\calA(\uparrow_x\uparrow_y\uparrow)$, to form a new ``long''
  strand labeled $z$, without touching all extra strands.
}
\lbl{fig:mxyz}
\end{figure}

\begin{lemma} \lbl{lemma:CyclicSymmetry} (the cyclic invariance lemma).
If $G\in\calA(\uparrow_x\uparrow_y\uparrow_E)$ is group-like and
$\sigma \mright^{xy}_z G$ is an integrable member of $\calB_n$, then $\sigma
\mright^{yx}_z G$ is also integrable and the two integrals are equal:
\[ \Gint\sigma \mright^{xy}_z G\,dzdE = \Gint\sigma \mright^{yx}_z  G\,dzdE. \]
\end{lemma}

\begin{proof} It is easy to verify that $G_1=\sigma\mright^{xy}_z G$ and
$G_2=\sigma\mright^{yx}_z G$ are both group-like and hence Gaussian, and that
they have the same covariance matrix (when we apply this lemma as in
proposition~\ref{prop:MoveDots}, in both cases the covariance matrix is
the linking matrix of the underlying link). Thus if one is integrable
so is the other, and we have to prove the equality of the integrals.

{\em The case of knots: } If $n=1$ then the fact that $\calA(\uparrow)$
is isomorphic to $\calA(\circlearrowleft)$ (namely, the commutativity
of $\calA(\uparrow)$, see~\cite{Bar-Natan:Vassiliev}) implies that
$\mright^{xy}_z=\mright^{yx}_z$ and there's nothing to prove.

{\em The lucky case: } If $G_{1,2}$ are integrable with respect to
$E$ we can use proposition~\ref{prop:IteratedIntegral}
and compute the integrals with respect to those variables first. The
results are diagrams labeled by just one variable ($z$) (namely,
functions of just one variable), and we are back in the previous case.

{\em The ugly case: } If $G_{1,2}$ are not integrable with respect to
$E$, we can perturb them a bit by multiplying by some
$\exp\sum_{i,j}\epsilon_{ij}\strutn{e_i}{e_j}$ to get
$G^\epsilon_{1,2}$. The integrals of $G^\epsilon_{1,2}$ (with respect
to all variables) depend polynomially on the $\epsilon$'s in any given
degree. For generic $\epsilon$'s we get $G^\epsilon_{1,2}$ that are
integrable with respect to $E$, and we fall back to the
lucky case. Thus the integrals of $G^\epsilon_{1,2}$ are equal as power
series in the $\epsilon$'s, and in particular they are equal at
$\epsilon_{ij=0}$.
\end{proof}

Every (framed) pure tangle $L$ defines a class of associated (framed)
dotted Morse links, obtained by picking a specific Morse representative
of $L$, marking dots at the tops of all strands, and closing to a link
in some specific way making sure that the down-going strands used in
the closure are very far ($d$ miles away) from the original pure tangle.
An example is in figure~\ref{fig:DottedMorseLinks}.  What's very far? In
the infinite limit; meaning that whenever we refer to an associated
(framed) dotted Morse link, we really mean ``a sequence of such,
with $d\to\infty$''. To remind ourselves of that, we add the phrase
``(at limit)'' to the statements that are true only when this (or a
similar, see below) limit is taken. If one is ready to sacrifice some
simplicity, all of these statements can be formulated without limits
if the technology of q-tangles (\cite{LeMurakami:Universal}) (or, what
is nearly the same, non associative tangles (\cite{Bar-Natan:NAT}))
is used instead of using specific Morse embeddings. Readers familiar
with~\cite{LeMurakami:Universal} and/or \cite{Bar-Natan:NAT} should have
no difficulty translating our language to the more precise language of
those papers.

\begin{figure}[htpb]
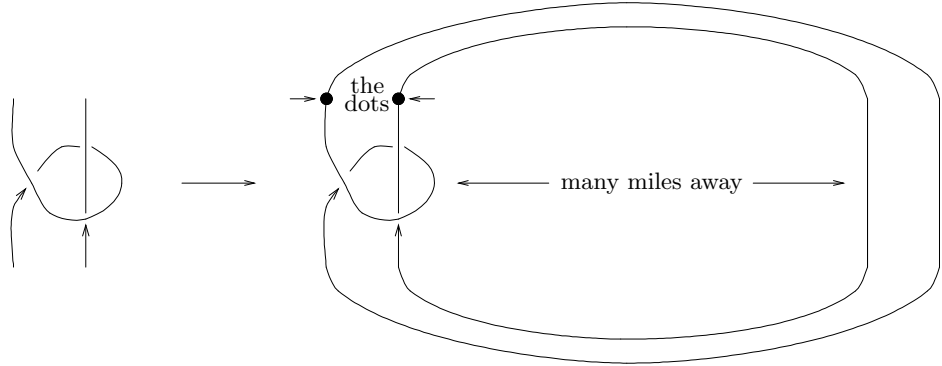

\[ \printname{DottedMorseLinks}
	\setlength{\unitlength}{0.5\standardunitlength}
	\begin{array}{c}  \hspace{-1.7mm}
        	\raisebox{-8pt}{\input draws/DottedMorseLinks.tex }
        	\hspace{-1.9mm}
	\end{array}
 \]
\caption{A pure tangle and an associated dotted Morse link.}
\lbl{fig:DottedMorseLinks}
\end{figure}

\begin{proposition} \lbl{prop:RPTisRDML} (at limit)
If $L$ is a pure tangle and $L^\bullet$ is an associated dotted Morse
link, then $\Zcheck(L)=\Zcheck(L^\bullet)$.
\end{proposition}

\par\noindent
\parbox{5in}{
  \begin{proof} (at limit)
  The dotted Morse link $L^\bullet$ is obtained from $L$ by sticking
  $L$ within a ``closure element'' $C_X$, shown on the right (for
  $|X|=3$).  Let $C'_X$ be $C_X$ with the two boxes at its ends removed.
  These two boxes denote ``adapters'' $A$ and $A^{-1}$ that only change
  the strand spacings to be uniform, from a possibly non-uniform
  spacings in $C'_X$.
  \def\qed{}
  \end{proof}
}$\quad\printname{Closer}
	\setlength{\unitlength}{0.5\standardunitlength}
	\begin{array}{c}  \hspace{-1.7mm}
        	\raisebox{-8pt}{\input draws/Closer.tex }
        	\hspace{-1.9mm}
	\end{array}
$

Inspecting the definitions of $\Zcheck$ for pure tangles
(see~\cite[definition~\ref{I-def:Zcheck}]{Aarhus:I}) and for dotted Morse
links (see~\cite{LeMurakami2Ohtsuki:3Manifold}), we see that we only
need to show that $Z(C_X)=\Delta_X\nu$ in the space $\calA(\uparrow_X)$
(check~\cite[definition~\ref{I-def:Zcheck}]{Aarhus:I} for the definition
of $\Delta_X$). Here $C_X$ is itself regarded as a dotted Morse link
(with the dots at the space allotted for $L$, which is assumed to
be small relative to the size of $C_X$ itself) and $Z$ denotes
the Kontsevich integral in its standard normalization. Clearly
$Z(C_{\{x\}})=\nu=\Delta_{\{x\}}\nu$, as $C_{\{x\}}$ is the dotted
unknot and $\nu$ is by definition the Kontsevich integral of
the unknot.  Theorem~4.1 of~\cite{LeMurakami:Parallel}, rephrased
for dotted Morse links, says that doubling a component (so that
the two daughter components are parallel and very close) and then
computing $Z$ is equal to $Z$ followed by $\Delta$. In other words,
$Z(C_{\{x,y\}})=\Delta_{\{x,y\}}(\nu)$. Iterating this argument, we find
that $Z(C'_X)=\Delta_X\nu$, for some {\em specific} (at limit) choice of
strand spacings in $C'_X$. But $\Delta_X\nu$ is central and hence $Z(C_X)
= Z(A^{-1})Z(C'_X)Z(A) = Z(A)^{-1}(\Delta_X\nu)Z(A) = \Delta_X\nu$. \qed

\subsection{$\Aa_0$ is invariant under the second Kirby move.}
\lbl{subsec:KirbyII}

\begin{definition}
A tight Kirby move $L_1\to L_2$ is a move between two framed dotted Morse
links $L_1$ and $L_2$ as in figure~\ref{fig:DetailedKirbyII}, in which
\par\noindent
\parbox{5in}{
  \begin{itemize}
  \item Before the move the two parallel strands in the domain $S$ are
    ``tight''. Namely, they are very close to each other relative to the
    distance between them and any other feature of the link.
  \item The doubling of the $y$ component is done in a ``very tight''
    fashion.  Namely, the distance between the the copies of $y$ produced
    is very small relative to the scale in which the rest of the link
    is drawn, even much smaller than the original distance between the
    $x$ and $y$ components.
  \item The dots on the $x$ and $y$ components are inside the domain $S$ both
    before and after surgery, and they are placed as in the picture on the
    right.
  \end{itemize}
}\qquad$\printname{DotsInS}
	\setlength{\unitlength}{0.75\standardunitlength}
	\begin{array}{c}  \hspace{-1.7mm}
        	\raisebox{-8pt}{\input draws/DotsInS.tex }
        	\hspace{-1.9mm}
	\end{array}
$
\par\noindent
We extend the notion of ``at limit'' to mean that ``tightness'' is also
increased ad infinitum.
\end{definition}

\begin{figure}[htpb]
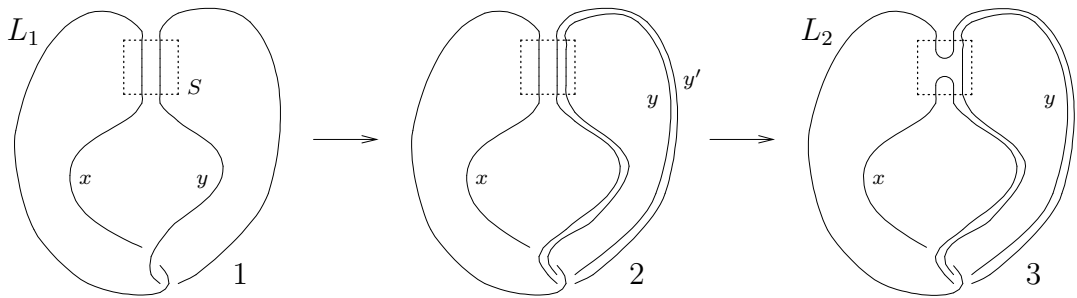

\[ \printname{DetailedKirbyII}
	\setlength{\unitlength}{0.75\standardunitlength}
	\begin{array}{c}  \hspace{-1.7mm}
        	\raisebox{-8pt}{\input draws/DetailedKirbyII.tex }
        	\hspace{-1.9mm}
	\end{array}
 \]
\caption{
  The second Kirby move $L_1\to L_2$: (1) Some domain $S$ in space
  in which some two components of the link (denoted $x$ and $y$)
  are adjacent and nearly parallel is specified.  (2) The component
  $y$ is doubled (using its framing), getting a new component $y'$.
  (3) A surgery is performed in $S$ combining $y'$ into $x$, so that
  now the $x$ component runs parallel to the $y$ component in addition
  to running its own course. We say that the component $x$ ``slides''
  over the component $y$.
}
\lbl{fig:DetailedKirbyII}
\end{figure}

The following proposition is due to Le, H.~Murakami, J.~Murakami, and
Ohtsuki~\cite{LeMurakami2Ohtsuki:3Manifold}. It holds for $\Zcheck$ and
not for $Z$, and it is the reason why
\cite{LeMurakami2Ohtsuki:3Manifold} introduced $\Zcheck$. We present the
``at limit'' version, which is equivalent to the ``q-tangle'' version
proven in~\cite{LeMurakami2Ohtsuki:3Manifold}.

\begin{proposition} \lbl{prop:LMMOResult}
(at limit, proof in~\cite{LeMurakami2Ohtsuki:3Manifold}) Let $L_1\to
L_2$ be a tight Kirby move between two framed dotted Morse links $L_1$
and $L_2$ marked as in figure~\ref{fig:DetailedKirbyII}. Then
\[ \Zcheck(L_2) = \Upsright\Zcheck(L_1), \]
where $\Upsright=\mright^{xy'}_x\circ\Delta^y_{yy'}$ and $\Delta^y_{yy'}$
denotes the diagram-level operation of doubling the $y$ strand (lifting
all vertices on it in all possible ways, and calling the double $y'$).
(Compare with~\cite[equation~\eqref{I-Uglevel}]{Aarhus:I}).
\end{proposition}

\begin{proof}[Proof of proposition~\ref{prop:KirbyII}]
After propositions \ref{prop:Insensitive}, \ref{prop:MoveDots}
and~\ref{prop:RPTisRDML} have been proven, all that remains is to
show that $\Aa_0$ is invariant under tight Kirby moves of framed
dotted Morse links. (Notice that every Kirby move between links has a
presentation as a tight Kirby move between dotted Morse links).  Using
proposition~\ref{prop:LMMOResult} we find that it is enough to show that
whenever $G$ is a non-degenerate Gaussian (think $G=\sigma\Zcheck(L_1)$),
\begin{equation} \lbl{eq:Upsilon}
  \Gint G \,dE = \Gint\Upsright G\,dE,
\end{equation}
where we re-use the symbol $\Upsright$ to denote the same operation
on the level of uni-trivalent diagrams.

Let $\Upsilon$ be the same as $\Upsright$, only with $m^{xy'}_x$
replacing $\mright^{xy'}_x$, where $m^{xy'}_xG:=G/(y'\to x)$. The
operation $\Upsilon$ is a substitution operation of the form discussed
in section~\ref{sec:calculus}; $\Upsilon G= G/(y\to x+y)$.
Exercise~\ref{ex:ReparametrizeGint} shows that
equation~\eqref{eq:Upsilon} holds if $\Upsright$ is replaced by
$\Upsilon$. So we only need to analyze the difference
$\Upsright-\Upsilon$. The difference $\mright^{xy'}_x-m^{xy'}_x$ is
given by gluing a certain sum $D'$ of forests whose roots are labeled
$x$ and whose leaves are labeled $\partial_x$ and $\partial_{y'}$,
followed by the substitution $(y'\to x)$. Hence,
\[
  (\Upsright-\Upsilon)G
  = \left(D'\apply[G/(y\to y+y')]\right)/(y'\to x)
  = (D\apply G)/(y\to x+y),
\]
where $D$ is $D'$ with every $\partial_{y'}$ replaced by a
$\partial_y$. (A precise formula for $D$ can be derived from the
results of section~\ref{subsec:BCH}, but we don't need it here).
Clearly, $\operatorname{div}_{y}D=0$; the coefficients of $D$ are
independent of $y$ and every term in $D$ is of positive degree in
$\partial_y$. Now
\begin{alignat*}{2}
  \Gint(\Upsright-\Upsilon)G\,dE
    & = \Gint (D\apply G)/(y\to x+y)\,dE & \\
    & = \Gint D\apply G\,dE
      & \text{by exercise~\ref{ex:ReparametrizeGint}} \\
    & = 0
      & \text{by proposition~\ref{prop:Divergence}.}
\end{alignat*}
\end{proof}

\subsection{$\Aa$ and invariance under the first Kirby move.}
\lbl{subsec:FirstKirby}
\begin{proof}[Proof of theorem~\ref{thm:Invariance}]
Flipping the orientation of a component negates all linking numbers
between it and any other component, and hence the linking matrix changes
by a similarity transformation. The second Kirby moves adds all linking
numbers involving the $y$-component (see figure~\ref{fig:DetailedKirbyII})
to the corresponding ones with the $x$-component. This again is a
similarity transformation. Similarity transformations do not change the
numbers $\sigma_\pm$ of positive/negative eigenvalues. Thus $\Aa$ is
invariant under orientation flips and under the second Kirby move.

All that's left is to show that $\Aa$ is invariant under the first Kirby
move. Namely, that it is invariant under taking the disjoint union of a
link with $U_\pm$, the unknot with framing $\pm 1$.

Let $L$ be an $n$-component regular link. Adding a far-away $U_+$
component to $L$ multiplies $\sigma\Zcheck$ by $\sigma\Zcheck(U_+)$
(using the disjoint union product). The new linking matrix is block
diagonal, with an additional $+1$ entry on the diagonal, and the same
holds for the new inverse linking matrix. Thus the $(n+1)$-variable
Gaussian integral of $\sigma\Zcheck(L\mathbin\udot U_+)$ factors as the
$n$-variable integral of $\sigma\Zcheck(L)$ times the $1$-variable
integral of $\sigma\Zcheck(U_+)$. We find that $\Aa_0(L\mathbin\udot
U_+)=\Aa_0(L)\mathbin\udot\Aa_0(U_+)$, and as $\sigma_+$ also increases
by $1$, $\Aa(L\mathbin\udot U_+)=\Aa(L)$ as required. A similar argument
works in the case of $U_-$.
\end{proof}

%% file: draws/mxyz.tex
\begingroup\makeatletter\ifx\SetFigFont\undefined
\def\x#1#2#3#4#5#6#7\relax{\def\x{#1#2#3#4#5#6}}%
\expandafter\x\fmtname xxxxxx\relax \def\y{splain}%
\ifx\x\y   
\gdef\SetFigFont#1#2#3{%
  \ifnum #1<17\tiny\else \ifnum #1<20\small\else
  \ifnum #1<24\normalsize\else \ifnum #1<29\large\else
  \ifnum #1<34\Large\else \ifnum #1<41\LARGE\else
     \huge\fi\fi\fi\fi\fi\fi
  \csname #3\endcsname}%
\else
\gdef\SetFigFont#1#2#3{\begingroup
  \count@#1\relax \ifnum 25<\count@\count@25\fi
  \def\x{\endgroup\@setsize\SetFigFont{#2pt}}%
  \expandafter\x
    \csname \romannumeral\the\count@ pt\expandafter\endcsname
    \csname @\romannumeral\the\count@ pt\endcsname
  \csname #3\endcsname}%
\fi
\fi\endgroup
{\renewcommand{\dashlinestretch}{30}
\begin{picture}(5746,1703)(0,-10)
\put(5112.000,619.000){\arc{450.000}{4.7124}{7.8540}}
\put(5074.500,956.500){\arc{679.154}{4.8230}{7.7433}}
\path(12,244)(12,1444)
\path(42.000,1324.000)(12.000,1444.000)(-18.000,1324.000)
\path(612,244)(612,1444)
\path(642.000,1324.000)(612.000,1444.000)(582.000,1324.000)
\path(1212,244)(1212,1444)
\path(1242.000,1324.000)(1212.000,1444.000)(1182.000,1324.000)
\path(12,1144)(612,1144)
\path(612,844)(1212,844)
\path(12,544)(612,544)
\path(1512,844)(2412,844)
\path(2292.000,814.000)(2412.000,844.000)(2292.000,874.000)
\path(2712,244)(2712,1444)
\path(2742.000,1324.000)(2712.000,1444.000)(2682.000,1324.000)
\path(3612,244)(3612,1444)
\path(3642.000,1324.000)(3612.000,1444.000)(3582.000,1324.000)
\path(4212,244)(4212,1444)
\path(4242.000,1324.000)(4212.000,1444.000)(4182.000,1324.000)
\path(5112,244)(5112,1444)
\path(5142.000,1324.000)(5112.000,1444.000)(5082.000,1324.000)
\path(5712,244)(5712,1444)
\path(5742.000,1324.000)(5712.000,1444.000)(5682.000,1324.000)
\path(3612,844)(4212,844)
\path(2712,1144)(3612,1144)
\path(2712,544)(3612,544)
\path(5112,1069)(5712,1069)
\path(2712,1444)	(2712.000,1492.986)
	(2712.000,1519.000)

\path(2712,1519)	(2707.418,1555.495)
	(2712.000,1594.000)

\path(2712,1594)	(2775.870,1656.602)
	(2818.355,1673.085)
	(2862.000,1669.000)

\path(2862,1669)	(2914.539,1644.018)
	(2960.539,1613.891)
	(3034.662,1539.576)
	(3063.658,1496.073)
	(3087.858,1448.798)
	(3107.698,1398.092)
	(3123.614,1344.301)
	(3136.042,1287.765)
	(3145.418,1228.829)
	(3152.179,1167.834)
	(3156.759,1105.125)
	(3159.595,1041.044)
	(3161.123,975.935)
	(3161.780,910.139)
	(3162.000,844.000)
	(3162.220,777.861)
	(3162.877,712.065)
	(3164.405,646.956)
	(3167.241,582.875)
	(3171.821,520.166)
	(3178.582,459.171)
	(3187.958,400.235)
	(3200.386,343.699)
	(3216.302,289.908)
	(3236.142,239.202)
	(3260.342,191.927)
	(3289.338,148.424)
	(3363.461,74.109)
	(3409.461,43.982)
	(3462.000,19.000)

\path(3462,19)	(3505.645,14.915)
	(3548.130,31.398)
	(3612.000,94.000)

\path(3612,94)	(3616.582,132.505)
	(3612.000,169.000)

\path(3612,169)	(3612.000,195.014)
	(3612.000,244.000)

\put(87,169){\makebox(0,0)[lb]{\smash{{\mathmode{x}}}}}
\put(687,169){\makebox(0,0)[lb]{\smash{{\mathmode{y}}}}}
\put(837,919){\makebox(0,0)[lb]{\smash{{\mathmode{\scriptstyle 2}}}}}
\put(237,1219){\makebox(0,0)[lb]{\smash{{\mathmode{\scriptstyle 3}}}}}
\put(237,619){\makebox(0,0)[lb]{\smash{{\mathmode{\scriptstyle 1}}}}}
\put(4587,769){\makebox(0,0)[lb]{\smash{{\mathmode{=}}}}}
\put(3837,919){\makebox(0,0)[lb]{\smash{{\mathmode{\scriptstyle 2}}}}}
\put(2862,1219){\makebox(0,0)[lb]{\smash{{\mathmode{\scriptstyle 3}}}}}
\put(2862,619){\makebox(0,0)[lb]{\smash{{\mathmode{\scriptstyle 1}}}}}
\put(5412,769){\makebox(0,0)[lb]{\smash{{\mathmode{\scriptstyle 3}}}}}
\put(5337,394){\makebox(0,0)[lb]{\smash{{\mathmode{\scriptstyle 1}}}}}
\put(5562,1144){\makebox(0,0)[lb]{\smash{{\mathmode{\scriptstyle 2}}}}}
\put(2787,169){\makebox(0,0)[lb]{\smash{{\mathmode{z}}}}}
\put(1662,994){\makebox(0,0)[lb]{\smash{{\mathmode{\mright^{xy}_z}}}}}
\put(5187,169){\makebox(0,0)[lb]{\smash{{\mathmode{z}}}}}
\end{picture}
}

%% file: draws/DottedMorseLinks.tex
\begingroup\makeatletter\ifx\SetFigFont\undefined
\def\x#1#2#3#4#5#6#7\relax{\def\x{#1#2#3#4#5#6}}%
\expandafter\x\fmtname xxxxxx\relax \def\y{splain}%
\ifx\x\y   
\gdef\SetFigFont#1#2#3{%
  \ifnum #1<17\tiny\else \ifnum #1<20\small\else
  \ifnum #1<24\normalsize\else \ifnum #1<29\large\else
  \ifnum #1<34\Large\else \ifnum #1<41\LARGE\else
     \huge\fi\fi\fi\fi\fi\fi
  \csname #3\endcsname}%
\else
\gdef\SetFigFont#1#2#3{\begingroup
  \count@#1\relax \ifnum 25<\count@\count@25\fi
  \def\x{\endgroup\@setsize\SetFigFont{#2pt}}%
  \expandafter\x
    \csname \romannumeral\the\count@ pt\expandafter\endcsname
    \csname @\romannumeral\the\count@ pt\endcsname
  \csname #3\endcsname}%
\fi
\fi\endgroup
{\renewcommand{\dashlinestretch}{30}
\begin{picture}(11601,4535)(0,-10)
\put(3939,3310){\blacken\ellipse{150}{150}}
\put(3939,3310){\ellipse{150}{150}}
\put(4839,3310){\blacken\ellipse{150}{150}}
\put(4839,3310){\ellipse{150}{150}}
\path(10689,3310)(10689,1210)
\path(10689,3310)(10689,1210)
\path(11589,3310)(11589,1210)
\path(11589,3310)(11589,1210)
\path(6714,2260)(5589,2260)
\path(5709.000,2290.000)(5589.000,2260.000)(5709.000,2230.000)
\path(5289,3310)(4989,3310)
\path(5109.000,3340.000)(4989.000,3310.000)(5109.000,3280.000)
\path(3489,3310)(3789,3310)
\path(3669.000,3280.000)(3789.000,3310.000)(3669.000,3340.000)
\path(2139,2260)(3039,2260)
\path(2919.000,2230.000)(3039.000,2260.000)(2919.000,2290.000)
\path(939,1210)(939,1735)
\path(976.500,1600.000)(939.000,1735.000)(901.500,1600.000)
\path(939,1885)(939,3310)
\path(4839,1210)(4839,1735)
\path(4876.500,1600.000)(4839.000,1735.000)(4801.500,1600.000)
\path(4839,1885)(4839,3310)
\path(9264,2260)(10389,2260)
\path(10269.000,2230.000)(10389.000,2260.000)(10269.000,2290.000)
\path(4839,3310)	(4859.449,3373.207)
	(4878.533,3427.311)
	(4896.693,3473.192)
	(4914.368,3511.727)
	(4950.018,3570.280)
	(4989.000,3610.000)

\path(4989,3610)	(5048.794,3653.822)
	(5113.045,3695.575)
	(5181.464,3735.303)
	(5253.765,3773.051)
	(5291.282,3791.197)
	(5329.662,3808.865)
	(5368.868,3826.060)
	(5408.866,3842.789)
	(5449.619,3859.058)
	(5491.091,3874.871)
	(5533.247,3890.234)
	(5576.050,3905.153)
	(5619.465,3919.634)
	(5663.456,3933.682)
	(5707.987,3947.304)
	(5753.022,3960.504)
	(5798.526,3973.288)
	(5844.461,3985.662)
	(5890.793,3997.632)
	(5937.486,4009.203)
	(5984.503,4020.381)
	(6031.809,4031.171)
	(6079.368,4041.580)
	(6127.145,4051.613)
	(6175.102,4061.275)
	(6223.205,4070.572)
	(6271.417,4079.510)
	(6319.703,4088.095)
	(6368.026,4096.332)
	(6416.351,4104.227)
	(6464.642,4111.785)
	(6512.864,4119.012)
	(6560.979,4125.914)
	(6608.953,4132.496)
	(6656.749,4138.765)
	(6704.331,4144.725)
	(6751.665,4150.383)
	(6798.713,4155.744)
	(6845.440,4160.813)
	(6891.810,4165.597)
	(6937.787,4170.101)
	(6983.336,4174.331)
	(7028.420,4178.292)
	(7073.003,4181.989)
	(7117.051,4185.430)
	(7160.526,4188.619)
	(7203.392,4191.561)
	(7245.615,4194.264)
	(7287.158,4196.731)
	(7327.985,4198.970)
	(7368.061,4200.985)
	(7407.349,4202.782)
	(7445.814,4204.368)
	(7483.420,4205.746)
	(7555.910,4207.907)
	(7624.532,4209.310)
	(7689.000,4210.000)

\path(7689,4210)	(7756.889,4210.050)
	(7829.159,4209.429)
	(7866.843,4208.851)
	(7905.509,4208.085)
	(7945.120,4207.124)
	(7985.637,4205.962)
	(8027.024,4204.591)
	(8069.242,4203.006)
	(8112.254,4201.198)
	(8156.022,4199.162)
	(8200.508,4196.891)
	(8245.675,4194.378)
	(8291.485,4191.615)
	(8337.901,4188.597)
	(8384.884,4185.317)
	(8432.397,4181.767)
	(8480.403,4177.941)
	(8528.863,4173.832)
	(8577.739,4169.434)
	(8626.995,4164.739)
	(8676.593,4159.741)
	(8726.494,4154.433)
	(8776.662,4148.808)
	(8827.058,4142.860)
	(8877.644,4136.581)
	(8928.384,4129.965)
	(8979.239,4123.005)
	(9030.172,4115.694)
	(9081.145,4108.026)
	(9132.120,4099.994)
	(9183.060,4091.590)
	(9233.926,4082.809)
	(9284.682,4073.643)
	(9335.290,4064.085)
	(9385.711,4054.130)
	(9435.908,4043.769)
	(9485.844,4032.996)
	(9535.481,4021.805)
	(9584.780,4010.189)
	(9633.706,3998.140)
	(9682.218,3985.653)
	(9730.281,3972.719)
	(9777.856,3959.334)
	(9824.906,3945.489)
	(9871.393,3931.177)
	(9917.279,3916.393)
	(9962.527,3901.130)
	(10007.098,3885.379)
	(10050.956,3869.136)
	(10094.062,3852.393)
	(10136.379,3835.142)
	(10177.869,3817.379)
	(10218.495,3799.094)
	(10258.218,3780.283)
	(10297.001,3760.938)
	(10334.807,3741.052)
	(10407.335,3699.630)
	(10475.500,3655.964)
	(10539.000,3610.000)

\path(10539,3610)	(10578.263,3570.342)
	(10614.007,3511.810)
	(10631.658,3473.269)
	(10649.748,3427.373)
	(10668.715,3373.243)
	(10689.000,3310.000)

\path(11589,1210)	(11568.302,1146.846)
	(11549.039,1092.779)
	(11530.772,1046.921)
	(11513.062,1008.393)
	(11477.555,949.810)
	(11439.000,910.000)

\path(11439,910)	(11398.664,878.921)
	(11356.729,848.475)
	(11313.247,818.658)
	(11268.265,789.464)
	(11221.834,760.890)
	(11174.002,732.931)
	(11124.818,705.583)
	(11074.332,678.840)
	(11022.594,652.699)
	(10969.652,627.155)
	(10915.556,602.204)
	(10860.355,577.841)
	(10804.098,554.061)
	(10746.834,530.861)
	(10688.613,508.235)
	(10629.484,486.179)
	(10569.497,464.689)
	(10508.700,443.760)
	(10447.142,423.387)
	(10384.874,403.567)
	(10321.944,384.294)
	(10258.402,365.565)
	(10194.296,347.374)
	(10129.677,329.718)
	(10064.593,312.591)
	(9999.093,295.990)
	(9933.227,279.909)
	(9867.045,264.345)
	(9800.595,249.292)
	(9733.926,234.747)
	(9667.089,220.704)
	(9600.131,207.160)
	(9533.103,194.110)
	(9466.054,181.549)
	(9399.033,169.472)
	(9332.089,157.877)
	(9265.271,146.757)
	(9198.629,136.108)
	(9132.212,125.926)
	(9066.069,116.207)
	(9000.250,106.946)
	(8934.804,98.138)
	(8869.780,89.779)
	(8805.227,81.864)
	(8741.194,74.390)
	(8677.731,67.351)
	(8614.888,60.743)
	(8552.713,54.561)
	(8491.255,48.802)
	(8430.564,43.460)
	(8370.690,38.531)
	(8311.680,34.010)
	(8253.586,29.894)
	(8196.455,26.177)
	(8140.337,22.856)
	(8085.282,19.925)
	(8031.339,17.380)
	(7978.556,15.217)
	(7926.984,13.431)
	(7876.671,12.017)
	(7827.667,10.972)
	(7780.021,10.290)
	(7733.782,9.968)
	(7689.000,10.000)

\path(7689,10)	(7645.911,10.381)
	(7601.424,11.125)
	(7555.587,12.235)
	(7508.445,13.716)
	(7460.048,15.570)
	(7410.442,17.802)
	(7359.675,20.414)
	(7307.794,23.410)
	(7254.847,26.793)
	(7200.881,30.567)
	(7145.944,34.735)
	(7090.083,39.301)
	(7033.345,44.268)
	(6975.778,49.639)
	(6917.430,55.418)
	(6858.348,61.609)
	(6798.579,68.215)
	(6738.171,75.239)
	(6677.171,82.685)
	(6615.626,90.556)
	(6553.585,98.855)
	(6491.094,107.587)
	(6428.201,116.754)
	(6364.954,126.360)
	(6301.400,136.409)
	(6237.586,146.903)
	(6173.559,157.846)
	(6109.368,169.242)
	(6045.060,181.095)
	(5980.681,193.406)
	(5916.281,206.181)
	(5851.905,219.423)
	(5787.602,233.134)
	(5723.418,247.318)
	(5659.402,261.979)
	(5595.601,277.121)
	(5532.062,292.746)
	(5468.833,308.858)
	(5405.960,325.460)
	(5343.493,342.557)
	(5281.477,360.151)
	(5219.961,378.246)
	(5158.992,396.845)
	(5098.617,415.952)
	(5038.884,435.569)
	(4979.840,455.702)
	(4921.533,476.352)
	(4864.010,497.524)
	(4807.318,519.221)
	(4751.506,541.446)
	(4696.620,564.203)
	(4642.708,587.496)
	(4589.817,611.326)
	(4537.995,635.699)
	(4487.290,660.618)
	(4437.748,686.085)
	(4389.417,712.105)
	(4342.345,738.681)
	(4296.579,765.816)
	(4252.167,793.513)
	(4209.155,821.777)
	(4167.592,850.611)
	(4127.524,880.017)
	(4089.000,910.000)

\path(4089,910)	(4050.682,949.861)
	(4015.253,1008.460)
	(3997.523,1046.984)
	(3979.197,1092.829)
	(3959.836,1146.875)
	(3939.000,1210.000)

\path(11589,3310)	(11568.302,3373.154)
	(11549.039,3427.221)
	(11530.772,3473.079)
	(11513.062,3511.608)
	(11477.555,3570.190)
	(11439.000,3610.000)

\path(11439,3610)	(11398.664,3641.079)
	(11356.729,3671.525)
	(11313.247,3701.342)
	(11268.265,3730.536)
	(11221.834,3759.110)
	(11174.002,3787.069)
	(11124.818,3814.417)
	(11074.332,3841.160)
	(11022.594,3867.301)
	(10969.652,3892.845)
	(10915.556,3917.796)
	(10860.355,3942.159)
	(10804.098,3965.939)
	(10746.834,3989.139)
	(10688.613,4011.765)
	(10629.484,4033.821)
	(10569.497,4055.311)
	(10508.700,4076.240)
	(10447.142,4096.613)
	(10384.874,4116.433)
	(10321.944,4135.706)
	(10258.402,4154.435)
	(10194.296,4172.626)
	(10129.677,4190.282)
	(10064.593,4207.409)
	(9999.093,4224.010)
	(9933.227,4240.091)
	(9867.045,4255.655)
	(9800.595,4270.708)
	(9733.926,4285.253)
	(9667.089,4299.296)
	(9600.131,4312.840)
	(9533.103,4325.890)
	(9466.054,4338.451)
	(9399.033,4350.528)
	(9332.089,4362.123)
	(9265.271,4373.243)
	(9198.629,4383.892)
	(9132.212,4394.074)
	(9066.069,4403.793)
	(9000.250,4413.054)
	(8934.804,4421.862)
	(8869.780,4430.221)
	(8805.227,4438.136)
	(8741.194,4445.610)
	(8677.731,4452.649)
	(8614.888,4459.257)
	(8552.713,4465.439)
	(8491.255,4471.198)
	(8430.564,4476.540)
	(8370.690,4481.469)
	(8311.680,4485.990)
	(8253.586,4490.106)
	(8196.455,4493.823)
	(8140.337,4497.144)
	(8085.282,4500.075)
	(8031.339,4502.620)
	(7978.556,4504.783)
	(7926.984,4506.569)
	(7876.671,4507.983)
	(7827.667,4509.028)
	(7780.021,4509.710)
	(7733.782,4510.032)
	(7689.000,4510.000)

\path(7689,4510)	(7645.911,4509.619)
	(7601.424,4508.875)
	(7555.587,4507.765)
	(7508.445,4506.284)
	(7460.048,4504.430)
	(7410.442,4502.198)
	(7359.675,4499.586)
	(7307.794,4496.590)
	(7254.847,4493.207)
	(7200.881,4489.433)
	(7145.944,4485.265)
	(7090.083,4480.699)
	(7033.345,4475.732)
	(6975.778,4470.361)
	(6917.430,4464.582)
	(6858.348,4458.391)
	(6798.579,4451.785)
	(6738.171,4444.761)
	(6677.171,4437.315)
	(6615.626,4429.444)
	(6553.585,4421.145)
	(6491.094,4412.413)
	(6428.201,4403.246)
	(6364.954,4393.640)
	(6301.400,4383.591)
	(6237.586,4373.097)
	(6173.559,4362.154)
	(6109.368,4350.758)
	(6045.060,4338.905)
	(5980.681,4326.594)
	(5916.281,4313.819)
	(5851.905,4300.577)
	(5787.602,4286.866)
	(5723.418,4272.682)
	(5659.402,4258.021)
	(5595.601,4242.879)
	(5532.062,4227.254)
	(5468.833,4211.142)
	(5405.960,4194.540)
	(5343.493,4177.443)
	(5281.477,4159.849)
	(5219.961,4141.754)
	(5158.992,4123.155)
	(5098.617,4104.048)
	(5038.884,4084.431)
	(4979.840,4064.298)
	(4921.533,4043.648)
	(4864.010,4022.476)
	(4807.318,4000.779)
	(4751.506,3978.554)
	(4696.620,3955.797)
	(4642.708,3932.504)
	(4589.817,3908.674)
	(4537.995,3884.301)
	(4487.290,3859.382)
	(4437.748,3833.915)
	(4389.417,3807.895)
	(4342.345,3781.319)
	(4296.579,3754.184)
	(4252.167,3726.487)
	(4209.155,3698.223)
	(4167.592,3669.389)
	(4127.524,3639.983)
	(4089.000,3610.000)

\path(4089,3610)	(4050.682,3570.139)
	(4015.253,3511.540)
	(3997.523,3473.016)
	(3979.197,3427.171)
	(3959.836,3373.125)
	(3939.000,3310.000)

\path(4839,1210)	(4859.541,1146.813)
	(4878.691,1092.723)
	(4896.890,1046.850)
	(4914.578,1008.318)
	(4950.175,949.754)
	(4989.000,910.000)

\path(4989,910)	(5046.947,867.248)
	(5109.243,826.450)
	(5175.611,787.565)
	(5245.770,750.553)
	(5319.440,715.374)
	(5357.504,698.459)
	(5396.342,681.987)
	(5435.917,665.952)
	(5476.196,650.351)
	(5517.142,635.177)
	(5558.722,620.425)
	(5600.900,606.092)
	(5643.641,592.170)
	(5686.911,578.657)
	(5730.674,565.545)
	(5774.895,552.831)
	(5819.540,540.509)
	(5864.573,528.574)
	(5909.959,517.021)
	(5955.664,505.845)
	(6001.653,495.041)
	(6047.891,484.604)
	(6094.342,474.529)
	(6140.972,464.810)
	(6187.745,455.443)
	(6234.628,446.423)
	(6281.584,437.744)
	(6328.579,429.401)
	(6375.578,421.390)
	(6422.546,413.705)
	(6469.448,406.342)
	(6516.249,399.295)
	(6562.915,392.558)
	(6609.409,386.128)
	(6655.698,379.999)
	(6701.746,374.165)
	(6747.518,368.623)
	(6792.979,363.366)
	(6838.095,358.390)
	(6882.831,353.689)
	(6927.150,349.259)
	(6971.020,345.094)
	(7014.404,341.190)
	(7057.267,337.541)
	(7099.575,334.143)
	(7141.293,330.989)
	(7182.385,328.076)
	(7222.818,325.398)
	(7262.555,322.949)
	(7301.562,320.726)
	(7339.804,318.722)
	(7413.852,315.355)
	(7484.421,312.805)
	(7551.230,311.034)
	(7614.000,310.000)

\path(7614,310)	(7683.611,309.554)
	(7757.715,309.759)
	(7796.356,310.124)
	(7836.004,310.675)
	(7876.620,311.417)
	(7918.168,312.360)
	(7960.606,313.509)
	(8003.898,314.872)
	(8048.005,316.458)
	(8092.887,318.272)
	(8138.507,320.322)
	(8184.826,322.616)
	(8231.805,325.161)
	(8279.405,327.964)
	(8327.589,331.033)
	(8376.317,334.375)
	(8425.550,337.997)
	(8475.252,341.907)
	(8525.381,346.112)
	(8575.901,350.619)
	(8626.773,355.435)
	(8677.957,360.569)
	(8729.416,366.027)
	(8781.111,371.816)
	(8833.002,377.944)
	(8885.053,384.419)
	(8937.223,391.247)
	(8989.475,398.437)
	(9041.770,405.994)
	(9094.069,413.928)
	(9146.333,422.244)
	(9198.525,430.950)
	(9250.605,440.054)
	(9302.536,449.563)
	(9354.277,459.484)
	(9405.792,469.825)
	(9457.040,480.593)
	(9507.984,491.795)
	(9558.585,503.439)
	(9608.805,515.532)
	(9658.604,528.081)
	(9707.945,541.093)
	(9756.788,554.577)
	(9805.095,568.539)
	(9852.828,582.986)
	(9899.948,597.927)
	(9946.416,613.368)
	(9992.194,629.316)
	(10037.243,645.779)
	(10081.525,662.765)
	(10125.000,680.280)
	(10167.631,698.332)
	(10209.379,716.929)
	(10250.205,736.077)
	(10290.071,755.784)
	(10328.938,776.057)
	(10403.520,818.332)
	(10473.644,862.960)
	(10539.000,910.000)

\path(10539,910)	(10578.398,949.630)
	(10614.188,1008.152)
	(10631.827,1046.696)
	(10649.883,1092.599)
	(10668.794,1146.741)
	(10689.000,1210.000)

\path(39,1210)	(35.278,1249.349)
	(31.796,1287.328)
	(25.553,1359.314)
	(20.271,1426.232)
	(15.949,1488.358)
	(12.587,1545.965)
	(10.186,1599.329)
	(8.745,1648.724)
	(8.265,1694.425)
	(8.745,1736.706)
	(10.186,1775.843)
	(15.949,1845.780)
	(25.553,1906.433)
	(39.000,1960.000)

\path(39,1960)	(54.931,1998.876)
	(82.897,2045.819)
	(126.415,2106.102)
	(155.105,2142.895)
	(189.000,2185.000)

\path(132.806,2056.651)(189.000,2185.000)(74.660,2104.021)
\path(339,2410)	(390.476,2476.718)
	(436.202,2533.330)
	(477.055,2580.715)
	(513.915,2619.753)
	(579.171,2676.299)
	(639.000,2710.000)

\path(639,2710)	(674.073,2719.799)
	(719.749,2723.065)
	(781.300,2719.799)
	(819.677,2715.716)
	(864.000,2710.000)

\path(1014,2710)	(1065.861,2678.001)
	(1113.291,2647.288)
	(1156.428,2617.696)
	(1195.409,2589.060)
	(1261.452,2533.996)
	(1312.519,2480.778)
	(1349.708,2428.087)
	(1374.120,2374.607)
	(1386.850,2319.017)
	(1389.000,2260.000)

\path(1389,2260)	(1384.104,2216.621)
	(1373.767,2175.120)
	(1358.492,2135.530)
	(1338.780,2097.885)
	(1288.056,2028.561)
	(1225.612,1967.418)
	(1155.468,1914.722)
	(1081.639,1870.742)
	(1008.144,1835.745)
	(939.000,1810.000)

\path(939,1810)	(887.951,1798.748)
	(830.528,1794.705)
	(769.174,1797.212)
	(706.335,1805.613)
	(644.454,1819.249)
	(585.975,1837.464)
	(533.342,1859.600)
	(489.000,1885.000)

\path(489,1885)	(452.399,1917.200)
	(418.128,1960.202)
	(386.233,2010.699)
	(356.760,2065.390)
	(329.754,2120.969)
	(305.262,2174.133)
	(283.329,2221.578)
	(264.000,2260.000)

\path(264,2260)	(239.622,2303.362)
	(209.482,2356.074)
	(175.927,2415.300)
	(141.300,2478.209)
	(107.948,2541.966)
	(78.217,2603.740)
	(54.453,2660.695)
	(39.000,2710.000)

\path(39,2710)	(30.771,2752.174)
	(24.892,2800.212)
	(21.366,2855.870)
	(20.190,2920.908)
	(20.484,2957.492)
	(21.366,2997.081)
	(22.835,3039.894)
	(24.892,3086.149)
	(27.538,3136.068)
	(30.771,3189.870)
	(34.591,3247.774)
	(39.000,3310.000)

\path(3939,1210)	(3935.278,1249.349)
	(3931.796,1287.328)
	(3925.553,1359.314)
	(3920.271,1426.232)
	(3915.949,1488.358)
	(3912.587,1545.965)
	(3910.186,1599.329)
	(3908.745,1648.724)
	(3908.265,1694.425)
	(3908.745,1736.706)
	(3910.186,1775.843)
	(3915.949,1845.780)
	(3925.553,1906.433)
	(3939.000,1960.000)

\path(3939,1960)	(3954.931,1998.876)
	(3982.898,2045.819)
	(4026.415,2106.102)
	(4055.105,2142.895)
	(4089.000,2185.000)

\path(4032.806,2056.651)(4089.000,2185.000)(3974.660,2104.021)
\path(4239,2410)	(4290.476,2476.718)
	(4336.202,2533.330)
	(4377.055,2580.715)
	(4413.915,2619.753)
	(4479.171,2676.299)
	(4539.000,2710.000)

\path(4539,2710)	(4574.073,2719.799)
	(4619.749,2723.065)
	(4681.300,2719.799)
	(4719.677,2715.716)
	(4764.000,2710.000)

\path(4914,2710)	(4965.861,2678.001)
	(5013.291,2647.288)
	(5056.428,2617.696)
	(5095.409,2589.060)
	(5161.452,2533.996)
	(5212.519,2480.778)
	(5249.708,2428.087)
	(5274.120,2374.607)
	(5286.850,2319.017)
	(5289.000,2260.000)

\path(5289,2260)	(5284.104,2216.621)
	(5273.767,2175.120)
	(5258.492,2135.530)
	(5238.780,2097.885)
	(5188.056,2028.561)
	(5125.613,1967.418)
	(5055.468,1914.722)
	(4981.639,1870.742)
	(4908.144,1835.745)
	(4839.000,1810.000)

\path(4839,1810)	(4787.951,1798.748)
	(4730.528,1794.705)
	(4669.174,1797.212)
	(4606.335,1805.613)
	(4544.454,1819.249)
	(4485.975,1837.464)
	(4433.342,1859.600)
	(4389.000,1885.000)

\path(4389,1885)	(4352.399,1917.200)
	(4318.128,1960.202)
	(4286.233,2010.699)
	(4256.760,2065.390)
	(4229.754,2120.969)
	(4205.262,2174.133)
	(4183.329,2221.578)
	(4164.000,2260.000)

\path(4164,2260)	(4139.622,2303.362)
	(4109.483,2356.074)
	(4075.927,2415.300)
	(4041.300,2478.209)
	(4007.948,2541.966)
	(3978.218,2603.740)
	(3954.453,2660.695)
	(3939.000,2710.000)

\path(3939,2710)	(3930.771,2752.174)
	(3924.892,2800.212)
	(3921.366,2855.870)
	(3920.190,2920.908)
	(3920.484,2957.492)
	(3921.366,2997.081)
	(3922.835,3039.894)
	(3924.892,3086.149)
	(3927.538,3136.068)
	(3930.771,3189.870)
	(3934.591,3247.774)
	(3939.000,3310.000)

\put(4164,3160){\makebox(0,0)[lb]{\smash{{{\rm\scriptsize dots}}}}}
\put(4239,3385){\makebox(0,0)[lb]{\smash{{{\rm\scriptsize the}}}}}
\put(6864,2185){\makebox(0,0)[lb]{\smash{{{\rm\scriptsize many miles away}}}}}
\end{picture}
}

%% file: draws/Closer.tex
\begingroup\makeatletter\ifx\SetFigFont\undefined
\def\x#1#2#3#4#5#6#7\relax{\def\x{#1#2#3#4#5#6}}%
\expandafter\x\fmtname xxxxxx\relax \def\y{splain}%
\ifx\x\y   
\gdef\SetFigFont#1#2#3{%
  \ifnum #1<17\tiny\else \ifnum #1<20\small\else
  \ifnum #1<24\normalsize\else \ifnum #1<29\large\else
  \ifnum #1<34\Large\else \ifnum #1<41\LARGE\else
     \huge\fi\fi\fi\fi\fi\fi
  \csname #3\endcsname}%
\else
\gdef\SetFigFont#1#2#3{\begingroup
  \count@#1\relax \ifnum 25<\count@\count@25\fi
  \def\x{\endgroup\@setsize\SetFigFont{#2pt}}%
  \expandafter\x
    \csname \romannumeral\the\count@ pt\expandafter\endcsname
    \csname @\romannumeral\the\count@ pt\endcsname
  \csname #3\endcsname}%
\fi
\fi\endgroup
{\renewcommand{\dashlinestretch}{30}
\begin{picture}(2645,2207)(0,-10)
\path(87,871)(87,946)
\path(12,871)(612,871)(612,796)
	(12,796)(12,871)
\path(12,1396)(612,1396)(612,1321)
	(12,1321)(12,1396)
\path(537,871)(537,946)
\path(312,871)(312,946)
\path(312,1246)(312,1321)
\path(537,1246)(537,1321)
\path(87,1246)(87,1321)
\path(537,1396)(537,1321)
\path(462,1396)(312,1321)
\path(87,1396)(87,1321)
\path(87,871)(87,796)
\path(312,871)(462,796)
\path(537,871)(537,796)
\path(537,1396)	(555.270,1454.990)
	(572.689,1497.310)
	(612.000,1546.000)

\path(612,1546)	(680.845,1589.885)
	(717.813,1610.764)
	(756.337,1630.873)
	(796.312,1650.167)
	(837.630,1668.599)
	(880.187,1686.125)
	(923.877,1702.698)
	(968.594,1718.273)
	(1014.231,1732.804)
	(1060.684,1746.246)
	(1107.846,1758.552)
	(1155.611,1769.678)
	(1203.873,1779.577)
	(1252.527,1788.204)
	(1301.468,1795.514)
	(1350.587,1801.460)
	(1399.782,1805.997)
	(1448.944,1809.079)
	(1497.969,1810.661)
	(1546.750,1810.697)
	(1595.182,1809.142)
	(1643.158,1805.949)
	(1690.574,1801.073)
	(1737.323,1794.468)
	(1783.299,1786.089)
	(1828.396,1775.891)
	(1872.509,1763.826)
	(1915.532,1749.851)
	(1957.359,1733.918)
	(1997.883,1715.983)
	(2037.000,1696.000)

\path(2037,1696)	(2104.915,1647.944)
	(2164.727,1584.929)
	(2216.143,1510.587)
	(2238.611,1470.302)
	(2258.869,1428.546)
	(2276.881,1385.773)
	(2292.611,1342.437)
	(2306.021,1298.992)
	(2317.076,1255.890)
	(2325.738,1213.586)
	(2331.970,1172.534)
	(2335.736,1133.188)
	(2337.000,1096.000)

\path(2337,1096)	(2335.736,1058.812)
	(2331.970,1019.466)
	(2325.738,978.414)
	(2317.076,936.110)
	(2306.021,893.008)
	(2292.611,849.563)
	(2276.881,806.227)
	(2258.869,763.454)
	(2238.611,721.698)
	(2216.143,681.413)
	(2164.727,607.071)
	(2104.915,544.056)
	(2037.000,496.000)

\path(2037,496)	(1997.882,476.016)
	(1957.357,458.080)
	(1915.530,442.147)
	(1872.507,428.171)
	(1828.393,416.107)
	(1783.295,405.908)
	(1737.319,397.529)
	(1690.570,390.925)
	(1643.154,386.049)
	(1595.177,382.856)
	(1546.745,381.301)
	(1497.964,381.337)
	(1448.940,382.919)
	(1399.777,386.002)
	(1350.584,390.540)
	(1301.464,396.486)
	(1252.524,403.796)
	(1203.870,412.424)
	(1155.608,422.324)
	(1107.843,433.450)
	(1060.681,445.757)
	(1014.229,459.198)
	(968.592,473.730)
	(923.876,489.305)
	(880.186,505.878)
	(837.629,523.403)
	(796.311,541.836)
	(756.337,561.129)
	(717.813,581.238)
	(680.845,602.116)
	(612.000,646.000)

\path(612,646)	(572.689,694.690)
	(555.270,737.010)
	(537.000,796.000)

\path(462,1396)	(484.802,1444.477)
	(505.567,1485.821)
	(542.745,1549.746)
	(612.000,1621.000)

\path(612,1621)	(680.116,1659.206)
	(716.824,1677.690)
	(755.150,1695.686)
	(794.984,1713.137)
	(836.217,1729.987)
	(878.736,1746.178)
	(922.433,1761.654)
	(967.196,1776.358)
	(1012.915,1790.234)
	(1059.480,1803.224)
	(1106.779,1815.272)
	(1154.704,1826.321)
	(1203.142,1836.314)
	(1251.984,1845.195)
	(1301.119,1852.908)
	(1350.436,1859.394)
	(1399.826,1864.597)
	(1449.178,1868.462)
	(1498.381,1870.930)
	(1547.324,1871.945)
	(1595.898,1871.451)
	(1643.992,1869.391)
	(1691.495,1865.707)
	(1738.297,1860.344)
	(1784.288,1853.244)
	(1829.356,1844.351)
	(1873.392,1833.608)
	(1916.284,1820.958)
	(1957.924,1806.344)
	(1998.199,1789.711)
	(2037.000,1771.000)

\path(2037,1771)	(2077.234,1747.214)
	(2115.873,1718.779)
	(2187.811,1649.846)
	(2220.836,1610.288)
	(2251.712,1567.964)
	(2280.304,1523.343)
	(2306.471,1476.895)
	(2330.077,1429.091)
	(2350.984,1380.402)
	(2369.052,1331.297)
	(2384.145,1282.247)
	(2396.125,1233.721)
	(2404.852,1186.192)
	(2410.190,1140.128)
	(2412.000,1096.000)

\path(2412,1096)	(2410.190,1051.872)
	(2404.852,1005.809)
	(2396.125,958.279)
	(2384.145,909.755)
	(2369.052,860.705)
	(2350.984,811.601)
	(2330.077,762.912)
	(2306.471,715.109)
	(2280.304,668.662)
	(2251.712,624.041)
	(2220.836,581.716)
	(2187.811,542.158)
	(2115.873,473.224)
	(2077.234,444.788)
	(2037.000,421.000)

\path(2037,421)	(1998.198,402.288)
	(1957.922,385.654)
	(1916.282,371.040)
	(1873.389,358.389)
	(1829.353,347.646)
	(1784.284,338.752)
	(1738.293,331.652)
	(1691.491,326.289)
	(1643.988,322.605)
	(1595.894,320.544)
	(1547.320,320.050)
	(1498.376,321.066)
	(1449.174,323.534)
	(1399.822,327.398)
	(1350.432,332.602)
	(1301.115,339.089)
	(1251.980,346.801)
	(1203.139,355.682)
	(1154.701,365.676)
	(1106.777,376.726)
	(1059.477,388.774)
	(1012.913,401.764)
	(967.194,415.640)
	(922.432,430.345)
	(878.735,445.821)
	(836.216,462.012)
	(794.984,478.862)
	(755.149,496.313)
	(716.823,514.310)
	(680.116,532.794)
	(612.000,571.000)

\path(612,571)	(542.745,642.254)
	(505.567,706.179)
	(484.802,747.523)
	(462.000,796.000)

\path(87,1396)	(126.068,1455.511)
	(162.909,1510.565)
	(197.686,1561.355)
	(230.564,1608.071)
	(261.709,1650.908)
	(291.284,1690.057)
	(346.388,1758.059)
	(397.192,1813.616)
	(445.017,1858.267)
	(491.180,1893.549)
	(537.000,1921.000)

\path(537,1921)	(608.092,1955.694)
	(646.536,1972.747)
	(686.752,1989.515)
	(728.615,2005.931)
	(772.006,2021.927)
	(816.801,2037.436)
	(862.880,2052.392)
	(910.119,2066.728)
	(958.399,2080.376)
	(1007.596,2093.270)
	(1057.588,2105.342)
	(1108.255,2116.525)
	(1159.475,2126.753)
	(1211.124,2135.957)
	(1263.082,2144.072)
	(1315.228,2151.031)
	(1367.438,2156.765)
	(1419.591,2161.209)
	(1471.566,2164.294)
	(1523.241,2165.955)
	(1574.493,2166.123)
	(1625.201,2164.733)
	(1675.244,2161.716)
	(1724.499,2157.007)
	(1772.845,2150.537)
	(1820.159,2142.240)
	(1866.321,2132.048)
	(1911.207,2119.896)
	(1954.697,2105.715)
	(1996.669,2089.439)
	(2037.000,2071.000)

\path(2037,2071)	(2096.428,2037.855)
	(2154.578,1998.184)
	(2211.098,1952.593)
	(2265.635,1901.684)
	(2317.835,1846.063)
	(2367.347,1786.334)
	(2413.818,1723.100)
	(2456.895,1656.966)
	(2496.225,1588.537)
	(2531.456,1518.416)
	(2562.234,1447.207)
	(2588.208,1375.515)
	(2609.024,1303.944)
	(2624.330,1233.099)
	(2633.773,1163.583)
	(2637.000,1096.000)

\path(2637,1096)	(2633.776,1028.418)
	(2624.336,958.902)
	(2609.032,888.056)
	(2588.218,816.486)
	(2562.245,744.795)
	(2531.467,673.587)
	(2496.237,603.467)
	(2456.906,535.037)
	(2413.829,468.904)
	(2367.357,405.671)
	(2317.844,345.941)
	(2265.642,290.320)
	(2211.104,239.411)
	(2154.582,193.819)
	(2096.430,154.147)
	(2037.000,121.000)

\path(2037,121)	(1996.669,102.561)
	(1954.697,86.284)
	(1911.207,72.102)
	(1866.321,59.949)
	(1820.159,49.758)
	(1772.845,41.460)
	(1724.499,34.991)
	(1675.244,30.281)
	(1625.201,27.265)
	(1574.493,25.874)
	(1523.241,26.043)
	(1471.566,27.704)
	(1419.591,30.790)
	(1367.438,35.234)
	(1315.228,40.969)
	(1263.082,47.927)
	(1211.124,56.043)
	(1159.475,65.248)
	(1108.255,75.476)
	(1057.588,86.660)
	(1007.596,98.732)
	(958.399,111.626)
	(910.119,125.274)
	(862.880,139.610)
	(816.801,154.567)
	(772.006,170.076)
	(728.615,186.072)
	(686.752,202.487)
	(646.536,219.255)
	(608.092,236.307)
	(537.000,271.000)

\path(537,271)	(491.180,298.451)
	(445.017,333.733)
	(397.192,378.384)
	(346.388,433.941)
	(291.284,501.943)
	(261.709,541.092)
	(230.564,583.929)
	(197.686,630.645)
	(162.909,681.435)
	(126.068,736.489)
	(87.000,796.000)

\path(987,946)	(944.895,979.895)
	(908.102,1008.585)
	(847.819,1052.102)
	(800.876,1080.069)
	(762.000,1096.000)

\path(762,1096)	(716.312,1105.219)
	(655.770,1108.293)
	(617.732,1107.524)
	(573.343,1105.219)
	(521.726,1101.378)
	(462.000,1096.000)

\path(578.587,1137.321)(462.000,1096.000)(584.316,1077.596)
\put(1137,1033){\makebox(0,0)[lb]{\smash{{{\rm\scriptsize goes}}}}}
\put(1137,745){\makebox(0,0)[lb]{\smash{{{\rm\scriptsize here}}}}}
\put(1287,1321){\makebox(0,0)[lb]{\smash{{\mathmode{\scriptstyle L}}}}}
\end{picture}
}

%% file: draws/DotsInS.tex
\begingroup\makeatletter\ifx\SetFigFont\undefined
\def\x#1#2#3#4#5#6#7\relax{\def\x{#1#2#3#4#5#6}}%
\expandafter\x\fmtname xxxxxx\relax \def\y{splain}%
\ifx\x\y   
\gdef\SetFigFont#1#2#3{%
  \ifnum #1<17\tiny\else \ifnum #1<20\small\else
  \ifnum #1<24\normalsize\else \ifnum #1<29\large\else
  \ifnum #1<34\Large\else \ifnum #1<41\LARGE\else
     \huge\fi\fi\fi\fi\fi\fi
  \csname #3\endcsname}%
\else
\gdef\SetFigFont#1#2#3{\begingroup
  \count@#1\relax \ifnum 25<\count@\count@25\fi
  \def\x{\endgroup\@setsize\SetFigFont{#2pt}}%
  \expandafter\x
    \csname \romannumeral\the\count@ pt\expandafter\endcsname
    \csname @\romannumeral\the\count@ pt\endcsname
  \csname #3\endcsname}%
\fi
\fi\endgroup
{\renewcommand{\dashlinestretch}{30}
\begin{picture}(1322,3516)(0,-10)
\put(575.500,400.500){\arc{525.595}{3.1892}{6.2356}}
\put(575.500,1200.500){\arc{530.330}{0.1419}{2.9997}}
\put(912,2739){\blacken\ellipse{150}{150}}
\put(912,2739){\ellipse{150}{150}}
\put(312,3264){\blacken\ellipse{150}{150}}
\put(312,3264){\ellipse{150}{150}}
\put(988,788){\blacken\ellipse{150}{150}}
\put(988,788){\ellipse{150}{150}}
\put(313,1313){\blacken\ellipse{150}{150}}
\put(313,1313){\ellipse{150}{150}}
\path(312,3489)(312,2289)
\path(342.000,3369.000)(312.000,3489.000)(282.000,3369.000)
\path(912,3489)(912,2289)
\path(882.000,2409.000)(912.000,2289.000)(942.000,2409.000)
\dottedline{45}(12,3489)(1287,3489)(1287,2289)
	(12,2289)(12,3489)
\path(988,1538)(988,338)
\path(958.000,458.000)(988.000,338.000)(1018.000,458.000)
\path(313,1538)(313,1163)
\path(343.000,1418.000)(313.000,1538.000)(283.000,1418.000)
\path(838,1538)(838,1163)
\path(838,413)(838,338)
\path(808.000,458.000)(838.000,338.000)(868.000,458.000)
\path(313,413)(313,338)
\dottedline{45}(13,1538)(1288,1538)(1288,338)
	(13,338)(13,1538)
\put(462,3189){\makebox(0,0)[lb]{\smash{{\mathmode{\scriptstyle x}}}}}
\put(1062,2664){\makebox(0,0)[lb]{\smash{{\mathmode{\scriptstyle y}}}}}
\put(1138,713){\makebox(0,0)[lb]{\smash{{\mathmode{\scriptstyle y}}}}}
\put(463,1238){\makebox(0,0)[lb]{\smash{{\mathmode{\scriptstyle x}}}}}
\put(12,1989){\makebox(0,0)[lb]{\smash{{{\rm\scriptsize before surgery}}}}}
\put(87,39){\makebox(0,0)[lb]{\smash{{{\rm\scriptsize after surgery}}}}}
\end{picture}
}

%% file: draws/DetailedKirbyII.tex
\begingroup\makeatletter\ifx\SetFigFont\undefined
\def\x#1#2#3#4#5#6#7\relax{\def\x{#1#2#3#4#5#6}}%
\expandafter\x\fmtname xxxxxx\relax \def\y{splain}%
\ifx\x\y   
\gdef\SetFigFont#1#2#3{%
  \ifnum #1<17\tiny\else \ifnum #1<20\small\else
  \ifnum #1<24\normalsize\else \ifnum #1<29\large\else
  \ifnum #1<34\Large\else \ifnum #1<41\LARGE\else
     \huge\fi\fi\fi\fi\fi\fi
  \csname #3\endcsname}%
\else
\gdef\SetFigFont#1#2#3{\begingroup
  \count@#1\relax \ifnum 25<\count@\count@25\fi
  \def\x{\endgroup\@setsize\SetFigFont{#2pt}}%
  \expandafter\x
    \csname \romannumeral\the\count@ pt\expandafter\endcsname
    \csname @\romannumeral\the\count@ pt\endcsname
  \csname #3\endcsname}%
\fi
\fi\endgroup
{\renewcommand{\dashlinestretch}{30}
\begin{picture}(8895,2433)(0,-10)
\put(7800.000,1763.000){\arc{150.000}{3.1416}{6.2832}}
\put(7800.000,2063.000){\arc{150.000}{6.2832}{9.4248}}
\dottedline{45}(975,2138)(1425,2138)(1425,1688)
	(975,1688)(975,2138)
\dottedline{45}(7577,2136)(8027,2136)(8027,1686)
	(7577,1686)(7577,2136)
\dottedline{45}(4277,2138)(4727,2138)(4727,1688)
	(4277,1688)(4277,2138)
\path(2550,1313)(3075,1313)
\path(2955.000,1283.000)(3075.000,1313.000)(2955.000,1343.000)
\path(5850,1313)(6375,1313)
\path(6255.000,1283.000)(6375.000,1313.000)(6255.000,1343.000)
\path(1125,413)	(1083.687,431.835)
	(1044.104,450.390)
	(970.020,486.770)
	(902.527,522.358)
	(841.406,557.375)
	(786.438,592.041)
	(737.402,626.574)
	(694.080,661.196)
	(656.250,696.125)
	(596.191,767.785)
	(573.523,804.956)
	(555.469,843.312)
	(541.809,883.076)
	(532.324,924.465)
	(526.794,967.700)
	(525.000,1013.000)

\path(525,1013)	(530.776,1071.798)
	(547.190,1124.856)
	(572.867,1172.796)
	(606.435,1216.239)
	(646.522,1255.807)
	(691.754,1292.120)
	(740.758,1325.801)
	(792.161,1357.471)
	(844.591,1387.751)
	(896.673,1417.263)
	(947.036,1446.628)
	(994.306,1476.468)
	(1037.111,1507.403)
	(1074.077,1540.056)
	(1125.000,1613.000)

\path(1125,1613)	(1127.902,1651.100)
	(1125.000,1688.000)

\path(1125,1688)	(1125.000,1725.500)
	(1125.000,1763.000)

\path(1125,1763)	(1125.000,1800.500)
	(1125.000,1838.000)

\path(1125,1838)	(1125.000,1875.500)
	(1125.000,1913.000)

\path(1125,1913)	(1125.000,1950.500)
	(1125.000,1988.000)

\path(1125,1988)	(1125.000,2025.500)
	(1125.000,2063.000)

\path(1125,2063)	(1125.000,2100.500)
	(1125.000,2138.000)

\path(1125,2138)	(1130.452,2174.341)
	(1125.000,2213.000)

\path(1125,2213)	(1085.598,2258.389)
	(1037.628,2298.179)
	(982.968,2331.363)
	(923.494,2356.936)
	(861.083,2373.892)
	(797.612,2381.225)
	(734.959,2377.930)
	(675.000,2363.000)

\path(675,2363)	(607.374,2332.956)
	(544.362,2294.342)
	(485.836,2248.110)
	(431.667,2195.216)
	(381.728,2136.612)
	(335.892,2073.253)
	(294.029,2006.092)
	(256.012,1936.085)
	(221.714,1864.184)
	(191.007,1791.344)
	(163.762,1718.519)
	(139.852,1646.662)
	(119.148,1576.728)
	(101.524,1509.670)
	(86.850,1446.443)
	(75.000,1388.000)

\path(75,1388)	(68.857,1344.920)
	(65.656,1297.747)
	(65.209,1247.113)
	(67.330,1193.655)
	(71.832,1138.008)
	(78.527,1080.806)
	(87.229,1022.685)
	(97.751,964.280)
	(109.906,906.225)
	(123.507,849.156)
	(138.366,793.708)
	(154.297,740.515)
	(171.113,690.212)
	(188.627,643.436)
	(206.651,600.820)
	(225.000,563.000)

\path(225,563)	(259.737,505.699)
	(306.110,442.690)
	(361.373,376.976)
	(422.782,311.562)
	(487.594,249.454)
	(553.064,193.654)
	(616.447,147.168)
	(675.000,113.000)

\path(675,113)	(734.295,87.847)
	(804.954,64.582)
	(843.471,54.068)
	(883.539,44.517)
	(924.729,36.092)
	(966.611,28.959)
	(1008.756,23.280)
	(1050.732,19.219)
	(1092.111,16.939)
	(1132.463,16.606)
	(1171.358,18.381)
	(1208.365,22.430)
	(1275.000,38.000)

\path(1275,38)	(1321.065,66.144)
	(1350.000,113.000)

\path(1350,113)	(1338.604,176.690)
	(1314.207,215.260)
	(1275.000,263.000)

\path(8026,188)	(8066.371,217.229)
	(8105.188,245.624)
	(8178.290,300.036)
	(8245.552,351.483)
	(8307.220,400.213)
	(8363.544,446.473)
	(8414.769,490.510)
	(8461.143,532.571)
	(8502.914,572.904)
	(8540.327,611.755)
	(8573.632,649.372)
	(8628.900,721.893)
	(8670.697,792.443)
	(8701.000,863.000)

\path(8701,863)	(8725.808,937.755)
	(8737.765,977.569)
	(8749.308,1018.849)
	(8760.349,1061.477)
	(8770.803,1105.333)
	(8780.581,1150.299)
	(8789.597,1196.256)
	(8797.764,1243.086)
	(8804.994,1290.670)
	(8811.201,1338.890)
	(8816.298,1387.625)
	(8820.197,1436.759)
	(8822.811,1486.171)
	(8824.054,1535.744)
	(8823.839,1585.359)
	(8822.078,1634.896)
	(8818.684,1684.238)
	(8813.570,1733.266)
	(8806.650,1781.860)
	(8797.836,1829.902)
	(8787.041,1877.274)
	(8774.179,1923.857)
	(8759.161,1969.532)
	(8741.902,2014.180)
	(8722.313,2057.683)
	(8700.309,2099.921)
	(8675.801,2140.777)
	(8648.704,2180.132)
	(8618.929,2217.866)
	(8551.000,2288.000)

\path(8551,2288)	(8492.462,2327.440)
	(8426.612,2350.307)
	(8356.094,2359.252)
	(8283.550,2356.929)
	(8211.625,2345.989)
	(8142.962,2329.086)
	(8080.206,2308.872)
	(8026.000,2288.000)

\path(8026,2288)	(7982.695,2256.305)
	(7951.000,2213.000)

\path(7951,2213)	(7940.365,2176.224)
	(7951.000,2138.000)

\path(7951,2138)	(7951.000,2138.000)

\path(7951,2138)	(7951.000,2138.000)

\path(7951,2138)	(7958.733,2099.176)
	(7951.000,2063.000)

\path(7951,2063)	(7951.000,2025.500)
	(7951.000,1988.000)

\path(7951,1988)	(7951.000,1950.500)
	(7951.000,1913.000)

\path(7951,1913)	(7951.000,1875.500)
	(7951.000,1838.000)

\path(7951,1838)	(7951.000,1800.500)
	(7951.000,1763.000)

\path(7951,1763)	(7949.432,1725.729)
	(7951.000,1688.000)

\path(7951,1688)	(7963.096,1652.278)
	(7982.493,1609.925)
	(8026.000,1538.000)

\path(8026,1538)	(8059.992,1502.995)
	(8105.027,1465.852)
	(8157.423,1427.282)
	(8213.500,1388.000)
	(8269.577,1348.718)
	(8321.973,1310.148)
	(8367.008,1273.005)
	(8401.000,1238.000)

\path(8401,1238)	(8449.195,1170.196)
	(8469.178,1128.359)
	(8476.000,1088.000)

\path(8476,1088)	(8463.708,1025.674)
	(8436.886,962.024)
	(8399.199,898.581)
	(8354.309,836.874)
	(8305.880,778.433)
	(8257.577,724.787)
	(8213.062,677.466)
	(8176.000,638.000)

\path(8176,638)	(8139.534,608.029)
	(8090.398,579.223)
	(8033.440,549.923)
	(7973.507,518.473)
	(7915.449,483.213)
	(7864.113,442.486)
	(7824.348,394.635)
	(7801.000,338.000)

\path(7801,338)	(7810.319,272.930)
	(7835.235,234.705)
	(7876.000,188.000)

\path(7651,413)	(7579.018,454.431)
	(7513.163,494.313)
	(7453.242,532.865)
	(7399.064,570.307)
	(7350.437,606.858)
	(7307.168,642.739)
	(7235.940,713.368)
	(7183.844,783.950)
	(7149.346,856.244)
	(7138.216,893.583)
	(7130.909,932.008)
	(7127.235,971.741)
	(7127.000,1013.000)

\path(7127,1013)	(7134.551,1069.962)
	(7152.282,1121.658)
	(7178.866,1168.659)
	(7212.979,1211.540)
	(7253.295,1250.871)
	(7298.487,1287.225)
	(7347.231,1321.175)
	(7398.200,1353.294)
	(7450.069,1384.153)
	(7501.513,1414.326)
	(7551.205,1444.384)
	(7597.821,1474.901)
	(7640.034,1506.449)
	(7676.518,1539.599)
	(7727.000,1613.000)

\path(7727,1613)	(7729.903,1651.100)
	(7727.000,1688.000)

\path(7727,1688)	(7727.000,1714.014)
	(7727.000,1763.000)

\path(7877,1763)	(7877.000,1714.014)
	(7877.000,1688.000)

\path(7877,1688)	(7874.097,1651.100)
	(7877.000,1613.000)

\path(7877,1613)	(7904.867,1572.072)
	(7947.395,1530.995)
	(7992.225,1493.421)
	(8027.000,1463.000)

\path(8027,1463)	(8067.328,1430.494)
	(8120.563,1396.257)
	(8181.186,1359.102)
	(8243.679,1317.841)
	(8302.525,1271.288)
	(8352.205,1218.255)
	(8387.203,1157.555)
	(8402.000,1088.000)

\path(8402,1088)	(8399.764,1040.759)
	(8389.921,995.940)
	(8373.363,953.435)
	(8350.978,913.139)
	(8292.293,838.746)
	(8220.988,771.905)
	(8182.828,741.050)
	(8144.183,711.764)
	(8105.945,683.939)
	(8069.003,657.468)
	(8002.568,608.165)
	(7952.000,563.000)

\path(7952,563)	(7898.169,520.074)
	(7826.749,469.794)
	(7761.455,409.866)
	(7726.000,338.000)

\path(7726,338)	(7730.872,295.743)
	(7751.985,254.259)
	(7802.000,188.000)

\path(7802,188)	(7826.626,160.254)
	(7877.000,113.000)

\path(7727,2061)	(7727.000,2109.986)
	(7727.000,2136.000)

\path(7727,2136)	(7732.453,2172.341)
	(7727.000,2211.000)

\path(7727,2211)	(7687.598,2256.389)
	(7639.628,2296.179)
	(7584.968,2329.363)
	(7525.494,2354.936)
	(7463.083,2371.892)
	(7399.612,2379.225)
	(7336.959,2375.930)
	(7277.000,2361.000)

\path(7277,2361)	(7209.374,2330.956)
	(7146.362,2292.342)
	(7087.836,2246.110)
	(7033.667,2193.216)
	(6983.728,2134.612)
	(6937.892,2071.253)
	(6896.029,2004.092)
	(6858.012,1934.085)
	(6823.714,1862.184)
	(6793.007,1789.344)
	(6765.762,1716.519)
	(6741.852,1644.662)
	(6721.148,1574.728)
	(6703.524,1507.670)
	(6688.850,1444.443)
	(6677.000,1386.000)

\path(6677,1386)	(6670.857,1342.920)
	(6667.656,1295.747)
	(6667.209,1245.113)
	(6669.330,1191.655)
	(6673.832,1136.008)
	(6680.527,1078.806)
	(6689.229,1020.685)
	(6699.751,962.280)
	(6711.906,904.225)
	(6725.507,847.156)
	(6740.366,791.708)
	(6756.297,738.515)
	(6773.113,688.212)
	(6790.627,641.436)
	(6808.651,598.820)
	(6827.000,561.000)

\path(6827,561)	(6861.737,503.699)
	(6908.110,440.690)
	(6963.373,374.976)
	(7024.783,309.562)
	(7089.594,247.454)
	(7155.064,191.654)
	(7218.447,145.168)
	(7277.000,111.000)

\path(7277,111)	(7336.295,85.847)
	(7406.954,62.582)
	(7445.471,52.068)
	(7485.539,42.517)
	(7526.729,34.092)
	(7568.611,26.959)
	(7610.756,21.280)
	(7652.732,17.219)
	(7694.111,14.939)
	(7734.463,14.606)
	(7773.358,16.381)
	(7810.365,20.430)
	(7877.000,36.000)

\path(7877,36)	(7923.065,64.144)
	(7952.000,111.000)

\path(7952,111)	(7940.604,174.690)
	(7916.207,213.260)
	(7877.000,261.000)

\path(8027,111)	(8084.550,136.217)
	(8126.390,155.951)
	(8177.000,186.000)

\path(8177,186)	(8248.595,242.814)
	(8289.148,276.552)
	(8331.971,313.356)
	(8376.388,352.860)
	(8421.724,394.698)
	(8467.302,438.505)
	(8512.449,483.915)
	(8556.487,530.562)
	(8598.743,578.080)
	(8638.540,626.104)
	(8675.202,674.267)
	(8708.055,722.204)
	(8736.422,769.549)
	(8759.629,815.936)
	(8777.000,861.000)

\path(8777,861)	(8798.980,934.851)
	(8809.419,974.292)
	(8819.398,1015.249)
	(8828.846,1057.600)
	(8837.693,1101.225)
	(8845.870,1146.003)
	(8853.305,1191.812)
	(8859.929,1238.532)
	(8865.672,1286.042)
	(8870.462,1334.219)
	(8874.232,1382.945)
	(8876.909,1432.096)
	(8878.424,1481.553)
	(8878.707,1531.194)
	(8877.688,1580.899)
	(8875.295,1630.545)
	(8871.461,1680.013)
	(8866.113,1729.181)
	(8859.182,1777.928)
	(8850.599,1826.133)
	(8840.291,1873.675)
	(8828.191,1920.433)
	(8814.226,1966.286)
	(8798.328,2011.112)
	(8780.426,2054.792)
	(8760.449,2097.203)
	(8738.329,2138.225)
	(8713.994,2177.737)
	(8687.374,2215.617)
	(8627.000,2286.000)

\path(8627,2286)	(8565.698,2334.513)
	(8494.097,2368.741)
	(8455.432,2380.984)
	(8415.389,2390.237)
	(8374.368,2396.695)
	(8332.767,2400.551)
	(8290.987,2401.999)
	(8249.425,2401.234)
	(8208.481,2398.448)
	(8168.554,2393.835)
	(8130.043,2387.591)
	(8093.348,2379.907)
	(8027.000,2361.000)

\path(8027,2361)	(7982.623,2336.755)
	(7938.860,2298.495)
	(7901.667,2253.987)
	(7877.000,2211.000)

\path(7877,2211)	(7874.097,2172.900)
	(7877.000,2136.000)

\path(7877,2136)	(7877.000,2109.986)
	(7877.000,2061.000)

\path(4351,413)	(4279.018,454.431)
	(4213.163,494.313)
	(4153.242,532.865)
	(4099.064,570.307)
	(4050.437,606.858)
	(4007.168,642.739)
	(3935.940,713.368)
	(3883.844,783.950)
	(3849.346,856.244)
	(3838.216,893.583)
	(3830.909,932.008)
	(3827.235,971.741)
	(3827.000,1013.000)

\path(3827,1013)	(3834.551,1069.962)
	(3852.282,1121.658)
	(3878.866,1168.659)
	(3912.979,1211.540)
	(3953.295,1250.871)
	(3998.487,1287.225)
	(4047.231,1321.175)
	(4098.200,1353.294)
	(4150.069,1384.153)
	(4201.513,1414.326)
	(4251.205,1444.384)
	(4297.821,1474.901)
	(4340.034,1506.449)
	(4376.518,1539.599)
	(4427.000,1613.000)

\path(4427,1613)	(4429.903,1651.100)
	(4427.000,1688.000)

\path(4427,1688)	(4427.000,1725.500)
	(4427.000,1763.000)

\path(4427,1763)	(4427.000,1800.500)
	(4427.000,1838.000)

\path(4427,1838)	(4427.000,1875.500)
	(4427.000,1913.000)

\path(4427,1913)	(4427.000,1950.500)
	(4427.000,1988.000)

\path(4427,1988)	(4427.000,2025.500)
	(4427.000,2063.000)

\path(4427,2063)	(4427.000,2100.500)
	(4427.000,2138.000)

\path(4427,2138)	(4432.453,2174.341)
	(4427.000,2213.000)

\path(4427,2213)	(4387.598,2258.389)
	(4339.628,2298.179)
	(4284.968,2331.363)
	(4225.494,2356.936)
	(4163.083,2373.892)
	(4099.612,2381.225)
	(4036.959,2377.930)
	(3977.000,2363.000)

\path(3977,2363)	(3909.374,2332.956)
	(3846.362,2294.342)
	(3787.836,2248.110)
	(3733.667,2195.216)
	(3683.728,2136.612)
	(3637.892,2073.253)
	(3596.029,2006.092)
	(3558.012,1936.085)
	(3523.714,1864.184)
	(3493.007,1791.344)
	(3465.762,1718.519)
	(3441.852,1646.662)
	(3421.148,1576.728)
	(3403.524,1509.670)
	(3388.850,1446.443)
	(3377.000,1388.000)

\path(3377,1388)	(3370.857,1344.920)
	(3367.656,1297.747)
	(3367.209,1247.113)
	(3369.330,1193.655)
	(3373.832,1138.008)
	(3380.527,1080.806)
	(3389.229,1022.685)
	(3399.751,964.280)
	(3411.906,906.225)
	(3425.507,849.156)
	(3440.366,793.708)
	(3456.297,740.515)
	(3473.113,690.212)
	(3490.627,643.436)
	(3508.651,600.820)
	(3527.000,563.000)

\path(3527,563)	(3561.737,505.699)
	(3608.110,442.690)
	(3663.373,376.976)
	(3724.783,311.562)
	(3789.594,249.454)
	(3855.064,193.654)
	(3918.447,147.168)
	(3977.000,113.000)

\path(3977,113)	(4036.295,87.847)
	(4106.954,64.582)
	(4145.471,54.068)
	(4185.539,44.517)
	(4226.729,36.092)
	(4268.611,28.959)
	(4310.756,23.280)
	(4352.732,19.219)
	(4394.111,16.939)
	(4434.463,16.606)
	(4473.358,18.381)
	(4510.365,22.430)
	(4577.000,38.000)

\path(4577,38)	(4623.065,66.144)
	(4652.000,113.000)

\path(4652,113)	(4640.604,176.690)
	(4616.207,215.260)
	(4577.000,263.000)

\path(4726,188)	(4766.371,217.229)
	(4805.188,245.624)
	(4878.290,300.036)
	(4945.552,351.483)
	(5007.220,400.213)
	(5063.544,446.473)
	(5114.769,490.510)
	(5161.143,532.571)
	(5202.914,572.904)
	(5240.327,611.755)
	(5273.632,649.372)
	(5328.900,721.893)
	(5370.697,792.443)
	(5401.000,863.000)

\path(5401,863)	(5425.808,937.755)
	(5437.765,977.569)
	(5449.308,1018.849)
	(5460.349,1061.477)
	(5470.803,1105.333)
	(5480.581,1150.299)
	(5489.597,1196.256)
	(5497.764,1243.086)
	(5504.994,1290.670)
	(5511.201,1338.890)
	(5516.298,1387.625)
	(5520.197,1436.759)
	(5522.811,1486.171)
	(5524.054,1535.744)
	(5523.839,1585.359)
	(5522.078,1634.896)
	(5518.684,1684.238)
	(5513.570,1733.266)
	(5506.650,1781.860)
	(5497.836,1829.902)
	(5487.041,1877.274)
	(5474.179,1923.857)
	(5459.161,1969.532)
	(5441.902,2014.180)
	(5422.313,2057.683)
	(5400.309,2099.921)
	(5375.801,2140.777)
	(5348.704,2180.132)
	(5318.929,2217.866)
	(5251.000,2288.000)

\path(5251,2288)	(5192.462,2327.440)
	(5126.613,2350.307)
	(5056.094,2359.252)
	(4983.550,2356.929)
	(4911.625,2345.989)
	(4842.962,2329.086)
	(4780.206,2308.872)
	(4726.000,2288.000)

\path(4726,2288)	(4682.695,2256.305)
	(4651.000,2213.000)

\path(4651,2213)	(4640.365,2176.224)
	(4651.000,2138.000)

\path(4651,2138)	(4651.000,2138.000)

\path(4651,2138)	(4651.000,2138.000)

\path(4651,2138)	(4658.733,2099.176)
	(4651.000,2063.000)

\path(4651,2063)	(4651.000,2025.500)
	(4651.000,1988.000)

\path(4651,1988)	(4651.000,1950.500)
	(4651.000,1913.000)

\path(4651,1913)	(4651.000,1875.500)
	(4651.000,1838.000)

\path(4651,1838)	(4651.000,1800.500)
	(4651.000,1763.000)

\path(4651,1763)	(4649.432,1725.729)
	(4651.000,1688.000)

\path(4651,1688)	(4663.096,1652.278)
	(4682.493,1609.925)
	(4726.000,1538.000)

\path(4726,1538)	(4759.992,1502.995)
	(4805.027,1465.852)
	(4857.423,1427.282)
	(4913.500,1388.000)
	(4969.577,1348.718)
	(5021.973,1310.148)
	(5067.008,1273.005)
	(5101.000,1238.000)

\path(5101,1238)	(5149.195,1170.196)
	(5169.178,1128.359)
	(5176.000,1088.000)

\path(5176,1088)	(5163.708,1025.674)
	(5136.886,962.024)
	(5099.199,898.581)
	(5054.309,836.874)
	(5005.880,778.433)
	(4957.577,724.787)
	(4913.062,677.466)
	(4876.000,638.000)

\path(4876,638)	(4839.534,608.029)
	(4790.398,579.223)
	(4733.440,549.923)
	(4673.508,518.473)
	(4615.449,483.213)
	(4564.113,442.486)
	(4524.348,394.635)
	(4501.000,338.000)

\path(4501,338)	(4510.319,272.930)
	(4535.235,234.705)
	(4576.000,188.000)

\path(1425,113)	(1482.550,138.217)
	(1524.390,157.951)
	(1575.000,188.000)

\path(1575,188)	(1646.595,244.814)
	(1687.148,278.552)
	(1729.971,315.356)
	(1774.388,354.860)
	(1819.724,396.698)
	(1865.302,440.505)
	(1910.449,485.915)
	(1954.487,532.562)
	(1996.743,580.080)
	(2036.540,628.104)
	(2073.202,676.267)
	(2106.055,724.204)
	(2134.422,771.549)
	(2157.629,817.936)
	(2175.000,863.000)

\path(2175,863)	(2196.980,936.851)
	(2207.419,976.292)
	(2217.398,1017.249)
	(2226.846,1059.600)
	(2235.693,1103.225)
	(2243.870,1148.003)
	(2251.305,1193.812)
	(2257.929,1240.532)
	(2263.672,1288.042)
	(2268.462,1336.219)
	(2272.232,1384.945)
	(2274.909,1434.096)
	(2276.424,1483.553)
	(2276.707,1533.194)
	(2275.688,1582.899)
	(2273.295,1632.545)
	(2269.461,1682.013)
	(2264.113,1731.181)
	(2257.182,1779.928)
	(2248.599,1828.133)
	(2238.291,1875.675)
	(2226.191,1922.433)
	(2212.226,1968.286)
	(2196.328,2013.112)
	(2178.426,2056.792)
	(2158.449,2099.203)
	(2136.329,2140.225)
	(2111.994,2179.737)
	(2085.374,2217.617)
	(2025.000,2288.000)

\path(2025,2288)	(1963.698,2336.513)
	(1892.097,2370.741)
	(1853.432,2382.984)
	(1813.389,2392.237)
	(1772.368,2398.695)
	(1730.767,2402.551)
	(1688.987,2403.999)
	(1647.425,2403.234)
	(1606.481,2400.448)
	(1566.554,2395.835)
	(1528.043,2389.591)
	(1491.348,2381.907)
	(1425.000,2363.000)

\path(1425,2363)	(1380.623,2338.755)
	(1336.860,2300.495)
	(1299.667,2255.987)
	(1275.000,2213.000)

\path(1275,2213)	(1272.098,2174.900)
	(1275.000,2138.000)

\path(1275,2138)	(1275.000,2100.500)
	(1275.000,2063.000)

\path(1275,2063)	(1275.000,2025.500)
	(1275.000,1988.000)

\path(1275,1988)	(1275.000,1950.500)
	(1275.000,1913.000)

\path(1275,1913)	(1275.000,1875.500)
	(1275.000,1838.000)

\path(1275,1838)	(1275.000,1800.500)
	(1275.000,1763.000)

\path(1275,1763)	(1275.000,1725.500)
	(1275.000,1688.000)

\path(1275,1688)	(1272.098,1651.100)
	(1275.000,1613.000)

\path(1275,1613)	(1302.867,1572.072)
	(1345.395,1530.995)
	(1390.225,1493.421)
	(1425.000,1463.000)

\path(1425,1463)	(1465.328,1430.494)
	(1518.563,1396.257)
	(1579.186,1359.102)
	(1641.679,1317.841)
	(1700.525,1271.288)
	(1750.205,1218.255)
	(1785.203,1157.555)
	(1800.000,1088.000)

\path(1800,1088)	(1797.634,1040.874)
	(1787.436,996.367)
	(1770.349,954.329)
	(1747.315,914.606)
	(1687.174,841.496)
	(1614.551,775.816)
	(1575.915,745.382)
	(1536.985,716.347)
	(1462.012,661.868)
	(1397.171,611.159)
	(1350.000,563.000)

\path(1350,563)	(1313.572,516.187)
	(1268.378,456.406)
	(1226.495,393.672)
	(1200.000,338.000)

\path(1200,338)	(1190.096,263.480)
	(1191.932,222.128)
	(1200.000,188.000)

\path(1200,188)	(1222.564,156.729)
	(1275.000,113.000)

\path(4727,113)	(4784.550,138.217)
	(4826.390,157.951)
	(4877.000,188.000)

\path(4877,188)	(4948.595,244.814)
	(4989.148,278.552)
	(5031.971,315.356)
	(5076.388,354.860)
	(5121.724,396.698)
	(5167.302,440.505)
	(5212.449,485.915)
	(5256.487,532.562)
	(5298.743,580.080)
	(5338.540,628.104)
	(5375.202,676.267)
	(5408.055,724.204)
	(5436.422,771.549)
	(5459.629,817.936)
	(5477.000,863.000)

\path(5477,863)	(5498.980,936.851)
	(5509.419,976.292)
	(5519.398,1017.249)
	(5528.846,1059.600)
	(5537.693,1103.225)
	(5545.870,1148.003)
	(5553.305,1193.812)
	(5559.929,1240.532)
	(5565.672,1288.042)
	(5570.462,1336.219)
	(5574.232,1384.945)
	(5576.909,1434.096)
	(5578.424,1483.553)
	(5578.707,1533.194)
	(5577.688,1582.899)
	(5575.295,1632.545)
	(5571.461,1682.013)
	(5566.113,1731.181)
	(5559.182,1779.928)
	(5550.599,1828.133)
	(5540.291,1875.675)
	(5528.191,1922.433)
	(5514.226,1968.286)
	(5498.328,2013.112)
	(5480.426,2056.792)
	(5460.449,2099.203)
	(5438.329,2140.225)
	(5413.994,2179.737)
	(5387.374,2217.617)
	(5327.000,2288.000)

\path(5327,2288)	(5265.698,2336.513)
	(5194.097,2370.741)
	(5155.432,2382.984)
	(5115.389,2392.237)
	(5074.368,2398.695)
	(5032.767,2402.551)
	(4990.987,2403.999)
	(4949.425,2403.234)
	(4908.481,2400.448)
	(4868.554,2395.835)
	(4830.043,2389.591)
	(4793.348,2381.907)
	(4727.000,2363.000)

\path(4727,2363)	(4682.623,2338.755)
	(4638.860,2300.495)
	(4601.667,2255.987)
	(4577.000,2213.000)

\path(4577,2213)	(4574.097,2174.900)
	(4577.000,2138.000)

\path(4577,2138)	(4577.000,2100.500)
	(4577.000,2063.000)

\path(4577,2063)	(4577.000,2025.500)
	(4577.000,1988.000)

\path(4577,1988)	(4577.000,1950.500)
	(4577.000,1913.000)

\path(4577,1913)	(4577.000,1875.500)
	(4577.000,1838.000)

\path(4577,1838)	(4577.000,1800.500)
	(4577.000,1763.000)

\path(4577,1763)	(4577.000,1725.500)
	(4577.000,1688.000)

\path(4577,1688)	(4574.097,1651.100)
	(4577.000,1613.000)

\path(4577,1613)	(4604.867,1572.072)
	(4647.395,1530.995)
	(4692.225,1493.421)
	(4727.000,1463.000)

\path(4727,1463)	(4767.328,1430.494)
	(4820.563,1396.257)
	(4881.186,1359.102)
	(4943.679,1317.841)
	(5002.525,1271.288)
	(5052.205,1218.255)
	(5087.203,1157.555)
	(5102.000,1088.000)

\path(5102,1088)	(5099.764,1040.759)
	(5089.921,995.940)
	(5073.363,953.435)
	(5050.978,913.139)
	(4992.293,838.746)
	(4920.988,771.905)
	(4882.828,741.050)
	(4844.183,711.764)
	(4805.945,683.939)
	(4769.003,657.468)
	(4702.568,608.165)
	(4652.000,563.000)

\path(4652,563)	(4598.169,520.074)
	(4526.749,469.794)
	(4461.455,409.866)
	(4426.000,338.000)

\path(4426,338)	(4430.872,295.743)
	(4451.985,254.259)
	(4502.000,188.000)

\path(4502,188)	(4526.626,160.254)
	(4577.000,113.000)

\put(5175,113){\makebox(0,0)[lb]{\smash{{\mathmode{2}}}}}
\put(1875,113){\makebox(0,0)[lb]{\smash{{\mathmode{1}}}}}
\put(8475,113){\makebox(0,0)[lb]{\smash{{\mathmode{3}}}}}
\put(600,938){\makebox(0,0)[lb]{\smash{{\mathmode{\scriptstyle x}}}}}
\put(1575,938){\makebox(0,0)[lb]{\smash{{\mathmode{\scriptstyle y}}}}}
\put(1500,1688){\makebox(0,0)[lb]{\smash{{\mathmode{\scriptstyle S}}}}}
\put(0,2138){\makebox(0,0)[lb]{\smash{{\mathmode{L_1}}}}}
\put(6600,2138){\makebox(0,0)[lb]{\smash{{\mathmode{L_2}}}}}
\put(3900,938){\makebox(0,0)[lb]{\smash{{\mathmode{\scriptstyle x}}}}}
\put(5625,1763){\makebox(0,0)[lb]{\smash{{\mathmode{\scriptstyle y'}}}}}
\put(7200,938){\makebox(0,0)[lb]{\smash{{\mathmode{\scriptstyle x}}}}}
\put(8625,1613){\makebox(0,0)[lb]{\smash{{\mathmode{\scriptstyle y}}}}}
\put(5325,1613){\makebox(0,0)[lb]{\smash{{\mathmode{\scriptstyle y}}}}}
\end{picture}
}

%% file: Universality.tex
\section{The Universality of the \Arhus{} Integral} \lbl{sec:Universality}

\subsection{What is universality?}
Let us first recall the definition of universality, as presented
in~\cite[section~\ref{I-subsubsec:Universal}]{Aarhus:I}.

\begin{definition} An invariant $U$ of integer homology spheres with
values in $\calA(\emptyset)$ is a ``universal Ohtsuki invariant'' if
\begin{enumerate}
\item The degree $m$ part $U^{(m)}$ of $U$ is of Ohtsuki type $3m$
  (\cite{Ohtsuki:IntegralHomology}).
\item If $OGL$ denotes the Ohtsuki-Garoufalidis-Le map, defined in
  figure~\ref{fig:OGLMap}, from manifold diagrams to formal
  linear combinations of unit framed algebraically split links in $S^3$,
  and $S$ denotes the surgery map from such links to integer homology
  spheres, then
  \[ (U\circ S\circ OGL)(D)=D+(\text{higher degree diagrams})\qquad
     (\text{in }\calA(\emptyset))
  \]
  whenever $D$ is a manifold diagram (we implicitly linearly extend $S$
  and $U$, to make this a meaningful equation).
\end{enumerate}
\end{definition}

\begin{figure}[htpb]
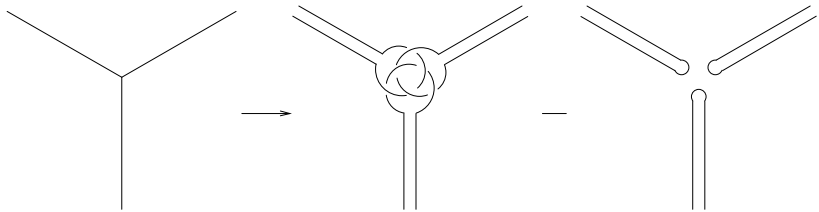

\[ \printname{OGLMap}
	\setlength{\unitlength}{0.5\standardunitlength}
	\begin{array}{c}  \hspace{-1.7mm}
        	\raisebox{-8pt}{\input draws/OGLMap.tex }
        	\hspace{-1.9mm}
	\end{array}
 \]
\caption{
  The $OGL$ map: Take a manifold diagram $D$, embed it in $S^3$ in some
  fixed way of your preference, double every edge, replace every vertex
  by the difference of the two local pictures shown here, and put a
  $+1$ framing on each link component you get. The result is a certain
  alternating sum of $2^v$ links with $e$ components each, where $v$
  and $e$ are the numbers of vertices and edges of $D$, respectively.
}
\lbl{fig:OGLMap}
\end{figure}

\begin{theorem} \lbl{thm:universal}
Restricted to integer homology spheres, $\Aa$ is a universal Ohtsuki
invariant.
\end{theorem}

Some consequences of this theorem were mentioned
in~\cite[corollary~\ref{I-cor:UniversallityMeans}]{Aarhus:I}.

\subsection{$\Aa$ is universal.} \lbl{subsec:Universal}
The proofs of the two properties in the definition of universality
are very similar and both depend on the same principle and
the same observation. Both ideas have been used previously;
see~\cite{Le:UniversalIHS, Bar-Natan:Vassiliev}.

The observation is that the degree $m$ part of $\Aa(L)$ comes from
the internal degree $m$ part of $\Zcheck^+(L)$, the strut-free part of
$\sigma\Zcheck(L)$ in $\calB_n$.  Formal Gaussian integration acts by
connecting all legs of a uni-trivalent diagram to each other using struts. All
univalent vertices disappear in this process, while the trivalent ones
are untouched. And so the degree $m$ part of $\Aa(L)$ is determined by
the internal degree $m$ part of $\Zcheck^+(L)$ (and the linking matrix).

The principle we use is a certain ``locality'' property of the
Kontsevich integral. Recall how the Kontsevich integral of a link $L$
is computed. One sprinkles the link in an arbitrary way with chords,
and takes the resulting chord diagrams with weights that are determined
by the positions of the end points of the chords sprinkled.  This means
that if a localized site on the link get modified, only the weights of
chord diagrams that have ends in that site can change.  Suppose one
marks $k$ localized sites, designates a modification to be made to the
link on each one of them, and computes the alternating sum of $\Zcheck$
evaluated on the $2^k$ links obtained by performing any subset of these
modifications. The result $Z$ must have a chord-end in each of the $k$
sites, and this bounds from below the complexity of any diagram
appearing in $Z$ and constrains the form of the diagrams of least
complexity that appear in $Z$. If more is known about the nature of the
modifications performed, more can be said about the parts of a diagram
$D$ in $Z$ that originate from the sites of the modifications, and thus
more can be said about $D$ altogether.

A very simple application of this principle is the proof of the
universality of the Kontsevich integral in,
say,~\cite{Bar-Natan:Vassiliev}. Two more applications prove
theorem~\ref{thm:universal}.

\begin{proof}[Proof of theorem~\ref{thm:universal}]

{\em $\Aa^{(m)}$ is of Ohtsuki type $3m$: } Take a unit framed
$(k+3m+1)$-component algebraically split link $L$. (That is,
the linking matrix of $L$ is a $(k+3m+1)$-dimensional diagonal
matrix with diagonal entries $\pm 1$). We think of the first $k$
components of $L$ as representing some ``background'' integral homology
sphere, and of the last $3m+1$ components as ``active'' components,
over which the alternating summation in Ohtsuki's definition of
finite-type~\cite{Ohtsuki:IntegralHomology} is performed. Let
$L^{\text{alt}}$ denote that alternating summation. Namely, it is
the alternating sum of the $2^{3m+1}$ sublinks of $L$ in which
some of the active components are removed. We have to show that
$\Aa^{(m)}(L^{\text{alt}})=0$. By the principle, every diagram
in $\Zcheck(L^{\text{alt}})$ must have a chord-end on every active
component of $L$. The map $\sigma$ never `disconnects' a diagram from
a component, and so every diagram $D$ in $\Zcheck^+(L^{\text{alt}})$
must have at least one leg per active component. But the linking
matrix is a diagonal matrix, and hence the struts that are glued in
the Gaussian integration are of form $\strutu{\partial_x}{\partial_x}$
(both ends labeled the same way). So for the Gaussian integration to be
non-trivial, there have to be at least two legs per active component
of $L$, bringing the total to at least $2(3m+1)=6m+2$ legs. Each such
leg must connect to some internal vertex, and there are at most three
legs connected to any internal vertex. So there must by at least $2m+1$
internal vertices, and so the internal degree of $D$ must be higher
than $m$. By the observation, this means that $\Aa(L^{\text{alt}})$
vanishes in degrees up to and including $m$.

{\em $\Aa\circ OGL$ is the identity mod higher degrees: } Let $D$ be a
manifold diagram. We aim to show that
\begin{equation} \lbl{eq:Aim}
  (\Aa\circ OGL)(D)=D+(\text{higher degree diagrams})\qquad
  (\text{in }\calA(\emptyset)).
\end{equation}
If $D$ is of degree $m$, it has $2m$ vertices and $L^{\text{alt}}:=
OGL(D)$ is an alternating summation over modifications in $2m$ sites. By
the principle, there must be a contribution to $\Zcheck(L^{\text{alt}})$
coming from each of those sites. Had there been just one such site,
we would have been looking at the difference $B$ between (a tangle
presentation of) the Borromean rings and a 3-component untangle. As
the Borromean linking numbers are equal to those of the untangle (both
are $0$), there are no struts in $\Zcheck(B)$, and the leading term
is proportional to a $Y$ diagram connecting the three components,
looking like $\hspace{-2mm}\printname{YDiagram}
	\setlength{\unitlength}{0.4\standardunitlength}
	\begin{array}{c}  \hspace{-1.7mm}
        	\raisebox{-8pt}{\input draws/YDiagram.tex }
        	\hspace{-1.9mm}
	\end{array}
\hspace{-2mm}$. A
simple computation shows that the constant of proportionality is $1$
(cf.~\cite{Le:UniversalIHS}).

Ergo, the leading term in $\Zcheck(L^{\text{alt}})$ has a $Y$ piece
corresponding with every vertex of $D$, and the overall coefficient is
$1$. The map $\sigma$ into uni-trivalent diagrams drops the loops
corresponding to the link components and replaces them by labels on the
univalent vertices thus created. (It also adds terms that come from
gluing trees; these terms have a higher internal degree, so, by the
observation, at lowest degree we can ignore them). Gaussian integration
(with an identity covariance matrix, as we have here) simply connects
legs with equal labels using struts, and the result is back again the
diagram $D$ we started with. This process is summarized in
figure~\ref{fig:LeadingTerm}. The renormalization
in~\eqref{eq:Renormalization} doesn't touch any of that, and hence
equation~\eqref{eq:Aim} holds.
\end{proof}

\begin{figure}[htpb]
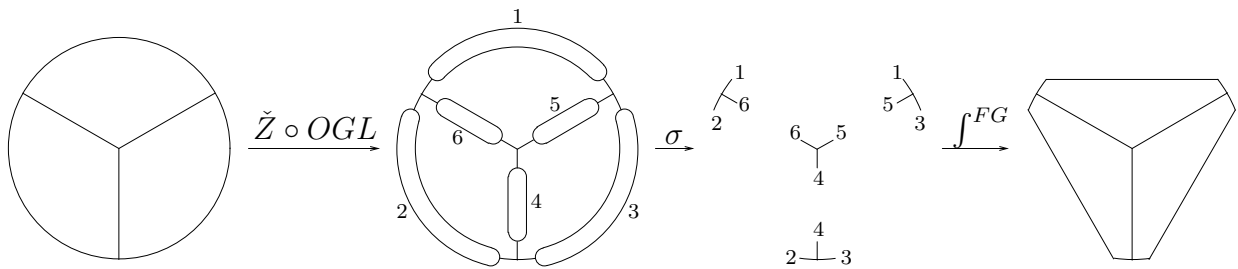

\[
  \def\ZOGL{{\Zcheck\circ OGL}}
  \hspace{-1.5mm}\printname{LeadingTerm}
	\setlength{\unitlength}{0.385\standardunitlength}
	\begin{array}{c}  \hspace{-1.7mm}
        	\raisebox{-8pt}{\input draws/LeadingTerm.tex }
        	\hspace{-1.9mm}
	\end{array}

\]
\caption{The computation of the leading order term in $\Aa_0(OGL(D))$.}
\lbl{fig:LeadingTerm}
\end{figure}

%% file: draws/OGLMap.tex
\begingroup\makeatletter\ifx\SetFigFont\undefined
\def\x#1#2#3#4#5#6#7\relax{\def\x{#1#2#3#4#5#6}}%
\expandafter\x\fmtname xxxxxx\relax \def\y{splain}%
\ifx\x\y   
\gdef\SetFigFont#1#2#3{%
  \ifnum #1<17\tiny\else \ifnum #1<20\small\else
  \ifnum #1<24\normalsize\else \ifnum #1<29\large\else
  \ifnum #1<34\Large\else \ifnum #1<41\LARGE\else
     \huge\fi\fi\fi\fi\fi\fi
  \csname #3\endcsname}%
\else
\gdef\SetFigFont#1#2#3{\begingroup
  \count@#1\relax \ifnum 25<\count@\count@25\fi
  \def\x{\endgroup\@setsize\SetFigFont{#2pt}}%
  \expandafter\x
    \csname \romannumeral\the\count@ pt\expandafter\endcsname
    \csname @\romannumeral\the\count@ pt\endcsname
  \csname #3\endcsname}%
\fi
\fi\endgroup
{\renewcommand{\dashlinestretch}{30}
\begin{picture}(10118,2580)(0,-10)
\put(5009.269,1526.423){\arc{664.088}{5.5962}{7.5267}}
\put(5024.929,1538.786){\arc{575.985}{3.3098}{5.0341}}
\put(5003.500,1474.500){\arc{530.330}{1.7127}{2.9997}}
\path(4966,1212)(4966,12)
\path(5116,1212)(5116,12)
\put(5177.159,1701.017){\arc{662.975}{3.5023}{5.4334}}
\put(5159.031,1708.028){\arc{575.878}{1.2162}{2.9387}}
\put(5224.344,1721.431){\arc{531.244}{5.9030}{7.1854}}
\path(5471,1820)(6510,2420)
\path(5396,1950)(6435,2550)
\put(4941.361,1758.389){\arc{661.861}{1.4049}{3.3412}}
\put(4944.438,1739.488){\arc{576.649}{5.4070}{7.1273}}
\put(4900.431,1789.344){\arc{531.244}{3.8102}{5.0926}}
\path(4692,1954)(3653,2553)
\path(4617,1824)(3578,2424)
\put(8641.000,1418.250){\arc{187.500}{2.4981}{6.9267}}
\path(8566,12)(8566,1362)
\path(8716,12)(8716,1362)
\put(8850.516,1784.740){\arc{187.929}{0.4017}{4.8349}}
\path(10106,2423)(8937,1748)
\path(10031,2553)(8862,1878)
\put(8428.484,1782.740){\arc{187.929}{4.5899}{9.0231}}
\path(7248,2551)(8417,1876)
\path(7173,2421)(8342,1746)
\path(1441,12)(1441,1662)
\path(2870,2487)(1441,1662)
\path(12,2487)(1441,1662)
\path(2941,1212)(3541,1212)
\path(3421.000,1182.000)(3541.000,1212.000)(3421.000,1242.000)
\path(6691,1212)(6991,1212)
\end{picture}
}

%% file: draws/YDiagram.tex
%
\begingroup\makeatletter\ifx\SetFigFont\undefined
\def\x#1#2#3#4#5#6#7\relax{\def\x{#1#2#3#4#5#6}}%
\expandafter\x\fmtname xxxxxx\relax \def\y{splain}%
\ifx\x\y   
\gdef\SetFigFont#1#2#3{%
  \ifnum #1<17\tiny\else \ifnum #1<20\small\else
  \ifnum #1<24\normalsize\else \ifnum #1<29\large\else
  \ifnum #1<34\Large\else \ifnum #1<41\LARGE\else
     \huge\fi\fi\fi\fi\fi\fi
  \csname #3\endcsname}%
\else
\gdef\SetFigFont#1#2#3{\begingroup
  \count@#1\relax \ifnum 25<\count@\count@25\fi
  \def\x{\endgroup\@setsize\SetFigFont{#2pt}}%
  \expandafter\x
    \csname \romannumeral\the\count@ pt\expandafter\endcsname
    \csname @\romannumeral\the\count@ pt\endcsname
  \csname #3\endcsname}%
\fi
\fi\endgroup
{\renewcommand{\dashlinestretch}{30}
\begin{picture}(947,837)(0,-10)
\put(83.590,684.158){\arc{299.303}{5.2328}{8.3835}}
\put(474.000,8.000){\arc{300.000}{3.1416}{6.2832}}
\put(863.410,684.158){\arc{299.304}{1.0413}{4.1920}}
\path(475,459)(735,609)
\path(473,459)(213,609)
\path(474,458)(474,158)
\end{picture}
}

%% file: draws/LeadingTerm.tex
\begingroup\makeatletter\ifx\SetFigFont\undefined
\def\x#1#2#3#4#5#6#7\relax{\def\x{#1#2#3#4#5#6}}%
\expandafter\x\fmtname xxxxxx\relax \def\y{splain}%
\ifx\x\y   
\gdef\SetFigFont#1#2#3{%
  \ifnum #1<17\tiny\else \ifnum #1<20\small\else
  \ifnum #1<24\normalsize\else \ifnum #1<29\large\else
  \ifnum #1<34\Large\else \ifnum #1<41\LARGE\else
     \huge\fi\fi\fi\fi\fi\fi
  \csname #3\endcsname}%
\else
\gdef\SetFigFont#1#2#3{\begingroup
  \count@#1\relax \ifnum 25<\count@\count@25\fi
  \def\x{\endgroup\@setsize\SetFigFont{#2pt}}%
  \expandafter\x
    \csname \romannumeral\the\count@ pt\expandafter\endcsname
    \csname @\romannumeral\the\count@ pt\endcsname
  \csname #3\endcsname}%
\fi
\fi\endgroup
{\renewcommand{\dashlinestretch}{30}
\begin{picture}(19921,4137)(0,-10)
\put(1808,1917){\ellipse{3600}{3600}}
\path(1808,1917)(3367,2817)
\path(1808,1917)(249,2817)
\path(1808,1917)(1808,117)
\put(18240.520,2153.402){\arc{4094.858}{1.4344}{1.7176}}
\put(17976.000,1765.776){\arc{4167.541}{5.6246}{5.9030}}
\put(18449.328,1764.977){\arc{4108.290}{3.5285}{3.8106}}
\path(16844,3036)(19629,3036)
\path(19909,2541)(18509,116)
\path(17934,126)(16549,2546)
\path(18234,1917)(18234,106)
\path(18234,1917)(19784,2822)
\path(18234,1917)(16684,2806)
\put(8278.521,2154.402){\arc{4094.858}{1.4344}{1.7176}}
\put(8014.000,1766.776){\arc{4167.541}{5.6246}{5.9030}}
\put(8487.329,1765.977){\arc{4108.292}{3.5285}{3.8106}}
\put(13139.520,2153.402){\arc{4094.858}{1.4344}{1.7176}}
\put(12875.000,1765.776){\arc{4167.541}{5.6246}{5.9030}}
\put(13348.328,1764.977){\arc{4108.290}{3.5285}{3.8106}}
\put(9569.938,3157.565){\arc{291.428}{5.4231}{8.7628}}
\put(10001.979,2409.301){\arc{290.983}{2.7591}{6.0986}}
\put(8264.505,1889.287){\arc{3347.756}{3.9254}{5.5011}}
\put(8268.554,1922.277){\arc{3887.213}{3.9100}{5.5170}}
\put(9438.158,2581.592){\arc{299.302}{4.1862}{7.3303}}
\put(8658.410,2131.158){\arc{299.304}{1.0413}{4.1920}}
\put(8250.906,1919.953){\arc{3349.705}{6.0198}{7.5949}}
\put(8278.764,1899.253){\arc{3886.520}{6.0044}{7.6112}}
\put(6534.021,2409.301){\arc{290.984}{3.3261}{6.6657}}
\put(6966.062,3157.565){\arc{291.429}{0.6620}{4.0017}}
\put(8286.384,1915.917){\arc{3351.327}{1.8318}{3.4061}}
\put(8252.838,1902.337){\arc{3885.996}{1.8150}{3.4225}}
\put(7097.841,2581.592){\arc{299.301}{2.0945}{5.2386}}
\put(7877.590,2131.158){\arc{299.303}{5.2328}{8.3835}}
\put(8268.000,555.000){\arc{300.000}{6.2832}{9.4248}}
\put(8268.000,1455.000){\arc{300.000}{3.1416}{6.2832}}
\put(8699.588,154.013){\arc{290.627}{1.2312}{4.5772}}
\put(7836.412,154.013){\arc{290.629}{4.8476}{8.1936}}
\path(3908,1917)(6008,1917)
\path(5888.000,1887.000)(6008.000,1917.000)(5888.000,1947.000)
\path(10508,1917)(11108,1917)
\path(10988.000,1887.000)(11108.000,1917.000)(10988.000,1947.000)
\path(15158,1917)(16208,1917)
\path(16088.000,1887.000)(16208.000,1917.000)(16088.000,1947.000)
\path(13128,1905)(12868,2055)
\path(11829,2655)(11569,2805)
\path(13130,1905)(13390,2055)
\path(14429,2655)(14689,2805)
\path(13129,1904)(13129,1604)
\path(13129,404)(13129,104)
\path(8734,2002)(9513,2452)
\path(8584,2261)(9363,2711)
\path(8269,1906)(8529,2056)
\path(9568,2656)(9828,2806)
\path(7952,2261)(7173,2711)
\path(7802,2002)(7023,2452)
\path(8267,1906)(8007,2056)
\path(6968,2656)(6708,2806)
\path(8118,1455)(8118,555)
\path(8418,1455)(8418,555)
\path(8268,1905)(8268,1605)
\path(8268,405)(8268,105)
\put(8183,3942){\makebox(0,0)[lb]{\smash{{\mathmode{\scriptstyle 1}}}}}
\put(8783,2517){\makebox(0,0)[lb]{\smash{{\mathmode{\scriptstyle 5}}}}}
\put(10058,792){\makebox(0,0)[lb]{\smash{{\mathmode{\scriptstyle 3}}}}}
\put(8483,942){\makebox(0,0)[lb]{\smash{{\mathmode{\scriptstyle 4}}}}}
\put(13058,492){\makebox(0,0)[lb]{\smash{{\mathmode{\scriptstyle 4}}}}}
\put(11783,3042){\makebox(0,0)[lb]{\smash{{\mathmode{\scriptstyle 1}}}}}
\put(14333,3042){\makebox(0,0)[lb]{\smash{{\mathmode{\scriptstyle 1}}}}}
\put(13508,42){\makebox(0,0)[lb]{\smash{{\mathmode{\scriptstyle 3}}}}}
\put(13433,2067){\makebox(0,0)[lb]{\smash{{\mathmode{\scriptstyle 5}}}}}
\put(12683,2067){\makebox(0,0)[lb]{\smash{{\mathmode{\scriptstyle 6}}}}}
\put(11858,2517){\makebox(0,0)[lb]{\smash{{\mathmode{\scriptstyle 6}}}}}
\put(7208,1992){\makebox(0,0)[lb]{\smash{{\mathmode{\scriptstyle 6}}}}}
\put(10658,1992){\makebox(0,0)[lb]{\smash{{\mathmode{\sigma}}}}}
\put(15308,2067){\makebox(0,0)[lb]{\smash{{\mathmode{\Gint}}}}}
\put(13058,1317){\makebox(0,0)[lb]{\smash{{\mathmode{\scriptstyle 4}}}}}
\put(3983,2067){\makebox(0,0)[lb]{\smash{{\mathmode{\ZOGL}}}}}
\put(14183,2517){\makebox(0,0)[lb]{\smash{{\mathmode{\scriptstyle 5}}}}}
\put(11408,2217){\makebox(0,0)[lb]{\smash{{\mathmode{\scriptstyle 2}}}}}
\put(14708,2217){\makebox(0,0)[lb]{\smash{{\mathmode{\scriptstyle 3}}}}}
\put(12608,42){\makebox(0,0)[lb]{\smash{{\mathmode{\scriptstyle 2}}}}}
\put(6308,792){\makebox(0,0)[lb]{\smash{{\mathmode{\scriptstyle 2}}}}}
\end{picture}
}

%% file: odds.tex
\section{Odds and ends} \lbl{sec:odds}

\subsection{Homology spheres with embedded links.} \lbl{subsec:EmbeddedLinks}
Everything said in the invariance section of this paper
(section~\ref{sec:Invariance}) holds (or has an obvious counterpart)
in the case of rational homology spheres with embedded links. Most
changes required are completely superficial --- wherever ``components''
are mentioned, of links, string links, dotted Morse links, skeletons
of chord diagrams, etc., one has to label some of the components
as ``surgery components'' (indexed by some set $Y$) and the rest as
``embedded link'' components (indexed by $X$).  Surgeries are performed
only on the components so labeled, $\sigma$ is only applied on those
components, and Gaussian integrations is carried out only with respect
to the variables corresponding to the surgery components. Only the
surgery components count for the purpose of determining $\sigma_\pm$
in~\eqref{eq:Renormalization}.  The embedded link components correspond to
the embedded link in the post-surgery manifold. The only action taken on
embedded link skeleton components (of chord diagrams in $\Zcheck(L)$)
is to take their closures. The target space of the link-enhanced
\Arhus{} integral is a mixture $\calA'(\circlearrowleft_X)$ of
the space $\calA(\udot\circlearrowleft_X)$ of chord diagrams
(mod $4T/STU$) whose skeleton is a disjoint union of $X$-marked
circles (see~\cite[figure~\ref{I-fig:LinkDiagrams}]{Aarhus:I}) and
the space $\calA(\emptyset)$ of manifold diagrams (modulo $AS$ and
$IHX$). The diagrams in $\calA'(\circlearrowleft_X)$ are the disjoint
unions of diagrams in $\calA(\circlearrowleft_X)$ and diagrams in
$\calA(\emptyset)$, and the relations are all the relations mentioned
above.

The only (slight) difficulty is that one should also prove invariance
under the second Kirby move (figure~\ref{fig:DetailedKirbyII}) in the
case where an embedded link component $x$ slides over a surgery
component $y$. A careful reading of the proof of
proposition~\ref{prop:KirbyII} shows that it covers this case as well,
as it uses only the integration with respect to $y$, the surgery
component.

While the link-enhanced target space $\calA'(\circlearrowleft_X)$ suggests
what universality should be like in the case of invariants of integer
homology spheres with embedded links, the necessary preliminaries on
finite-type invariants of such objects where never worked out in
detail. So at this time we do not attempt to generalize the results of
section~\ref{sec:Universality} to the case where embedded links are
present.

\subsection{The link relation.} We (the authors) are not terribly
happy about section~\ref{subsec:DescendsToLinks}. Rather than showing
that $\Aa_0$ descends to regular links, we would have much preferred
to be able to define it directly on regular links. The problem is
that the Kontsevich invariant of $X$-component links is valued in the
space $\calA(\circlearrowleft_X)$ of chord diagrams (mod $4T/STU$)
whose skeleton is a disjoint union of circles marked by the elements
of $X$ (see~\cite[figure~\ref{I-fig:LinkDiagrams}]{Aarhus:I}). This
space is not isomorphic to $\calB(X)$, but rather to a quotient space
$\calB^{\text{links}}(X)$ thereof, and we don't know how to define $\Gint$
on $\calB^{\text{links}}(X)$. Let us write a few more words. First,
a description of $\calB^{\text{links}}(X)$:

\begin{definition} A ``link relation symbol'' is an $X$-marked
uni-trivalent diagram $R^\ast$ one of whose legs is singled out and carries an
additional $\ast$ mark. If the other mark on the special leg of $R^\ast$ is,
say, $x$, we say that $R^\ast$ is an ``$x$-flavored link relation symbol''.
The ``link relation'' $R$ corresponding to an $x$-flavored link relation
symbol $R^\ast$ is the sum of all ways of connecting the $\ast$-marked leg
near the ends of all other $x$-marked legs. It is an element of
$\calB(X)$. An example appears in figure~\ref{fig:LinkRelation}.
Finally, let $\calB^{\text{links}}(X)$ be the quotient of $\calB(X)$ by
all link relations.
\end{definition}

\begin{figure}[htpb]
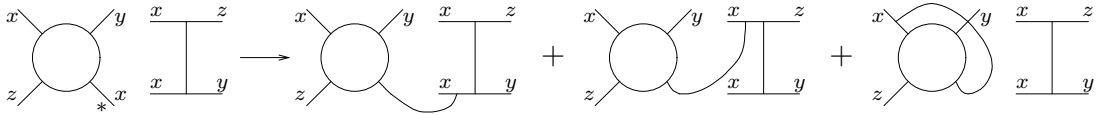

\[ \printname{LinkRelation}
	\setlength{\unitlength}{0.5\standardunitlength}
	\begin{array}{c}  \hspace{-1.7mm}
        	\raisebox{-8pt}{\input draws/LinkRelation.tex }
        	\hspace{-1.9mm}
	\end{array}
 \]
\caption{
  An $x$-flavored link relation symbol and the corresponding link relation.
}
\lbl{fig:LinkRelation}
\end{figure}

\begin{theorem} The isomorphism $\chi:\calB(X)\to\calA(\uparrow_X)$
descends to a well defined isomorphism
$\chi:\calB^{\text{links}}(X)\to\calA(\circlearrowleft_X)$.
\end{theorem}

\begin{proof} (sketch) The fact that the link relations get mapped to $0$
by $\chi$ composed with the projection on $\calA(\circlearrowleft_X)$
is easy --- after applying $\chi$, use the $STU$ relation near every
leg touched by the link relation. On a circular skeleton, the result is
an ouroboros\footnote{
  The medieval symbol of holism depicting a snake that bites its own
  tail. See e.g.{}
  http://\linebreak[0]shell4.\linebreak[0]ba.\linebreak[0]best.\linebreak[0]com/\linebreak[0]\~{}abacus/\linebreak[0]oro/\linebreak[0]ouroboros.html.
}
summation, namely, it is $0$. Suppose now you have a pair of diagrams in
$\calA(\uparrow_X)$ that get identified upon closing one of the
skeleton components, say $y$. Use $STU$ relations as here,
\[ \printname{LinkRelInA}
	\setlength{\unitlength}{0.5\standardunitlength}
	\begin{array}{c}  \hspace{-1.7mm}
        	\raisebox{-8pt}{\input draws/LinkRelInA.tex }
        	\hspace{-1.9mm}
	\end{array}
, \]
to turn their difference into a sum $S$ of diagrams with a lower number
of $y$-legs. Dropping the $y$ component of the skeleton and forgetting
the order of the $y$ legs, the result is a $y$-flavored link relation.
If trees are glued (as $\sigma=\chi^{-1}$ dictates) after the $y$
component of the skeleton is dropped, then using the $IHX$ relation one
can show that the result is still a link relation:
\[ \printname{StillLinkRel}
	\setlength{\unitlength}{0.5\standardunitlength}
	\begin{array}{c}  \hspace{-1.7mm}
        	\raisebox{-8pt}{\input draws/StillLinkRel.tex }
        	\hspace{-1.9mm}
	\end{array}
 \]
\end{proof}

\begin{problem} Is there a good definition for $\Gint$ on (a domain in)
$\calB^{\text{links}}(X)$?
\end{problem}

It makes no sense to ask if $\Gint$ is well defined modulo the link
relation; if $G$ is a Gaussian and $R$ is a link relation, $G+R$ would
no longer be a Gaussian. We are mostly interested in integrating
group-like $G$'s. Maybe there's a more restrictive ``group-like link
relation'' that relates any two group-like elements of $\calA(\uparrow_X)$ that
are equal modulo the usual link relation (namely, whose projections to
$\calA(\circlearrowleft_X)$ are the same)?

\begin{problem} If $G_{1,2}$ are group-like elements of
$\calA(\uparrow_z\uparrow_E)$ that are equal modulo the link relation
(applied only on the $z$ component), is there always a group-like
$G\in \calA(\uparrow_x\uparrow_y\uparrow_E)$ so that $G_1=\mright^{xy}_z G$
and $G_2=\mright^{yx}_z G$?  (notation as in definition~\ref{def:mxyz}).
\end{problem}

\subsection{An explicit formula for the map $\protect\mright^{xy}_z$ on
uni-trivalent diagrams.
} \lbl{subsec:BCH}
Let $x$ and $y$ be two elements in a  free associative (but
not-commutative) completed algebra.  The Baker-Campbel-Hausdorf (BCH)
formula (see e.g.~\cite{Jacobson:LieAlgebras}) measures the failure of
the identity $e^{x+y}=e^xe^y$ to hold, in terms of Lie elements, or, what
is the same, in terms of trees modulo the $IHX$ and $AS$ relation. The
first few terms in the BCH formula are:
\begin{equation} \lbl{eq:BCH}
  \def\scl{0.39}
  \def\neg{\!\!}
  \def\x{{\scriptstyle x}}
  \def\y{{\scriptstyle y}}
  \def\z{{}}
  \begin{array}{llllllll}
    \displaystyle \log e^xe^y
    &\displaystyle\neg = x
    &\displaystyle\neg + y
    &\displaystyle\neg + \frac12[x,y]
    &\displaystyle\neg + \frac1{12}[x,[x,y]]
    &\displaystyle\neg - \frac1{12}[y,[x,y]]
    &\displaystyle\neg - \frac1{24}[x,[y,[x,y]]]
    &\displaystyle\neg + \cdots \\
    \rule{0pt}{30pt}
    &\displaystyle\neg = \!\printname{x}
	\setlength{\unitlength}{\scl\standardunitlength}
	\begin{array}{c}  \hspace{-1.7mm}
        	\raisebox{-8pt}{\input draws/x.tex }
        	\hspace{-1.9mm}
	\end{array}

    &\displaystyle\neg + \!\printname{y}
	\setlength{\unitlength}{\scl\standardunitlength}
	\begin{array}{c}  \hspace{-1.7mm}
        	\raisebox{-8pt}{\input draws/y.tex }
        	\hspace{-1.9mm}
	\end{array}

    &\displaystyle\neg + \frac12\printname{xy}
	\setlength{\unitlength}{\scl\standardunitlength}
	\begin{array}{c}  \hspace{-1.7mm}
        	\raisebox{-8pt}{\input draws/xy.tex }
        	\hspace{-1.9mm}
	\end{array}

    &\displaystyle\neg + \frac1{12}\printname{xxy}
	\setlength{\unitlength}{\scl\standardunitlength}
	\begin{array}{c}  \hspace{-1.7mm}
        	\raisebox{-8pt}{\input draws/xxy.tex }
        	\hspace{-1.9mm}
	\end{array}

    &\displaystyle\neg - \frac1{12}\printname{yxy}
	\setlength{\unitlength}{\scl\standardunitlength}
	\begin{array}{c}  \hspace{-1.7mm}
        	\raisebox{-8pt}{\input draws/yxy.tex }
        	\hspace{-1.9mm}
	\end{array}

    &\displaystyle\neg - \frac1{24}\printname{xyxy}
	\setlength{\unitlength}{\scl\standardunitlength}
	\begin{array}{c}  \hspace{-1.7mm}
        	\raisebox{-8pt}{\input draws/xyxy.tex }
        	\hspace{-1.9mm}
	\end{array}

    &\displaystyle\neg + \cdots.
  \end{array}
\end{equation}

The proposition below states that as an operation on uni-trivalent diagrams,
the map
$\mright^{xy}_z:\calA(\uparrow_x\uparrow_y\uparrow_E)\to\calA\uparrow_z\uparrow_E)$
of definition~\ref{def:mxyz} is given by gluing the disjoint-union
exponential of the trees in the BCH formula~\eqref{eq:BCH}. Precisely,
let $\Lambda$ be the sum of trees in the BCH formula, only with
$\partial_x$ replacing $x$, with $\partial_y$ replacing $y$, and with
a $z$ marked on each root:
\begin{equation} \lbl{eq:BCHDef}
  \def\scl{0.39}
  \def\x{{\scriptstyle \partial_x}}
  \def\y{{\scriptstyle \partial_y}}
  \def\z{{\scriptstyle z}}
  \Lambda = \!\printname{x}
	\setlength{\unitlength}{\scl\standardunitlength}
	\begin{array}{c}  \hspace{-1.7mm}
        	\raisebox{-8pt}{\input draws/x.tex }
        	\hspace{-1.9mm}
	\end{array}
 + \!\printname{y}
	\setlength{\unitlength}{\scl\standardunitlength}
	\begin{array}{c}  \hspace{-1.7mm}
        	\raisebox{-8pt}{\input draws/y.tex }
        	\hspace{-1.9mm}
	\end{array}
 + \frac12\printname{xy}
	\setlength{\unitlength}{\scl\standardunitlength}
	\begin{array}{c}  \hspace{-1.7mm}
        	\raisebox{-8pt}{\input draws/xy.tex }
        	\hspace{-1.9mm}
	\end{array}

    + \frac1{12}\printname{xxy}
	\setlength{\unitlength}{\scl\standardunitlength}
	\begin{array}{c}  \hspace{-1.7mm}
        	\raisebox{-8pt}{\input draws/xxy.tex }
        	\hspace{-1.9mm}
	\end{array}
 - \frac1{12}\printname{yxy}
	\setlength{\unitlength}{\scl\standardunitlength}
	\begin{array}{c}  \hspace{-1.7mm}
        	\raisebox{-8pt}{\input draws/yxy.tex }
        	\hspace{-1.9mm}
	\end{array}

    - \frac1{24}\printname{xyxy}
	\setlength{\unitlength}{\scl\standardunitlength}
	\begin{array}{c}  \hspace{-1.7mm}
        	\raisebox{-8pt}{\input draws/xyxy.tex }
        	\hspace{-1.9mm}
	\end{array}
 + \cdots.
\end{equation}

\begin{proposition} \lbl{prop:ExpBCH}
For any $C\in\calB(\{x,y\}\udot E)$,
\begin{equation} \lbl{eq:ExpBCH}
  \sigma_n\mright^{xy}_z\chi_{n+1}C= \langle\exp_\udot\Lambda,C\rangle_{x,y},
\end{equation}
where $\chi_{n+1}$ denotes the standard isomorphism $\calB(\{x,y\}\udot
E)\to\calA(\uparrow_x\uparrow_y\uparrow_E)$ (whose inverse is
$\sigma_{n+1}$) and $\sigma_n$ denotes the standard isomorphism
$\calA(\uparrow_z\uparrow_E)\to\calB(\{z\}\udot E)$ (whose inverse
is $\chi_n$).
\end{proposition}

A noteworthy special case of this proposition is the case where $n=1$
and $C$ is a disjoint union $C_x\mathbin\udot C_y$ of a uni-trivalent
diagram
$C_x$ whose legs are labeled only by $x$ and a uni-trivalent diagram $C_y$
whose legs are labeled only by $y$. In this case $\mright^{xy}_z$ is
(up to leg labelings) the product $\times_A:\calB\otimes\calB\to\calB$
that $\calB$ inherits from $\calA$, and equation~\eqref{eq:ExpBCH}
becomes a specific formula for this product in terms of gluing
forests. The existence of such a formula is immediate from the
definition of $\sigma:\calA\to\calB$, and this existence was used
in several places before (see e.g.~\cite{Bar-Natan:Homotopy}),
but we are not aware of a previous place where this formula was
written explicitly. A similar formula is the ``wheeling formula''
of~\cite{Bar-NatanGaroufalidisRozanskyThurston:WheelsWheeling}.

\begin{proof}[Proof of proposition~\ref{prop:ExpBCH}]
Let $\calA^{xy}$ be the space of ``planted forests'' whose leaves are
labeled $\partial_x$ and $\partial_y$, modulo the usual $STU$ (and hence
$AS$ and $IHX$) relations. A planted forest is simply a forest in which the
roots of the trees are ``planted'' along a directed line:
\[
  \def\x{{\scriptstyle \partial_x}}
  \def\y{{\scriptstyle \partial_y}}
  \printname{PlantedForest}
	\setlength{\unitlength}{0.5\standardunitlength}
	\begin{array}{c}  \hspace{-1.7mm}
        	\raisebox{-8pt}{\input draws/PlantedForest.tex }
        	\hspace{-1.9mm}
	\end{array}
.
\]
$\calA^{xy}$ is similar in nature to $\calA$; in particular, it is an
algebra by the juxtaposition product $\times_A$ and it is graded, and
hence an exponential $\exp_A$ and a logarithm $\log_A$ can be defined
on it using power series.

Let $\xi$ and $\eta$ be the elements
$\def\x{{\scriptstyle \partial_x}}\setlength{\unitlength}{0.5\standardunitlength}
	\begin{array}{c}  \hspace{-1.7mm}
        	\raisebox{-2pt}{\input draws/xi.tex }
        	\hspace{-1.9mm}
	\end{array}
$ and
$\def\y{{\scriptstyle \partial_y}}\setlength{\unitlength}{0.5\standardunitlength}
	\begin{array}{c}  \hspace{-1.7mm}
        	\raisebox{-2pt}{\input draws/eta.tex }
        	\hspace{-1.9mm}
	\end{array}
$ of $\calA^{xy}$,
respectively. Clearly,
\[
  \mright^{xy}_z\chi_{n+1}C
  = \langle(\exp_A\xi)\times_A(\exp_A\eta),C\rangle_{x,y}
  = \left\langle
    \exp_A\log_A\left((\exp_A\xi)\times_A(\exp_A\eta)\right),C
  \right\rangle_{x,y}.
\]
But $\log_A\left((\exp_A\xi)\times_A(\exp_A\eta)\right)$
can be evaluated using the BCH formula~\eqref{eq:BCH}. The
result is $\chi_z^{xy}\Lambda$, where $\Lambda$ was defined in
equation~\eqref{eq:BCHDef} and $\chi_z^{xy}:\calB_z^{xy}\to\calA^{xy}$ is
the natural isomorphism (whose inverse is $\sigma_z^{xy}$) of the space
$\calB_z^{xy}$ of forests with trees as in equation~\eqref{eq:BCHDef}
(modulo $AS$ and $IHX$) and the space $\calA^{xy}$. Therefore
$\mright^{xy}_z\chi_{n+1}C=\langle\exp_A\chi_z^{xy}\Lambda,C\rangle_{x,y}$,
and hence
\begin{equation} \lbl{eq:SigmaExpChi}
  \sigma_n\mright^{xy}_z\chi_{n+1}C
  = \langle\sigma_z^{xy}\exp_A\chi_z^{xy}\Lambda,C\rangle_{x,y}.
\end{equation}
The only thing left to note is that $\Lambda$ is a sum of {\em
trees}, namely forests in which $z$ appears only once. On such forests
$\exp_A\circ\chi_z^{xy}=\chi_z^{xy}\circ\exp_\udot$, and we see that
equation~\ref{eq:SigmaExpChi} proves equation~\eqref{eq:ExpBCH}.
\end{proof}

\begin{corollary} \lbl{cor:BCH}
(Compare with~\cite[exercise~\ref{I-ex:BCH}]{Aarhus:I})
Let $d_{\BCH}$ be $\Lambda$ with the first two terms removed,
\[
  \def\scl{0.39}
  \def\x{{\scriptstyle \partial_x}}
  \def\y{{\scriptstyle \partial_y}}
  \def\z{{\scriptstyle z}}
  d_{\BCH} = \frac12\printname{xy}
	\setlength{\unitlength}{\scl\standardunitlength}
	\begin{array}{c}  \hspace{-1.7mm}
        	\raisebox{-8pt}{\input draws/xy.tex }
        	\hspace{-1.9mm}
	\end{array}
 + \frac1{12}\printname{xxy}
	\setlength{\unitlength}{\scl\standardunitlength}
	\begin{array}{c}  \hspace{-1.7mm}
        	\raisebox{-8pt}{\input draws/xxy.tex }
        	\hspace{-1.9mm}
	\end{array}

    - \frac1{12}\printname{yxy}
	\setlength{\unitlength}{\scl\standardunitlength}
	\begin{array}{c}  \hspace{-1.7mm}
        	\raisebox{-8pt}{\input draws/yxy.tex }
        	\hspace{-1.9mm}
	\end{array}
 - \frac1{24}\printname{xyxy}
	\setlength{\unitlength}{\scl\standardunitlength}
	\begin{array}{c}  \hspace{-1.7mm}
        	\raisebox{-8pt}{\input draws/xyxy.tex }
        	\hspace{-1.9mm}
	\end{array}
 + \cdots,
\]
and let $D_{\BCH}=\exp_\udot d_{\BCH}$. Then
\[ \sigma_n\mright^{xy}_z\chi_{n+1}C = (D_{\BCH}\apply C)/(x,y\to z). \]
\end{corollary}

\begin{proof} Gluing the exponentials of the struts $|^{\partial_x}_z$ and
$|^{\partial_y}_z$ is equivalent to applying the change of variables
$(x,y\to z)$.
\end{proof}

\subsection{A stronger form of the cyclic invariance lemma.} As stated
and proved in section~\ref{subsec:DescendsToLinks}, the cyclic invariance
lemma holds only when the integration is carried out over {\em all} the
available variables.  Otherwise, the argument given in the proof in the
``lucky case'' would not be a reduction to the case of knots. Here is
an alternative statement and proof, that apply even if integration is
carried out only over some subset $F$ of the relevant variables.

\begin{proposition} (Strong cyclic invariance).
If $G\in\calA(\uparrow_x\uparrow_y\uparrow_E)$ is group-like
and $\sigma\mright^{xy}_z G$ is integrable with respect to some set $F$ of
variables satisfying $z\in F\subset \{z\}\udot E$, then
$\sigma\mright^{yx}_z G$ is also integrable and the two integrals
are equal:
\[
  \Gint\sigma \mright^{xy}_z G\,dF
  = \Gint\sigma \mright^{yx}_z  G\,dF.
\]
\end{proposition}

\begin{proof} Let us describe the idea of the proof before getting into
the details. Recall
from~\cite[section~\ref{I-subsec:DiagrammaticCase}]{Aarhus:I} that we
like to think of $\calA(\uparrow_x\uparrow_y\uparrow_E)$ as a parallel of
$\hat{U}(\frakg)^{\otimes n+1}$ and of $\calB(\{x,y\}\udot E)$ as a parallel of
$\hat{S}(\frakg)^{\otimes n+1}$ or of the function space
$F(\frakg^\star\oplus\cdots(n+1)\cdots\oplus\frakg^\star)$. In this
model, $\mright^{xy}_z$ and $\mright^{yx}_z $ become maps $\hat{U}(\frakg)^{\otimes
n+1}\to\hat{U}(\frakg)^{\otimes n}$, defined by multiplying the first
two tensor factors in the two possible orders, using the product
$\times_U$ of $U(\frakg)$. If we had done the same with
$\hat{S}(\frakg)$, using its product $\times_S$, the picture would have
been a lot simpler. The product $\times_S$ is commutative, and
$m^{xy}_z$ and $m^{yx}_z$ are the same. In fact, in the function space
picture both maps become the ``evaluation on the diagonal'' map $m$:

\begin{eqnarray*}
  F(\frakg^\star\oplus\cdots(n+1)\cdots\oplus\frakg^\star)
  & \to &
  F(\frakg^\star\oplus\cdots(n)\cdots\oplus\frakg^\star) \\
  m:G(x,y,\ldots)
  & \mapsto &
  G(z,z,\ldots).
\end{eqnarray*}

The difference between the two products $\times_U$ and $\times_S$ was
discussed in~\cite[claim~\ref{I-DifferByPureDiff}]{Aarhus:I}. From that
discussion it follows that
\[ \mright^{xy}_z-\mright^{yx}_z  = m\circ D, \]
where $D$ is some differential operator acting on the variables $x$
and $y$.

We wish to study the integral with respect to $z$ of $mDG$; that
is, the integral of $DG$ on the diagonal $x=y$. We do this below by
changing coordinates to the parallel coordinate $\alpha=(x+y)/2$ and
the transverse coordinate $\beta=(x-y)/2$. In these coordinates the
diagonal is given by $\beta=0$, and $z$-integration can be replaced
by $\alpha$-integration at $\beta=0$. We show that this
$\alpha$-integral vanishes by showing that the divergence (with respect
to $\alpha$) of $D$ vanishes.

Let us turn to the details now. Let $\sigma G$ denote the (sum
of) uni-trivalent diagrams corresponding to $G$ via the isomorphism
$\sigma:\calA(\uparrow_x\uparrow_y\uparrow_E)\to\calB(\{x,y\}\udot
E)$.  By abuse of notation, we re-use the symbols $\mright^{xy}_z$
and $\mright^{yx}_z $ to denote the maps $\calB(\{x,y\}\udot
E)\to\calB(\{z\}\udot E)$ corresponding to the original
$\mright^{xy}_z$ and $\mright^{yx}_z $ via the isomorphism
$\sigma$ (really, $\sigma$ and its homonymic brother
$\sigma:\calA(\uparrow_z\uparrow_E)\to\calB(\{z\}\udot E)$).

Corollary~\ref{cor:BCH} implies that $\mright^{xy}_z$ and
$\mright^{yx}_z $ both
act on $\sigma G$ in the following manner:
\begin{itemize}
\item Glue some forest of non-trivial trees\footnote{
    That is, a forest in which each tree has at least one internal vertex
  }
  whose leaves are labeled $\partial_x$ and $\partial_y$ and whose roots
  are labeled $z$ to the $x$- and $y$-labeled legs of $\sigma G$,
  gluing only $\partial_x$'s to $x$'s and $\partial_y$'s to $y$'s, and
  making sure that all leaves get glued.
\item Relabel all remaining $x$- and $y$-labeled legs with $z$.
\end{itemize}

Gluing a forest of non-trivial trees does not touch the quadratic part of a
group-like element, and hence this description implies that $\mright^{xy}_z\sigma
G$ and $\mright^{yx}_z \sigma G$ have the same quadratic part. Hence if one is
integrable so is the other, and they can be integrated together under
the same formal integral sign, and we need to prove that
\[ \Gint(\mright^{xy}_z-\mright^{yx}_z )\sigma G \,dF = 0. \]
Using proposition~\ref{prop:ExpBCH}, we see that
$(\mright^{xy}_z-\mright^{yx}_z)(\sigma G)=\langle D_1,\sigma
G\rangle_{x,y}$, where $D_1$ is sum of forests of the general form
\[ D_1 = \cdots+\printname{xyForest}
	\setlength{\unitlength}{0.5\standardunitlength}
	\begin{array}{c}  \hspace{-1.7mm}
        	\raisebox{-8pt}{\input draws/xyForest.tex }
        	\hspace{-1.9mm}
	\end{array}
+\cdots. \]

We now change variables to $\alpha=(x+y)/2$, $\beta=(x-y)/2$,
$\partial_\alpha=\partial_x+\partial_y$, and
$\partial_\beta=\partial_x-\partial_y$ (see
exercises~\ref{ex:Reparametrization} and~\ref{ex:ReparametrizeGint}),
and rename the dummy integration variable $z$ to be $\alpha$. In the
new variables, $D_1$ becomes a sum of forests $D_2$ of the general form
\[ D_2 = \cdots+\printname{abForest}
	\setlength{\unitlength}{0.5\standardunitlength}
	\begin{array}{c}  \hspace{-1.7mm}
        	\raisebox{-8pt}{\input draws/abForest.tex }
        	\hspace{-1.9mm}
	\end{array}
+\cdots, \]
and the statement we need to prove is
\[ \Gint\langle D_2,\sigma G\rangle_{x,y}\,d\alpha d(F-\{z\}) = 0. \]
Notice that struts labeled $\partial_x$ or $\partial_y$ appear in
$D_1$ in a symmetric role, and hence $D_2$ has only struts labeled by
$\partial_\alpha$ and no struts labeled $\partial_\beta$. This implies
that if $D_3$ is the strutless part of $D_2$, then $\langle D_2,\sigma
G\rangle_{x,y}=D_3\apply(\sigma G$).  Using integration by parts (see
section~\ref{subsec:IBP}), we see that our job will be done one we can
prove that $\operatorname{div}_\alpha D_3=0$.

Each component of each forest in $D_3$ must have at least one
$\partial_\alpha$ leaf, for otherwise it would have only $\partial_\beta$
leaves and it would vanish by the $AS$ relation. Forests in which some tree
has more than one $\partial_\alpha$ leaf contribute nothing to
$\operatorname{div}_\alpha D_3$, as their
$\partial_\alpha$ order is higher than their $\alpha$ degree.
Thus we only care about forests of the form
\[ \printname{ContributingForest}
	\setlength{\unitlength}{0.5\standardunitlength}
	\begin{array}{c}  \hspace{-1.7mm}
        	\raisebox{-8pt}{\input draws/ContributingForest.tex }
        	\hspace{-1.9mm}
	\end{array}
, \]
in which each component carries exactly one $\partial_\alpha$.

The $\alpha$-divergence of such a forest is obtained by having the
$\alpha$-derivatives act on the coefficients --- namely, by summing over
all possible ways of connecting the $\partial_\alpha$'s at the top of the
picture to the $\alpha$'s at the bottom. The result is a sum $W$ of
diagrams like
\[ \printname{BetaWheels}
	\setlength{\unitlength}{0.5\standardunitlength}
	\begin{array}{c}  \hspace{-1.7mm}
        	\raisebox{-8pt}{\input draws/BetaWheels.tex }
        	\hspace{-1.9mm}
	\end{array}
, \]
that are disjoint unions of $\partial_\beta$-wheels.

The next thing to notice is that $\mright^{xy}_z-\mright^{yx}_z $ is odd under reversal of
the roles of $x$ and $y$, and hence under flipping the sign of $\beta$.
This means that $\partial_\beta$ appears an odd number of times in each
forest in $D_1$, $D_2$, and $D_3$, and thus in each union of wheels in $W$.
But this means that every term in $W$ contains a wheel with an odd number
of legs, and such wheels vanish by the $AS$ relation. Hence $W=0$ and thus
$\operatorname{div}_\alpha D_3=0$, as required.
\end{proof}

%% file: draws/LinkRelation.tex
\begingroup\makeatletter\ifx\SetFigFont\undefined
\def\x#1#2#3#4#5#6#7\relax{\def\x{#1#2#3#4#5#6}}%
\expandafter\x\fmtname xxxxxx\relax \def\y{splain}%
\ifx\x\y   
\gdef\SetFigFont#1#2#3{%
  \ifnum #1<17\tiny\else \ifnum #1<20\small\else
  \ifnum #1<24\normalsize\else \ifnum #1<29\large\else
  \ifnum #1<34\Large\else \ifnum #1<41\LARGE\else
     \huge\fi\fi\fi\fi\fi\fi
  \csname #3\endcsname}%
\else
\gdef\SetFigFont#1#2#3{\begingroup
  \count@#1\relax \ifnum 25<\count@\count@25\fi
  \def\x{\endgroup\@setsize\SetFigFont{#2pt}}%
  \expandafter\x
    \csname \romannumeral\the\count@ pt\expandafter\endcsname
    \csname @\romannumeral\the\count@ pt\endcsname
  \csname #3\endcsname}%
\fi
\fi\endgroup
{\renewcommand{\dashlinestretch}{30}
\begin{picture}(13611,1406)(0,-10)
\put(750,708){\ellipse{848}{848}}
\put(4351,708){\ellipse{848}{848}}
\put(7951,708){\ellipse{848}{848}}
\put(11551,708){\ellipse{848}{848}}
\path(2925,708)(3525,708)
\path(3405.000,678.000)(3525.000,708.000)(3405.000,738.000)
\path(150,1308)(450,1008)
\path(450,408)(150,108)
\path(1050,1008)(1350,1308)
\path(1050,408)(1350,108)
\path(1801,1158)(2701,1158)
\path(2251,1158)(2251,258)
\path(1801,258)(2701,258)
\path(3751,1308)(4051,1008)
\path(4051,408)(3751,108)
\path(4651,1008)(4951,1308)
\path(5402,1158)(6302,1158)
\path(5852,1158)(5852,258)
\path(5402,258)(6302,258)
\path(7351,1308)(7651,1008)
\path(7651,408)(7351,108)
\path(8251,1008)(8551,1308)
\path(9002,1158)(9902,1158)
\path(9452,1158)(9452,258)
\path(9002,258)(9902,258)
\path(10951,1308)(11251,1008)
\path(11251,408)(10951,108)
\path(11851,1008)(12151,1308)
\path(12602,1158)(13502,1158)
\path(13052,1158)(13052,258)
\path(12602,258)(13502,258)
\path(4650,408)	(4704.541,349.876)
	(4745.378,307.845)
	(4800.000,258.000)

\path(4800,258)	(4834.387,231.221)
	(4876.276,199.289)
	(4923.699,164.588)
	(4974.690,129.502)
	(5027.282,96.417)
	(5079.509,67.715)
	(5129.404,45.781)
	(5175.000,33.000)

\path(5175,33)	(5216.227,28.560)
	(5263.567,28.167)
	(5314.769,31.771)
	(5367.585,39.323)
	(5419.765,50.771)
	(5469.061,66.067)
	(5513.222,85.160)
	(5550.000,108.000)

\path(5550,108)	(5587.676,157.058)
	(5605.503,199.285)
	(5625.000,258.000)

\path(8250,408)	(8290.315,340.865)
	(8326.410,295.830)
	(8400.000,258.000)

\path(8400,258)	(8456.171,254.779)
	(8512.415,258.516)
	(8568.452,268.614)
	(8624.001,284.480)
	(8678.778,305.517)
	(8732.504,331.130)
	(8784.897,360.725)
	(8835.675,393.705)
	(8884.557,429.476)
	(8931.261,467.442)
	(8975.507,507.007)
	(9017.012,547.578)
	(9055.495,588.558)
	(9090.675,629.351)
	(9122.271,669.364)
	(9150.000,708.000)

\path(9150,708)	(9183.971,775.154)
	(9196.603,816.772)
	(9206.625,865.507)
	(9214.256,922.680)
	(9219.716,989.607)
	(9221.701,1027.141)
	(9223.224,1067.608)
	(9224.315,1111.173)
	(9225.000,1158.000)

\path(11850,408)	(11886.097,337.102)
	(11920.785,290.812)
	(11957.581,265.617)
	(12000.000,258.000)

\path(12000,258)	(12059.339,265.194)
	(12113.315,285.345)
	(12161.488,316.303)
	(12203.423,355.920)
	(12238.679,402.047)
	(12266.819,452.535)
	(12287.406,505.236)
	(12300.000,558.000)

\path(12300,558)	(12303.837,599.591)
	(12302.320,641.846)
	(12295.953,684.506)
	(12285.235,727.311)
	(12270.670,770.002)
	(12252.759,812.320)
	(12232.003,854.007)
	(12208.905,894.803)
	(12183.966,934.448)
	(12157.688,972.684)
	(12103.122,1043.893)
	(12049.221,1106.356)
	(12000.000,1158.000)

\path(12000,1158)	(11960.411,1193.821)
	(11911.485,1231.585)
	(11855.511,1269.221)
	(11794.778,1304.659)
	(11731.571,1335.826)
	(11668.181,1360.653)
	(11606.895,1377.068)
	(11550.000,1383.000)

\path(11550,1383)	(11509.112,1380.253)
	(11466.472,1371.574)
	(11420.760,1356.303)
	(11370.660,1333.781)
	(11314.852,1303.349)
	(11252.018,1264.348)
	(11180.841,1216.118)
	(11141.710,1188.336)
	(11100.000,1158.000)

\put(1350,183){\makebox(0,0)[lb]{\smash{{\mathmode{\scriptstyle x}}}}}
\put(1125,33){\makebox(0,0)[lb]{\smash{{\mathmode{\scriptstyle \ast}}}}}
\put(1350,1158){\makebox(0,0)[lb]{\smash{{\mathmode{\scriptstyle y}}}}}
\put(0,183){\makebox(0,0)[lb]{\smash{{\mathmode{\scriptstyle z}}}}}
\put(0,1158){\makebox(0,0)[lb]{\smash{{\mathmode{\scriptstyle x}}}}}
\put(1801,1233){\makebox(0,0)[lb]{\smash{{\mathmode{\scriptstyle x}}}}}
\put(2626,1233){\makebox(0,0)[lb]{\smash{{\mathmode{\scriptstyle z}}}}}
\put(1801,333){\makebox(0,0)[lb]{\smash{{\mathmode{\scriptstyle x}}}}}
\put(2626,333){\makebox(0,0)[lb]{\smash{{\mathmode{\scriptstyle y}}}}}
\put(4951,1158){\makebox(0,0)[lb]{\smash{{\mathmode{\scriptstyle y}}}}}
\put(3601,183){\makebox(0,0)[lb]{\smash{{\mathmode{\scriptstyle z}}}}}
\put(3601,1158){\makebox(0,0)[lb]{\smash{{\mathmode{\scriptstyle x}}}}}
\put(5402,1233){\makebox(0,0)[lb]{\smash{{\mathmode{\scriptstyle x}}}}}
\put(6227,1233){\makebox(0,0)[lb]{\smash{{\mathmode{\scriptstyle z}}}}}
\put(5402,333){\makebox(0,0)[lb]{\smash{{\mathmode{\scriptstyle x}}}}}
\put(6227,333){\makebox(0,0)[lb]{\smash{{\mathmode{\scriptstyle y}}}}}
\put(8551,1158){\makebox(0,0)[lb]{\smash{{\mathmode{\scriptstyle y}}}}}
\put(7201,183){\makebox(0,0)[lb]{\smash{{\mathmode{\scriptstyle z}}}}}
\put(7201,1158){\makebox(0,0)[lb]{\smash{{\mathmode{\scriptstyle x}}}}}
\put(9002,1233){\makebox(0,0)[lb]{\smash{{\mathmode{\scriptstyle x}}}}}
\put(9827,1233){\makebox(0,0)[lb]{\smash{{\mathmode{\scriptstyle z}}}}}
\put(9002,333){\makebox(0,0)[lb]{\smash{{\mathmode{\scriptstyle x}}}}}
\put(9827,333){\makebox(0,0)[lb]{\smash{{\mathmode{\scriptstyle y}}}}}
\put(12151,1158){\makebox(0,0)[lb]{\smash{{\mathmode{\scriptstyle y}}}}}
\put(10801,183){\makebox(0,0)[lb]{\smash{{\mathmode{\scriptstyle z}}}}}
\put(10801,1158){\makebox(0,0)[lb]{\smash{{\mathmode{\scriptstyle x}}}}}
\put(12602,1233){\makebox(0,0)[lb]{\smash{{\mathmode{\scriptstyle x}}}}}
\put(13427,1233){\makebox(0,0)[lb]{\smash{{\mathmode{\scriptstyle z}}}}}
\put(12602,333){\makebox(0,0)[lb]{\smash{{\mathmode{\scriptstyle x}}}}}
\put(13427,333){\makebox(0,0)[lb]{\smash{{\mathmode{\scriptstyle y}}}}}
\put(6675,633){\makebox(0,0)[lb]{\smash{{\mathmode{+}}}}}
\put(10275,633){\makebox(0,0)[lb]{\smash{{\mathmode{+}}}}}
\end{picture}
}

%% file: draws/LinkRelInA.tex
\begingroup\makeatletter\ifx\SetFigFont\undefined
\def\x#1#2#3#4#5#6#7\relax{\def\x{#1#2#3#4#5#6}}%
\expandafter\x\fmtname xxxxxx\relax \def\y{splain}%
\ifx\x\y   
\gdef\SetFigFont#1#2#3{%
  \ifnum #1<17\tiny\else \ifnum #1<20\small\else
  \ifnum #1<24\normalsize\else \ifnum #1<29\large\else
  \ifnum #1<34\Large\else \ifnum #1<41\LARGE\else
     \huge\fi\fi\fi\fi\fi\fi
  \csname #3\endcsname}%
\else
\gdef\SetFigFont#1#2#3{\begingroup
  \count@#1\relax \ifnum 25<\count@\count@25\fi
  \def\x{\endgroup\@setsize\SetFigFont{#2pt}}%
  \expandafter\x
    \csname \romannumeral\the\count@ pt\expandafter\endcsname
    \csname @\romannumeral\the\count@ pt\endcsname
  \csname #3\endcsname}%
\fi
\fi\endgroup
{\renewcommand{\dashlinestretch}{30}
\begin{picture}(8725,954)(0,-10)
\put(612,477){\ellipse{450}{450}}
\path(12,27)(12,927)
\path(42.000,807.000)(12.000,927.000)(-18.000,807.000)
\path(1212,27)(1212,927)
\path(1242.000,807.000)(1212.000,927.000)(1182.000,807.000)
\path(387,477)(12,477)
\path(462,627)(12,627)
\path(1212,477)(837,477)
\path(1212,627)(762,627)
\path(1212,327)(762,327)
\path(462,327)(12,327)
\put(87,27){\makebox(0,0)[lb]{\smash{{\mathmode{\scriptstyle x}}}}}
\put(987,27){\makebox(0,0)[lb]{\smash{{\mathmode{\scriptstyle y}}}}}
\put(2413,477){\ellipse{450}{450}}
\path(1813,27)(1813,927)
\path(1843.000,807.000)(1813.000,927.000)(1783.000,807.000)
\path(3013,27)(3013,927)
\path(3043.000,807.000)(3013.000,927.000)(2983.000,807.000)
\path(2188,477)(1813,477)
\path(2263,627)(1813,627)
\path(3013,477)(2638,477)
\path(3013,627)(2563,627)
\path(3013,327)(2563,327)
\path(2263,327)(1813,327)
\put(1888,27){\makebox(0,0)[lb]{\smash{{\mathmode{\scriptstyle x}}}}}
\put(2788,27){\makebox(0,0)[lb]{\smash{{\mathmode{\scriptstyle y}}}}}
\put(4513,477){\ellipse{450}{450}}
\path(3913,27)(3913,927)
\path(3943.000,807.000)(3913.000,927.000)(3883.000,807.000)
\path(5113,27)(5113,927)
\path(5143.000,807.000)(5113.000,927.000)(5083.000,807.000)
\path(4288,477)(3913,477)
\path(4363,627)(3913,627)
\path(5113,477)(4738,477)
\path(5113,627)(4663,627)
\path(5113,327)(4663,327)
\path(4363,327)(3913,327)
\put(3988,27){\makebox(0,0)[lb]{\smash{{\mathmode{\scriptstyle x}}}}}
\put(4888,27){\makebox(0,0)[lb]{\smash{{\mathmode{\scriptstyle y}}}}}
\put(6313,477){\ellipse{450}{450}}
\path(5713,27)(5713,927)
\path(5743.000,807.000)(5713.000,927.000)(5683.000,807.000)
\path(6913,27)(6913,927)
\path(6943.000,807.000)(6913.000,927.000)(6883.000,807.000)
\path(6088,477)(5713,477)
\path(6163,627)(5713,627)
\path(6913,477)(6538,477)
\path(6913,627)(6463,627)
\path(6913,327)(6463,327)
\path(6163,327)(5713,327)
\put(5788,27){\makebox(0,0)[lb]{\smash{{\mathmode{\scriptstyle x}}}}}
\put(6688,27){\makebox(0,0)[lb]{\smash{{\mathmode{\scriptstyle y}}}}}
\put(8113,477){\ellipse{450}{450}}
\path(7513,27)(7513,927)
\path(7543.000,807.000)(7513.000,927.000)(7483.000,807.000)
\path(8713,27)(8713,927)
\path(8743.000,807.000)(8713.000,927.000)(8683.000,807.000)
\path(7888,477)(7513,477)
\path(7963,627)(7513,627)
\path(8713,477)(8338,477)
\path(8713,627)(8263,627)
\path(8713,327)(8263,327)
\path(7963,327)(7513,327)
\put(7588,27){\makebox(0,0)[lb]{\smash{{\mathmode{\scriptstyle x}}}}}
\put(8488,27){\makebox(0,0)[lb]{\smash{{\mathmode{\scriptstyle y}}}}}
\put(1362,401){\makebox(0,0)[lb]{\smash{{\mathmode{-}}}}}
\path(612,702)	(603.060,749.141)
	(612.000,777.000)

\path(612,777)	(675.896,824.123)
	(715.802,837.609)
	(758.250,845.970)
	(801.166,850.397)
	(842.479,852.082)
	(912.000,852.000)

\path(912,852)	(959.094,848.191)
	(1020.105,835.594)
	(1058.026,825.450)
	(1102.063,812.449)
	(1153.095,796.373)
	(1212.000,777.000)

\path(2412,702)	(2405.122,749.471)
	(2412.000,777.000)

\path(2412,777)	(2476.661,834.424)
	(2519.245,850.259)
	(2562.000,852.000)

\path(2562,852)	(2611.265,830.547)
	(2653.793,790.744)
	(2687.924,744.069)
	(2712.000,702.000)

\path(2712,702)	(2726.721,657.439)
	(2733.312,603.167)
	(2735.075,542.370)
	(2735.310,478.234)
	(2737.319,413.943)
	(2744.403,352.683)
	(2759.863,297.640)
	(2787.000,252.000)

\path(2787,252)	(2818.355,224.870)
	(2862.791,204.499)
	(2925.582,189.128)
	(2965.508,182.768)
	(3012.000,177.000)

\path(4512,702)	(4503.060,749.141)
	(4512.000,777.000)

\path(4512,777)	(4573.203,829.269)
	(4610.753,847.260)
	(4651.069,859.695)
	(4692.751,866.481)
	(4734.400,867.523)
	(4774.616,862.728)
	(4812.000,852.000)

\path(4812,852)	(4854.543,827.074)
	(4891.380,784.781)
	(4926.027,719.848)
	(4943.628,677.243)
	(4962.000,627.000)

\path(6312,702)	(6303.060,749.141)
	(6312.000,777.000)

\path(6312,777)	(6372.444,831.318)
	(6409.329,851.101)
	(6449.044,865.159)
	(6490.378,872.884)
	(6532.122,873.670)
	(6573.066,866.911)
	(6612.000,852.000)

\path(6612,852)	(6674.759,804.765)
	(6718.335,732.379)
	(6734.028,684.008)
	(6746.243,626.054)
	(6755.421,557.417)
	(6759.008,518.750)
	(6762.000,477.000)

\path(8112,702)	(8103.060,749.141)
	(8112.000,777.000)

\path(8112,777)	(8172.113,832.440)
	(8208.710,853.206)
	(8248.163,868.151)
	(8289.345,876.391)
	(8331.130,877.036)
	(8372.391,869.202)
	(8412.000,852.000)

\path(8412,852)	(8458.335,819.889)
	(8495.909,781.279)
	(8525.163,734.632)
	(8546.535,678.409)
	(8560.465,611.072)
	(8564.777,572.756)
	(8567.393,531.084)
	(8568.369,485.865)
	(8567.758,436.906)
	(8565.617,384.015)
	(8562.000,327.000)

\put(3312,402){\makebox(0,0)[lb]{\smash{{\mathmode{=}}}}}
\put(5262,402){\makebox(0,0)[lb]{\smash{{\mathmode{+}}}}}
\put(7062,402){\makebox(0,0)[lb]{\smash{{\mathmode{+}}}}}
\end{picture}
}

%% file: draws/StillLinkRel.tex
\begingroup\makeatletter\ifx\SetFigFont\undefined
\def\x#1#2#3#4#5#6#7\relax{\def\x{#1#2#3#4#5#6}}%
\expandafter\x\fmtname xxxxxx\relax \def\y{splain}%
\ifx\x\y   
\gdef\SetFigFont#1#2#3{%
  \ifnum #1<17\tiny\else \ifnum #1<20\small\else
  \ifnum #1<24\normalsize\else \ifnum #1<29\large\else
  \ifnum #1<34\Large\else \ifnum #1<41\LARGE\else
     \huge\fi\fi\fi\fi\fi\fi
  \csname #3\endcsname}%
\else
\gdef\SetFigFont#1#2#3{\begingroup
  \count@#1\relax \ifnum 25<\count@\count@25\fi
  \def\x{\endgroup\@setsize\SetFigFont{#2pt}}%
  \expandafter\x
    \csname \romannumeral\the\count@ pt\expandafter\endcsname
    \csname @\romannumeral\the\count@ pt\endcsname
  \csname #3\endcsname}%
\fi
\fi\endgroup
{\renewcommand{\dashlinestretch}{30}
\begin{picture}(12786,1008)(0,-10)
\put(3762,399){\makebox(0,0)[lb]{\smash{{\mathmode{+}}}}}
\put(8263,399){\makebox(0,0)[lb]{\smash{{\mathmode{+}}}}}
\put(612,474){\ellipse{450}{450}}
\put(2714,474){\ellipse{450}{450}}
\put(4813,474){\ellipse{450}{450}}
\put(7214,474){\ellipse{450}{450}}
\put(9314,474){\ellipse{450}{450}}
\put(11715,474){\ellipse{450}{450}}
\path(12,24)(12,924)
\path(42.000,804.000)(12.000,924.000)(-18.000,804.000)
\path(387,474)(12,474)
\path(462,624)(12,624)
\path(1212,474)(837,474)
\path(1212,624)(762,624)
\path(1212,324)(762,324)
\path(462,324)(12,324)
\path(2114,24)(2114,924)
\path(2144.000,804.000)(2114.000,924.000)(2084.000,804.000)
\path(2489,474)(2114,474)
\path(2564,624)(2114,624)
\path(3314,474)(2939,474)
\path(3314,624)(2864,624)
\path(3314,324)(2864,324)
\path(2564,324)(2114,324)
\path(1210,624)(1285,549)(1210,474)
\path(1285,549)(1510,549)
\path(1210,324)(1510,324)
\path(3310,624)(3385,549)(3310,474)
\path(3385,549)(3610,549)
\path(3310,324)(3610,324)
\path(4213,24)(4213,924)
\path(4243.000,804.000)(4213.000,924.000)(4183.000,804.000)
\path(4588,474)(4213,474)
\path(4663,624)(4213,624)
\path(5413,474)(5038,474)
\path(5413,624)(4963,624)
\path(5413,324)(4963,324)
\path(4663,324)(4213,324)
\path(5410,624)(5485,549)
\path(5410,474)(5485,549)(5710,549)
\path(5410,324)(5710,324)
\path(6614,24)(6614,924)
\path(6644.000,804.000)(6614.000,924.000)(6584.000,804.000)
\path(6989,474)(6614,474)
\path(7064,624)(6614,624)
\path(7814,474)(7439,474)
\path(7814,624)(7364,624)
\path(7814,324)(7364,324)
\path(7064,324)(6614,324)
\path(7811,624)(7886,549)
\path(7811,474)(7886,549)(8111,549)
\path(7811,324)(8111,324)
\path(10510,474)(10810,474)
\path(10630.000,504.000)(10510.000,474.000)(10630.000,444.000)
\path(10690.000,444.000)(10810.000,474.000)(10690.000,504.000)
\path(8714,24)(8714,924)
\path(8744.000,804.000)(8714.000,924.000)(8684.000,804.000)
\path(9089,474)(8714,474)
\path(9164,624)(8714,624)
\path(9914,474)(9539,474)
\path(9914,624)(9464,624)
\path(9914,324)(9464,324)
\path(9164,324)(8714,324)
\path(9911,624)(9986,549)
\path(9911,474)(9986,549)(10211,549)
\path(9911,324)(10211,324)
\path(11115,24)(11115,924)
\path(11145.000,804.000)(11115.000,924.000)(11085.000,804.000)
\path(11490,474)(11115,474)
\path(11565,624)(11115,624)
\path(12160,474)(11940,474)
\path(12160,624)(11865,624)
\path(12460,324)(11865,324)
\path(11565,324)(11115,324)
\path(12160,474)(12235,549)(12460,549)
\path(12160,624)(12235,549)
\path(611,699)	(602.060,746.141)
	(611.000,774.000)

\path(611,774)	(672.203,826.269)
	(709.753,844.260)
	(750.069,856.695)
	(791.751,863.481)
	(833.400,864.523)
	(873.616,859.728)
	(911.000,849.000)

\path(911,849)	(953.543,824.074)
	(990.380,781.781)
	(1025.027,716.848)
	(1042.628,674.243)
	(1061.000,624.000)

\path(2713,699)	(2704.060,746.141)
	(2713.000,774.000)

\path(2713,774)	(2773.444,828.318)
	(2810.329,848.101)
	(2850.044,862.159)
	(2891.378,869.884)
	(2933.122,870.670)
	(2974.066,863.911)
	(3013.000,849.000)

\path(3013,849)	(3075.759,801.765)
	(3119.335,729.379)
	(3135.028,681.008)
	(3147.243,623.054)
	(3156.421,554.417)
	(3160.008,515.750)
	(3163.000,474.000)

\path(4812,699)	(4803.060,746.141)
	(4812.000,774.000)

\path(4812,774)	(4872.113,829.440)
	(4908.710,850.206)
	(4948.163,865.151)
	(4989.345,873.391)
	(5031.130,874.036)
	(5072.391,866.202)
	(5112.000,849.000)

\path(5112,849)	(5158.335,816.889)
	(5195.909,778.279)
	(5225.163,731.632)
	(5246.535,675.409)
	(5260.465,608.072)
	(5264.777,569.756)
	(5267.393,528.084)
	(5268.369,482.865)
	(5267.758,433.906)
	(5265.617,381.015)
	(5262.000,324.000)

\path(7213,699)	(7202.650,746.029)
	(7213.000,774.000)

\path(7213,774)	(7281.211,836.210)
	(7319.297,861.796)
	(7359.569,883.674)
	(7401.672,901.821)
	(7445.256,916.217)
	(7489.968,926.838)
	(7535.455,933.664)
	(7581.365,936.672)
	(7627.347,935.840)
	(7673.047,931.148)
	(7718.114,922.572)
	(7762.195,910.091)
	(7804.938,893.683)
	(7845.990,873.327)
	(7885.000,849.000)

\path(7885,849)	(7927.779,806.681)
	(7953.366,747.263)
	(7960.263,708.943)
	(7963.521,663.712)
	(7963.360,610.691)
	(7960.000,549.000)

\path(9313,699)	(9302.800,746.036)
	(9313.000,774.000)

\path(9313,774)	(9372.991,829.456)
	(9442.204,872.070)
	(9479.481,888.510)
	(9518.118,901.680)
	(9557.801,911.558)
	(9598.214,918.124)
	(9639.042,921.358)
	(9679.970,921.240)
	(9720.683,917.749)
	(9760.867,910.866)
	(9800.205,900.570)
	(9838.383,886.840)
	(9910.000,849.000)

\path(9910,849)	(9952.223,815.911)
	(9986.861,776.603)
	(10014.353,729.536)
	(10035.138,673.174)
	(10049.655,605.977)
	(10054.701,567.835)
	(10058.345,526.408)
	(10060.642,481.503)
	(10061.647,432.928)
	(10061.415,380.491)
	(10060.000,324.000)

\path(11710,699)	(11703.122,746.471)
	(11710.000,774.000)

\path(11710,774)	(11779.405,825.701)
	(11822.582,840.821)
	(11860.000,849.000)

\path(11860,849)	(11931.339,852.481)
	(11974.003,849.214)
	(12023.088,842.314)
	(12079.911,831.562)
	(12145.792,816.739)
	(12182.541,807.732)
	(12222.049,797.625)
	(12264.480,786.390)
	(12310.000,774.000)

\put(87,24){\makebox(0,0)[lb]{\smash{{\mathmode{\scriptstyle x}}}}}
\put(1663,399){\makebox(0,0)[lb]{\smash{{\mathmode{+}}}}}
\put(2189,24){\makebox(0,0)[lb]{\smash{{\mathmode{\scriptstyle x}}}}}
\put(4288,24){\makebox(0,0)[lb]{\smash{{\mathmode{\scriptstyle x}}}}}
\put(6689,24){\makebox(0,0)[lb]{\smash{{\mathmode{\scriptstyle x}}}}}
\put(6010,399){\makebox(0,0)[lb]{\smash{{\mathmode{=}}}}}
\put(8789,24){\makebox(0,0)[lb]{\smash{{\mathmode{\scriptstyle x}}}}}
\put(11190,24){\makebox(0,0)[lb]{\smash{{\mathmode{\scriptstyle x}}}}}
\put(12310,849){\makebox(0,0)[lb]{\smash{{\mathmode{\scriptstyle \ast y}}}}}
\put(12385,399){\makebox(0,0)[lb]{\smash{{\mathmode{\scriptstyle y}}}}}
\put(10210,399){\makebox(0,0)[lb]{\smash{{\mathmode{\scriptstyle y}}}}}
\put(7960,399){\makebox(0,0)[lb]{\smash{{\mathmode{\scriptstyle y}}}}}
\put(5560,399){\makebox(0,0)[lb]{\smash{{\mathmode{\scriptstyle y}}}}}
\put(3460,399){\makebox(0,0)[lb]{\smash{{\mathmode{\scriptstyle y}}}}}
\put(1360,399){\makebox(0,0)[lb]{\smash{{\mathmode{\scriptstyle y}}}}}
\end{picture}
}

%% file: draws/x.tex
%
\begingroup\makeatletter\ifx\SetFigFont\undefined
\def\x#1#2#3#4#5#6#7\relax{\def\x{#1#2#3#4#5#6}}%
\expandafter\x\fmtname xxxxxx\relax \def\y{splain}%
\ifx\x\y   
\gdef\SetFigFont#1#2#3{%
  \ifnum #1<17\tiny\else \ifnum #1<20\small\else
  \ifnum #1<24\normalsize\else \ifnum #1<29\large\else
  \ifnum #1<34\Large\else \ifnum #1<41\LARGE\else
     \huge\fi\fi\fi\fi\fi\fi
  \csname #3\endcsname}%
\else
\gdef\SetFigFont#1#2#3{\begingroup
  \count@#1\relax \ifnum 25<\count@\count@25\fi
  \def\x{\endgroup\@setsize\SetFigFont{#2pt}}%
  \expandafter\x
    \csname \romannumeral\the\count@ pt\expandafter\endcsname
    \csname @\romannumeral\the\count@ pt\endcsname
  \csname #3\endcsname}%
\fi
\fi\endgroup
{\renewcommand{\dashlinestretch}{30}
\begin{picture}(288,1425)(0,-10)
\path(75,1227)(75,27)
\put(0,1302){\makebox(0,0)[lb]{\smash{{\mathmode{\x}}}}}
\put(150,27){\makebox(0,0)[lb]{\smash{{\mathmode{\z}}}}}
\end{picture}
}

%% file: draws/y.tex
%
\begingroup\makeatletter\ifx\SetFigFont\undefined
\def\x#1#2#3#4#5#6#7\relax{\def\x{#1#2#3#4#5#6}}%
\expandafter\x\fmtname xxxxxx\relax \def\y{splain}%
\ifx\x\y   
\gdef\SetFigFont#1#2#3{%
  \ifnum #1<17\tiny\else \ifnum #1<20\small\else
  \ifnum #1<24\normalsize\else \ifnum #1<29\large\else
  \ifnum #1<34\Large\else \ifnum #1<41\LARGE\else
     \huge\fi\fi\fi\fi\fi\fi
  \csname #3\endcsname}%
\else
\gdef\SetFigFont#1#2#3{\begingroup
  \count@#1\relax \ifnum 25<\count@\count@25\fi
  \def\x{\endgroup\@setsize\SetFigFont{#2pt}}%
  \expandafter\x
    \csname \romannumeral\the\count@ pt\expandafter\endcsname
    \csname @\romannumeral\the\count@ pt\endcsname
  \csname #3\endcsname}%
\fi
\fi\endgroup
{\renewcommand{\dashlinestretch}{30}
\begin{picture}(288,1461)(0,-10)
\path(75,1227)(75,27)
\put(0,1302){\makebox(0,0)[lb]{\smash{{\mathmode{\y}}}}}
\put(150,27){\makebox(0,0)[lb]{\smash{{\mathmode{\z}}}}}
\end{picture}
}

%% file: draws/xy.tex
%
\begingroup\makeatletter\ifx\SetFigFont\undefined
\def\x#1#2#3#4#5#6#7\relax{\def\x{#1#2#3#4#5#6}}%
\expandafter\x\fmtname xxxxxx\relax \def\y{splain}%
\ifx\x\y   
\gdef\SetFigFont#1#2#3{%
  \ifnum #1<17\tiny\else \ifnum #1<20\small\else
  \ifnum #1<24\normalsize\else \ifnum #1<29\large\else
  \ifnum #1<34\Large\else \ifnum #1<41\LARGE\else
     \huge\fi\fi\fi\fi\fi\fi
  \csname #3\endcsname}%
\else
\gdef\SetFigFont#1#2#3{\begingroup
  \count@#1\relax \ifnum 25<\count@\count@25\fi
  \def\x{\endgroup\@setsize\SetFigFont{#2pt}}%
  \expandafter\x
    \csname \romannumeral\the\count@ pt\expandafter\endcsname
    \csname @\romannumeral\the\count@ pt\endcsname
  \csname #3\endcsname}%
\fi
\fi\endgroup
{\renewcommand{\dashlinestretch}{30}
\begin{picture}(738,1461)(0,-10)
\path(75,1227)(375,927)
\path(675,1227)(375,927)
\path(375,927)(375,27)
\put(0,1302){\makebox(0,0)[lb]{\smash{{\mathmode{\x}}}}}
\put(600,1302){\makebox(0,0)[lb]{\smash{{\mathmode{\y}}}}}
\put(450,27){\makebox(0,0)[lb]{\smash{{\mathmode{\z}}}}}
\end{picture}
}

%% file: draws/xxy.tex
%
\begingroup\makeatletter\ifx\SetFigFont\undefined
\def\x#1#2#3#4#5#6#7\relax{\def\x{#1#2#3#4#5#6}}%
\expandafter\x\fmtname xxxxxx\relax \def\y{splain}%
\ifx\x\y   
\gdef\SetFigFont#1#2#3{%
  \ifnum #1<17\tiny\else \ifnum #1<20\small\else
  \ifnum #1<24\normalsize\else \ifnum #1<29\large\else
  \ifnum #1<34\Large\else \ifnum #1<41\LARGE\else
     \huge\fi\fi\fi\fi\fi\fi
  \csname #3\endcsname}%
\else
\gdef\SetFigFont#1#2#3{\begingroup
  \count@#1\relax \ifnum 25<\count@\count@25\fi
  \def\x{\endgroup\@setsize\SetFigFont{#2pt}}%
  \expandafter\x
    \csname \romannumeral\the\count@ pt\expandafter\endcsname
    \csname @\romannumeral\the\count@ pt\endcsname
  \csname #3\endcsname}%
\fi
\fi\endgroup
{\renewcommand{\dashlinestretch}{30}
\begin{picture}(1338,1461)(0,-10)
\path(1275,1227)(675,627)
\path(975,927)(675,1227)
\path(75,1227)(675,627)
\path(675,627)(675,27)
\put(0,1302){\makebox(0,0)[lb]{\smash{{\mathmode{\x}}}}}
\put(600,1302){\makebox(0,0)[lb]{\smash{{\mathmode{\x}}}}}
\put(1200,1302){\makebox(0,0)[lb]{\smash{{\mathmode{\y}}}}}
\put(750,27){\makebox(0,0)[lb]{\smash{{\mathmode{\z}}}}}
\end{picture}
}

%% file: draws/yxy.tex
%
\begingroup\makeatletter\ifx\SetFigFont\undefined
\def\x#1#2#3#4#5#6#7\relax{\def\x{#1#2#3#4#5#6}}%
\expandafter\x\fmtname xxxxxx\relax \def\y{splain}%
\ifx\x\y   
\gdef\SetFigFont#1#2#3{%
  \ifnum #1<17\tiny\else \ifnum #1<20\small\else
  \ifnum #1<24\normalsize\else \ifnum #1<29\large\else
  \ifnum #1<34\Large\else \ifnum #1<41\LARGE\else
     \huge\fi\fi\fi\fi\fi\fi
  \csname #3\endcsname}%
\else
\gdef\SetFigFont#1#2#3{\begingroup
  \count@#1\relax \ifnum 25<\count@\count@25\fi
  \def\x{\endgroup\@setsize\SetFigFont{#2pt}}%
  \expandafter\x
    \csname \romannumeral\the\count@ pt\expandafter\endcsname
    \csname @\romannumeral\the\count@ pt\endcsname
  \csname #3\endcsname}%
\fi
\fi\endgroup
{\renewcommand{\dashlinestretch}{30}
\begin{picture}(1338,1461)(0,-10)
\path(1275,1227)(675,627)
\path(975,927)(675,1227)
\path(75,1227)(675,627)
\path(675,627)(675,27)
\put(0,1302){\makebox(0,0)[lb]{\smash{{\mathmode{\y}}}}}
\put(600,1302){\makebox(0,0)[lb]{\smash{{\mathmode{\x}}}}}
\put(1200,1302){\makebox(0,0)[lb]{\smash{{\mathmode{\y}}}}}
\put(750,27){\makebox(0,0)[lb]{\smash{{\mathmode{\z}}}}}
\end{picture}
}

%% file: draws/xyxy.tex
%
\begingroup\makeatletter\ifx\SetFigFont\undefined
\def\x#1#2#3#4#5#6#7\relax{\def\x{#1#2#3#4#5#6}}%
\expandafter\x\fmtname xxxxxx\relax \def\y{splain}%
\ifx\x\y   
\gdef\SetFigFont#1#2#3{%
  \ifnum #1<17\tiny\else \ifnum #1<20\small\else
  \ifnum #1<24\normalsize\else \ifnum #1<29\large\else
  \ifnum #1<34\Large\else \ifnum #1<41\LARGE\else
     \huge\fi\fi\fi\fi\fi\fi
  \csname #3\endcsname}%
\else
\gdef\SetFigFont#1#2#3{\begingroup
  \count@#1\relax \ifnum 25<\count@\count@25\fi
  \def\x{\endgroup\@setsize\SetFigFont{#2pt}}%
  \expandafter\x
    \csname \romannumeral\the\count@ pt\expandafter\endcsname
    \csname @\romannumeral\the\count@ pt\endcsname
  \csname #3\endcsname}%
\fi
\fi\endgroup
{\renewcommand{\dashlinestretch}{30}
\begin{picture}(1938,1461)(0,-10)
\path(75,1227)(975,327)
\path(675,1227)(1275,627)
\path(1275,1227)(1575,927)
\path(1875,1227)(975,327)
\path(975,327)(975,27)
\put(0,1302){\makebox(0,0)[lb]{\smash{{\mathmode{\x}}}}}
\put(600,1302){\makebox(0,0)[lb]{\smash{{\mathmode{\y}}}}}
\put(1200,1302){\makebox(0,0)[lb]{\smash{{\mathmode{\x}}}}}
\put(1800,1302){\makebox(0,0)[lb]{\smash{{\mathmode{\y}}}}}
\put(1050,27){\makebox(0,0)[lb]{\smash{{\mathmode{\z}}}}}
\end{picture}
}

%% file: draws/PlantedForest.tex
%
\begingroup\makeatletter\ifx\SetFigFont\undefined
\def\x#1#2#3#4#5#6#7\relax{\def\x{#1#2#3#4#5#6}}%
\expandafter\x\fmtname xxxxxx\relax \def\y{splain}%
\ifx\x\y   
\gdef\SetFigFont#1#2#3{%
  \ifnum #1<17\tiny\else \ifnum #1<20\small\else
  \ifnum #1<24\normalsize\else \ifnum #1<29\large\else
  \ifnum #1<34\Large\else \ifnum #1<41\LARGE\else
     \huge\fi\fi\fi\fi\fi\fi
  \csname #3\endcsname}%
\else
\gdef\SetFigFont#1#2#3{\begingroup
  \count@#1\relax \ifnum 25<\count@\count@25\fi
  \def\x{\endgroup\@setsize\SetFigFont{#2pt}}%
  \expandafter\x
    \csname \romannumeral\the\count@ pt\expandafter\endcsname
    \csname @\romannumeral\the\count@ pt\endcsname
  \csname #3\endcsname}%
\fi
\fi\endgroup
{\renewcommand{\dashlinestretch}{30}
\begin{picture}(6024,1446)(0,-10)
\path(1512,1212)(1812,912)(2112,1212)
\path(612,912)(1212,312)(1812,912)
\path(1212,312)(1212,12)
\path(2712,1212)(3012,912)(3312,1212)
\path(3012,912)(3012,12)
\path(3912,1212)(4812,312)
\path(4512,1212)(5112,612)
\path(5112,1212)(5412,912)
\path(5712,1212)(4812,312)(4812,12)
\path(12,12)(6012,12)
\path(5892.000,-18.000)(6012.000,12.000)(5892.000,42.000)
\path(312,1212)(612,912)(912,1212)
\put(237,1287){\makebox(0,0)[lb]{\smash{{\mathmode{\x}}}}}
\put(837,1287){\makebox(0,0)[lb]{\smash{{\mathmode{\y}}}}}
\put(1437,1287){\makebox(0,0)[lb]{\smash{{\mathmode{\y}}}}}
\put(2037,1287){\makebox(0,0)[lb]{\smash{{\mathmode{\x}}}}}
\put(2637,1287){\makebox(0,0)[lb]{\smash{{\mathmode{\x}}}}}
\put(3237,1287){\makebox(0,0)[lb]{\smash{{\mathmode{\y}}}}}
\put(3837,1287){\makebox(0,0)[lb]{\smash{{\mathmode{\x}}}}}
\put(4437,1287){\makebox(0,0)[lb]{\smash{{\mathmode{\x}}}}}
\put(5037,1287){\makebox(0,0)[lb]{\smash{{\mathmode{\x}}}}}
\put(5637,1287){\makebox(0,0)[lb]{\smash{{\mathmode{\y}}}}}
\end{picture}
}

%% file: draws/xi.tex
%
\begingroup\makeatletter\ifx\SetFigFont\undefined
\def\x#1#2#3#4#5#6#7\relax{\def\x{#1#2#3#4#5#6}}%
\expandafter\x\fmtname xxxxxx\relax \def\y{splain}%
\ifx\x\y   
\gdef\SetFigFont#1#2#3{%
  \ifnum #1<17\tiny\else \ifnum #1<20\small\else
  \ifnum #1<24\normalsize\else \ifnum #1<29\large\else
  \ifnum #1<34\Large\else \ifnum #1<41\LARGE\else
     \huge\fi\fi\fi\fi\fi\fi
  \csname #3\endcsname}%
\else
\gdef\SetFigFont#1#2#3{\begingroup
  \count@#1\relax \ifnum 25<\count@\count@25\fi
  \def\x{\endgroup\@setsize\SetFigFont{#2pt}}%
  \expandafter\x
    \csname \romannumeral\the\count@ pt\expandafter\endcsname
    \csname @\romannumeral\the\count@ pt\endcsname
  \csname #3\endcsname}%
\fi
\fi\endgroup
{\renewcommand{\dashlinestretch}{30}
\begin{picture}(624,510)(0,-10)
\path(12,12)(612,12)
\path(492.000,-18.000)(612.000,12.000)(492.000,42.000)
\path(312,312)(312,12)
\put(237,387){\makebox(0,0)[lb]{\smash{{\mathmode{\x}}}}}
\end{picture}
}

%% file: draws/eta.tex
%
\begingroup\makeatletter\ifx\SetFigFont\undefined
\def\x#1#2#3#4#5#6#7\relax{\def\x{#1#2#3#4#5#6}}%
\expandafter\x\fmtname xxxxxx\relax \def\y{splain}%
\ifx\x\y   
\gdef\SetFigFont#1#2#3{%
  \ifnum #1<17\tiny\else \ifnum #1<20\small\else
  \ifnum #1<24\normalsize\else \ifnum #1<29\large\else
  \ifnum #1<34\Large\else \ifnum #1<41\LARGE\else
     \huge\fi\fi\fi\fi\fi\fi
  \csname #3\endcsname}%
\else
\gdef\SetFigFont#1#2#3{\begingroup
  \count@#1\relax \ifnum 25<\count@\count@25\fi
  \def\x{\endgroup\@setsize\SetFigFont{#2pt}}%
  \expandafter\x
    \csname \romannumeral\the\count@ pt\expandafter\endcsname
    \csname @\romannumeral\the\count@ pt\endcsname
  \csname #3\endcsname}%
\fi
\fi\endgroup
{\renewcommand{\dashlinestretch}{30}
\begin{picture}(624,546)(0,-10)
\path(12,12)(612,12)
\path(492.000,-18.000)(612.000,12.000)(492.000,42.000)
\path(312,312)(312,12)
\put(237,387){\makebox(0,0)[lb]{\smash{{\mathmode{\y}}}}}
\end{picture}
}

%% file: draws/xyForest.tex
\begingroup\makeatletter\ifx\SetFigFont\undefined
\def\x#1#2#3#4#5#6#7\relax{\def\x{#1#2#3#4#5#6}}%
\expandafter\x\fmtname xxxxxx\relax \def\y{splain}%
\ifx\x\y   
\gdef\SetFigFont#1#2#3{%
  \ifnum #1<17\tiny\else \ifnum #1<20\small\else
  \ifnum #1<24\normalsize\else \ifnum #1<29\large\else
  \ifnum #1<34\Large\else \ifnum #1<41\LARGE\else
     \huge\fi\fi\fi\fi\fi\fi
  \csname #3\endcsname}%
\else
\gdef\SetFigFont#1#2#3{\begingroup
  \count@#1\relax \ifnum 25<\count@\count@25\fi
  \def\x{\endgroup\@setsize\SetFigFont{#2pt}}%
  \expandafter\x
    \csname \romannumeral\the\count@ pt\expandafter\endcsname
    \csname @\romannumeral\the\count@ pt\endcsname
  \csname #3\endcsname}%
\fi
\fi\endgroup
{\renewcommand{\dashlinestretch}{30}
\begin{picture}(5615,1612)(0,-10)
\path(1274,1378)(1874,778)
\path(1874,1378)(2174,1078)
\path(75,1378)(75,178)
\path(675,1377)(675,177)
\path(3075,1377)(3975,477)
\path(3675,1377)(4275,777)
\path(4275,1377)(4575,1077)
\path(3974,478)(4875,1377)
\path(2474,1378)(1874,778)(1874,178)
\path(3974,478)(3974,178)
\put(2399,1453){\makebox(0,0)[lb]{\smash{{\mathmode{\scriptstyle \partial_y}}}}}
\put(1199,1453){\makebox(0,0)[lb]{\smash{{\mathmode{\scriptstyle \partial_x}}}}}
\put(1799,1453){\makebox(0,0)[lb]{\smash{{\mathmode{\scriptstyle \partial_x}}}}}
\put(0,28){\makebox(0,0)[lb]{\smash{{\mathmode{\scriptstyle z}}}}}
\put(600,27){\makebox(0,0)[lb]{\smash{{\mathmode{\scriptstyle z}}}}}
\put(4800,1452){\makebox(0,0)[lb]{\smash{{\mathmode{\scriptstyle \partial_y}}}}}
\put(3000,1452){\makebox(0,0)[lb]{\smash{{\mathmode{\scriptstyle \partial_x}}}}}
\put(4199,1453){\makebox(0,0)[lb]{\smash{{\mathmode{\scriptstyle \partial_x}}}}}
\put(3599,1453){\makebox(0,0)[lb]{\smash{{\mathmode{\scriptstyle \partial_y}}}}}
\put(0,1453){\makebox(0,0)[lb]{\smash{{\mathmode{\scriptstyle \partial_x}}}}}
\put(600,1452){\makebox(0,0)[lb]{\smash{{\mathmode{\scriptstyle \partial_y}}}}}
\put(1799,28){\makebox(0,0)[lb]{\smash{{\mathmode{\scriptstyle z}}}}}
\put(3900,27){\makebox(0,0)[lb]{\smash{{\mathmode{\scriptstyle z}}}}}
\end{picture}
}

%% file: draws/abForest.tex
%
\begingroup\makeatletter\ifx\SetFigFont\undefined
\def\x#1#2#3#4#5#6#7\relax{\def\x{#1#2#3#4#5#6}}%
\expandafter\x\fmtname xxxxxx\relax \def\y{splain}%
\ifx\x\y   
\gdef\SetFigFont#1#2#3{%
  \ifnum #1<17\tiny\else \ifnum #1<20\small\else
  \ifnum #1<24\normalsize\else \ifnum #1<29\large\else
  \ifnum #1<34\Large\else \ifnum #1<41\LARGE\else
     \huge\fi\fi\fi\fi\fi\fi
  \csname #3\endcsname}%
\else
\gdef\SetFigFont#1#2#3{\begingroup
  \count@#1\relax \ifnum 25<\count@\count@25\fi
  \def\x{\endgroup\@setsize\SetFigFont{#2pt}}%
  \expandafter\x
    \csname \romannumeral\the\count@ pt\expandafter\endcsname
    \csname @\romannumeral\the\count@ pt\endcsname
  \csname #3\endcsname}%
\fi
\fi\endgroup
{\renewcommand{\dashlinestretch}{30}
\begin{picture}(4107,1620)(0,-10)
\path(76,1386)(76,186)
\path(676,1385)(676,185)
\path(1276,1385)(2176,485)
\path(1876,1385)(2476,785)
\path(2476,1385)(2776,1085)
\path(2175,486)(3076,1385)
\path(2175,486)(2175,186)
\put(1200,1461){\makebox(0,0)[lb]{\smash{{\mathmode{\scriptstyle \partial_\beta}}}}}
\put(0,1461){\makebox(0,0)[lb]{\smash{{\mathmode{\scriptstyle \partial_\alpha}}}}}
\put(600,1461){\makebox(0,0)[lb]{\smash{{\mathmode{\scriptstyle \partial_\alpha}}}}}
\put(3000,1461){\makebox(0,0)[lb]{\smash{{\mathmode{\scriptstyle \partial_\beta}}}}}
\put(1800,1461){\makebox(0,0)[lb]{\smash{{\mathmode{\scriptstyle \partial_\beta}}}}}
\put(2400,1461){\makebox(0,0)[lb]{\smash{{\mathmode{\scriptstyle \partial_\alpha}}}}}
\put(0,36){\makebox(0,0)[lb]{\smash{{\mathmode{\scriptstyle \alpha}}}}}
\put(600,36){\makebox(0,0)[lb]{\smash{{\mathmode{\scriptstyle \alpha}}}}}
\put(2100,36){\makebox(0,0)[lb]{\smash{{\mathmode{\scriptstyle \alpha}}}}}
\end{picture}
}

%% file: draws/ContributingForest.tex
\begingroup\makeatletter\ifx\SetFigFont\undefined
\def\x#1#2#3#4#5#6#7\relax{\def\x{#1#2#3#4#5#6}}%
\expandafter\x\fmtname xxxxxx\relax \def\y{splain}%
\ifx\x\y   
\gdef\SetFigFont#1#2#3{%
  \ifnum #1<17\tiny\else \ifnum #1<20\small\else
  \ifnum #1<24\normalsize\else \ifnum #1<29\large\else
  \ifnum #1<34\Large\else \ifnum #1<41\LARGE\else
     \huge\fi\fi\fi\fi\fi\fi
  \csname #3\endcsname}%
\else
\gdef\SetFigFont#1#2#3{\begingroup
  \count@#1\relax \ifnum 25<\count@\count@25\fi
  \def\x{\endgroup\@setsize\SetFigFont{#2pt}}%
  \expandafter\x
    \csname \romannumeral\the\count@ pt\expandafter\endcsname
    \csname @\romannumeral\the\count@ pt\endcsname
  \csname #3\endcsname}%
\fi
\fi\endgroup
{\renewcommand{\dashlinestretch}{30}
\begin{picture}(5984,1620)(0,-10)
\path(750,1386)(1050,1086)
\path(150,1386)(750,786)
\path(1350,1386)(750,786)(750,186)
\path(1950,1386)(2850,486)
\path(2550,1386)(3150,786)
\path(3150,1386)(3450,1086)
\path(3750,1386)(2850,486)
\path(4350,1386)(4650,1086)
\path(4650,1086)(4650,186)
\path(4950,1386)(4650,1086)
\drawline(2850,486)(2850,486)
\path(2850,486)(2850,186)
\put(1200,1461){\makebox(0,0)[lb]{\smash{{\mathmode{\scriptstyle \partial_\alpha}}}}}
\put(3600,1461){\makebox(0,0)[lb]{\smash{{\mathmode{\scriptstyle \partial_\alpha}}}}}
\put(4800,1461){\makebox(0,0)[lb]{\smash{{\mathmode{\scriptstyle \partial_\alpha}}}}}
\put(0,1461){\makebox(0,0)[lb]{\smash{{\mathmode{\scriptstyle \partial_\beta}}}}}
\put(600,1461){\makebox(0,0)[lb]{\smash{{\mathmode{\scriptstyle \partial_\beta}}}}}
\put(2400,1461){\makebox(0,0)[lb]{\smash{{\mathmode{\scriptstyle \partial_\beta}}}}}
\put(3000,1461){\makebox(0,0)[lb]{\smash{{\mathmode{\scriptstyle \partial_\beta}}}}}
\put(4200,1461){\makebox(0,0)[lb]{\smash{{\mathmode{\scriptstyle \partial_\beta}}}}}
\put(1800,1461){\makebox(0,0)[lb]{\smash{{\mathmode{\scriptstyle \partial_\beta}}}}}
\put(675,36){\makebox(0,0)[lb]{\smash{{\mathmode{\scriptstyle \alpha}}}}}
\put(2775,36){\makebox(0,0)[lb]{\smash{{\mathmode{\scriptstyle \alpha}}}}}
\put(4575,36){\makebox(0,0)[lb]{\smash{{\mathmode{\scriptstyle \alpha}}}}}
\end{picture}
}

%% file: draws/BetaWheels.tex
\begingroup\makeatletter\ifx\SetFigFont\undefined
\def\x#1#2#3#4#5#6#7\relax{\def\x{#1#2#3#4#5#6}}%
\expandafter\x\fmtname xxxxxx\relax \def\y{splain}%
\ifx\x\y   
\gdef\SetFigFont#1#2#3{%
  \ifnum #1<17\tiny\else \ifnum #1<20\small\else
  \ifnum #1<24\normalsize\else \ifnum #1<29\large\else
  \ifnum #1<34\Large\else \ifnum #1<41\LARGE\else
     \huge\fi\fi\fi\fi\fi\fi
  \csname #3\endcsname}%
\else
\gdef\SetFigFont#1#2#3{\begingroup
  \count@#1\relax \ifnum 25<\count@\count@25\fi
  \def\x{\endgroup\@setsize\SetFigFont{#2pt}}%
  \expandafter\x
    \csname \romannumeral\the\count@ pt\expandafter\endcsname
    \csname @\romannumeral\the\count@ pt\endcsname
  \csname #3\endcsname}%
\fi
\fi\endgroup
{\renewcommand{\dashlinestretch}{30}
\begin{picture}(5307,1746)(0,-10)
\path(750,1512)(1050,1212)
\path(150,1512)(750,912)
\path(1350,1512)(750,912)(750,312)
\path(1950,1512)(2850,612)
\path(2550,1512)(3150,912)
\path(3150,1512)(3450,1212)
\path(3750,1512)(2850,612)
\path(4350,1512)(4650,1212)
\path(4950,1512)(4650,1212)
\path(4650,912)(4650,612)
\path(2850,612)(2850,312)
\path(4650,1212)(4650,912)
\path(3750,1512)	(3797.276,1519.980)
	(3825.000,1512.000)

\path(3825,1512)	(3875.735,1467.864)
	(3916.646,1415.303)
	(3949.057,1355.509)
	(3974.295,1289.672)
	(3993.685,1218.985)
	(4008.552,1144.639)
	(4014.703,1106.466)
	(4020.220,1067.825)
	(4025.270,1028.864)
	(4030.016,989.734)
	(4034.626,950.582)
	(4039.265,911.557)
	(4044.098,872.809)
	(4049.291,834.487)
	(4061.421,759.713)
	(4076.979,688.428)
	(4097.291,621.823)
	(4123.681,561.089)
	(4157.476,507.418)
	(4200.000,462.000)

\path(4200,462)	(4236.013,434.421)
	(4278.120,408.855)
	(4324.825,386.895)
	(4374.634,370.136)
	(4426.051,360.173)
	(4477.581,358.600)
	(4527.729,367.011)
	(4575.000,387.000)

\path(4575,387)	(4633.579,460.691)
	(4646.180,524.007)
	(4649.078,564.589)
	(4650.000,612.000)

\path(1350,1512)	(1397.959,1523.565)
	(1425.000,1512.000)

\path(1425,1512)	(1463.055,1463.510)
	(1492.266,1411.245)
	(1513.307,1355.791)
	(1526.851,1297.737)
	(1533.574,1237.669)
	(1534.149,1176.174)
	(1529.249,1113.841)
	(1519.549,1051.256)
	(1505.722,989.007)
	(1488.443,927.681)
	(1468.385,867.866)
	(1446.222,810.148)
	(1422.628,755.115)
	(1398.277,703.354)
	(1373.843,655.454)
	(1350.000,612.000)

\path(1350,612)	(1311.374,550.218)
	(1262.864,484.770)
	(1205.555,419.618)
	(1140.533,358.725)
	(1068.881,306.054)
	(1030.908,284.040)
	(991.685,265.568)
	(951.347,251.132)
	(910.029,241.229)
	(867.869,236.353)
	(825.000,237.000)

\path(825,237)	(787.301,255.169)
	(750.000,312.000)

\path(4950,1512)	(4997.006,1522.793)
	(5025.000,1512.000)

\path(5025,1512)	(5075.814,1459.168)
	(5119.916,1402.000)
	(5157.426,1341.060)
	(5188.466,1276.912)
	(5213.154,1210.119)
	(5231.612,1141.245)
	(5243.960,1070.853)
	(5250.319,999.506)
	(5250.808,927.769)
	(5245.548,856.204)
	(5234.659,785.376)
	(5218.262,715.847)
	(5196.478,648.182)
	(5169.425,582.943)
	(5137.226,520.695)
	(5100.000,462.000)

\path(5100,462)	(5067.608,419.363)
	(5032.326,379.655)
	(4994.334,342.767)
	(4953.814,308.591)
	(4910.944,277.019)
	(4865.907,247.941)
	(4818.881,221.249)
	(4770.049,196.835)
	(4719.589,174.590)
	(4667.683,154.405)
	(4614.511,136.173)
	(4560.253,119.784)
	(4505.090,105.130)
	(4449.203,92.102)
	(4392.771,80.593)
	(4335.975,70.493)
	(4278.996,61.693)
	(4222.014,54.086)
	(4165.209,47.563)
	(4108.762,42.015)
	(4052.854,37.333)
	(3997.665,33.410)
	(3943.374,30.137)
	(3890.164,27.404)
	(3838.213,25.104)
	(3787.703,23.128)
	(3738.814,21.368)
	(3691.727,19.715)
	(3646.621,18.060)
	(3603.678,16.295)
	(3563.077,14.311)
	(3525.000,12.000)

\path(3525,12)	(3470.980,10.243)
	(3406.499,11.522)
	(3335.130,15.910)
	(3260.445,23.479)
	(3186.017,34.300)
	(3115.418,48.446)
	(3052.222,65.988)
	(3000.000,87.000)

\path(3000,87)	(2958.855,121.474)
	(2925.000,162.000)

\path(2925,162)	(2895.634,212.790)
	(2875.729,254.585)
	(2850.000,312.000)

\put(0,1587){\makebox(0,0)[lb]{\smash{{\mathmode{\scriptstyle \partial_\beta}}}}}
\put(600,1587){\makebox(0,0)[lb]{\smash{{\mathmode{\scriptstyle \partial_\beta}}}}}
\put(2400,1587){\makebox(0,0)[lb]{\smash{{\mathmode{\scriptstyle \partial_\beta}}}}}
\put(4200,1587){\makebox(0,0)[lb]{\smash{{\mathmode{\scriptstyle \partial_\beta}}}}}
\put(3000,1587){\makebox(0,0)[lb]{\smash{{\mathmode{\scriptstyle \partial_\beta}}}}}
\put(1800,1587){\makebox(0,0)[lb]{\smash{{\mathmode{\scriptstyle \partial_\beta}}}}}
\end{picture}
}

%% file: refs.tex
\ifx\undefined\bysame
        \newcommand{\bysame}{\leavevmode\hbox to3em{\hrulefill}\,}
\fi